\documentclass[a4paper, 10pt]{article}
\usepackage[T2A]{fontenc}
\usepackage[english]{babel}

\usepackage{amsmath, amsthm}
\usepackage{amssymb}
\usepackage{amsfonts}
\usepackage{graphicx}
\usepackage{srcltx, dsfont, color}

\theoremstyle{plain}
\newtheorem{theorem}{Theorem}[section]

\numberwithin{equation}{section}

\def\Tht{\Theta}
\def\tht{\theta}
\def\Om{\Omega}
\def\om{\omega}
\def\e{\varepsilon}
\def\g{\gamma}
\def\G{\Gamma}
\def\l{\lambda}
\def\p{\partial}
\def\D{\Delta}
\def\a{\alpha}

\def\Ups{\Upsilon}

\def\d{\delta}

\def\z{\zeta}
\def\vs{\varsigma}

\def\vk{\varkappa}

\def\Ho{\mathring{W}}

\def\hf{\mathfrak{h}}

\def\la{\langle}
\def\ra{\rangle}

\def\hg{\mathfrak{g}}

\def\iu{\mathrm{i}}

\def\Op{\mathcal{H}}

\def\fV{\mathfrak{V}}
\def\fM{\mathfrak{M}}

\def\cB{\mathcal{B}}
\def\cL{\mathcal{L}}
\def\cF{\mathcal{F}}
\def\cI{\mathcal{I}}

\def\cX{\mathcal{X}}

\def\Ho{\mathring{W}}

\DeclareMathOperator{\RE}{Re}
\DeclareMathOperator{\IM}{Im}

\DeclareMathOperator{\mes}{mes}

\DeclareMathOperator{\dist}{dist}

\DeclareMathOperator*{\esssup}{ess\,sup}

\allowdisplaybreaks

\begin{document}

\title{Homogenization for operators with \\ arbitrary perturbations in coefficients}

\author{D.I. Borisov$^{1,2,3}$
}

\date{\empty}

\maketitle

{\small
    \begin{quote}
1) Institute of Mathematics, Ufa Federal Research Center, Russian Academy of Sciences,  Chernyshevsky str. 112, Ufa, Russia, 450008
\\
2) Bashkir State  University, Zaki Validi str. 32, Ufa, Russia, 450076
\\
3) University of Hradec Kr\'alov\'e
62, Rokitansk\'eho, Hradec Kr\'alov\'e 50003, Czech Republic
\\
Emails: borisovdi@yandex.ru
\end{quote}

{\small
\begin{quote}
\noindent{\bf Abstract.}  We consider a general second order matrix operator in a multi-dimensional domain subject to a classical boundary condition. This operator is perturbed by a first order differential operator, the coefficients of which depend arbitrarily on a small multi-dimensional parameter. We study the existence of a limiting (homogenized) operator in the sense of the norm resolvent convergence for such perturbed operator. The first part of our  main results states that the norm resolvent convergence is equivalent to the convergence of the coefficients in the perturbing operator in certain space of multipliers. If this is the case, the resolvent of the perturbed operator possesses a complete asymptotic expansion, which converges uniformly to the resolvent. The second part of our results says that the convergence in the mentioned spaces of multipliers is equivalent to the convergence of certain local mean values over small pieces of the considered domains. These results are supported by series of examples. %including classical %cases like locally periodic fast oscillating %coefficients as well as by %some new examples.
 We also provide a series of ways of generating new non-periodically oscillating perturbations, which finally leads to a very wide class of perturbations, for which our results are applicable.

\medskip


\medskip

 	
\end{quote}
}

\section{Introduction}

Nowadays the homogenization theory is a big part of the modern mathematics  aimed on understanding the behavior of solutions to differential equations with various perturbations. The main feature  of the perturbations is their localizations on small closely spaced sets; a classical example is a fast oscillating coefficient of form $V(x/\e)$, where $V$ is a periodic function and $\e$ is a small parameter. Given a perturbed problem, the usual main question in the homogenization theory is to describe the form of a homogenized problem, a solution of which serves as a limit for a solution of the perturbed problem. Another question is to estimate the convergence rates. Not trying to mention all works in homogenization theory, we just cite a few classical books, \cite{BaPa}, \cite{BePa}, \cite{MK1}, \cite{OIS}, \cite{ZKO}.

During the last 20 years there appeared a new direction in the homogenization theory aimed on establishing so-called operator estimates. The matter is to prove that the solutions of the perturbed problems converge to ones of the limiting problems uniformly in the norm of the given right hand sides; for  linear operators this is equivalent to proving the norm resolvent convergence. A lot of various operator estimates were established for operators with fast oscillating coefficients, see, for instance, \cite{BS1}, \cite{DS}, \cite{Gr1}, \cite{Gr2}, \cite{Kenig}, \cite{ZhP2}, \cite{Sen2}, \cite{Sen3}, \cite{Sen4}, \cite{BS2}, \cite{ZhJi}, \cite{Zhi}, \cite{ZhP1}, the references therein and other works by these authors. For many cases, the convergence rates were estimated and many estimates were order sharp.
We also mention papers \cite{DNT}, \cite{KKN}, in which a hyperbolic  equation  with fast and non-periodically oscillating coefficients was studied; the non-periodicity was controlled by apriori bounds of the perturbed coefficients with respect to a small parameter. A close paper  is \cite{DN}, in which the perturbed coefficients depended on two fast variables controlled by two small parameters. The dependence on one fast variable was periodic, while with respect to the other fast variable the coefficients rapidly decay as it increased. The main results in these papers described the asymptotic behavior of the solutions to the considered problems, in general, non-uniformly in the norms of the right hand side.

The cited results for problems with fast oscillating coefficients stimulated similar studies for problems with geometric perturbations like perforations, frequent alternation of boundary conditions, oscillating boundary. Apart of periodic and locally periodic cases,   the cases of general non-periodic perturbations were treated and the operator estimates were obtained, see \cite{Post}, \cite{PRSE16}, \cite{MSB21}, \cite{ZAMP}, \cite{JST}, \cite{Khra2}.  These results on non-periodic perturbations motivate studying the possibility of obtaining operator estimates for other non-periodic perturbations from the homogenization theory. A first classical case is the perturbation described by non-periodically fast oscillating coefficients, which is a natural continuation of the aforementioned works \cite{BS1}, \cite{Gr1}, \cite{Gr2}, \cite{Kenig}, \cite{ZhP2}, \cite{Sen2}, \cite{BS2}, \cite{ZhJi}, \cite{Zhi}, \cite{ZhP1}. The problem can be even formulated in a more general setting: given coefficients arbitrarily depending on one or several small parameters, what are the conditions on this dependence which ensure the norm resolvent convergence to some homogenized operator? And how to determine the homogenized operator? These two questions are addressed in the present paper.

We consider a general linear second order matrix differential operator in an arbitrary multi-dimensional domain with a sufficiently regular boundary. On the boundary we impose either the Dirichlet or Robin condition. This operator is perturbed by a first order differential operator with the matrix coefficients depending on a multi-dimensional parameter. We show that the perturbation gives rise to a norm resolvent convergence if and only if the matrix coefficients in the perturbing operator converges in certain spaces of multipliers. Then we study the convergence in these spaces and prove that it is equivalent to the uniform convergence of mean values and averaged $L_2$-norms over certain small subsets in the considered domain, which is a rather simple explicit criterion. We show how this criterion can be used for finding the limits in the mentioned spaces of multipliers. An additional important result is that once the criterion for the norm resolvent convergence holds, we also obtain a complete asymptotic expansion for the perturbed resolvent. This expansion is shown to converge uniformly to the perturbed resolvent. The obtained general results are supported by a series of examples including classical and new ones.

The paper is organized as follows. In the next section we formulate the problem and then we present  and discuss our main results. In the third section we provide a series of examples and show how to produce new examples from given ones. The fourth section is devoted to proving the equivalence between the norm resolvent convergence and the convergence in the multipliers spaces.  The convergence in these spaces is studied in the fifth section.

\section{Problem and main results}\label{sec2}

\subsection{Problem}

Let $x=(x_1,\ldots,x_d)$ be Cartesian coordinates in $\mathds{R}^d$ and $\Om$ be  an arbitrary bounded or unbounded domain in $\mathds{R}^d$.
If the boundary of $\Om$ is  non-empty, we assume that  $\p\Om\in C^2$. By $\tau$ we denote a distance measured along the inward normal vector  to $\p\Om$. We suppose that the variable $\tau$ is well-defined at least on
\begin{equation}\label{2.0}
\Pi(\tau_0):=\{x\in\Om:\ \dist(x,\p\Om)\leqslant \tau_0\}
\end{equation}
for some fixed $\tau_0>0$. Denoting by $s$ a local variable on $\p\Om$, we also assume that on $\Pi(\tau_0)$, the Jacobians associated with the passing from variables $x$ to $(\tau,s)$ are bounded and separated from zero uniformly and all first derivatives of $x$ in $(\tau,s)$ are also uniformly bounded.

By $L_2(\Om;\mathds{C}^n)$ we denote the space of vector functions defined on $\Om$ with values in $\mathds{C}^n$ such that each component of these vector functions belongs to $L_2(\Om)$. This space is equipped with the natural norm
\begin{equation*}
\|u\|_{L_2(\Om;\mathds{C}^n)}^2:=\sum\limits_{j=1}^{n} \|u_j\|_{L_2(\Om)}^2,\qquad u=(u_1,\ldots,u_n).
\end{equation*}

Let $\mathds{M}_n$ be the space of all square $n\times n$ matrices with complex-valued entries. The symbol $L_\infty(\Om;\mathds{M}_n)$ stands for the space of $n\times n$ matrix functions defined on $\Om$ such that each entry of these matrix functions belongs to  $L_\infty(\Om)$. This space is equipped with the norm
\begin{equation*}
\|A\|_{L_\infty(\Om;\mathds{M}_n)}:=\||A|\|_{L_\infty(\Om)},\qquad |A|=\sum\limits_{i,j=1}^{n} |a_{ij}|,\qquad A=(a_{ij}).
\end{equation*}
Apart from $L_\infty(\Om;\mathds{M}_n)$ and $L_2(\Om;\mathds{C}^n)$, we shall also use similar notations for various functional spaces of vector and matrix functions.

We introduce  a second order matrix differential operator $\Op$ in $L_2(\Om;\mathds{C}^n)$ with the differential expression
\begin{equation*}%\%label{2.2}
-\sum\limits_{i,j=1}^{d} \frac{\p\ }{\p x_i} A_{ij} \frac{\p\ }{\p x_j} + \sum\limits_{j=1}^{d} A_j^+ \frac{\p\ }{\p x_j} + \sum\limits_{j=1}^{d} \frac{\p\ }{\p x_j}A_j^- + A_0,
\end{equation*}
where $A_{ij}=A_{ij}(x)$, $A_j^\pm=A_j^\pm(x)$, $A_0=A_0(x)$ are some complex-valued $n\times n$ matrix functions, $n\geqslant 1$, obeying the conditions
\begin{equation}\label{2.3}
A_{ij},\ A_j^\pm,\ A_0\in L_\infty(\Om;\mathds{M}_n), \qquad
\RE \sum\limits_{i,j=1}^{d} (A_{ij}(x)z_i,z_j)_{\mathds{C}^n} \geqslant c_1 \sum\limits_{j=1}^{d} |z_i|^2,
\qquad z_i\in\mathds{C}^n,\quad x\in\Om,
\end{equation}
where $c_1>0$ is some fixed constant independent of $x$ and $z_i$.
The boundary condition for the operator $\Op$ is chosen to be a classical one and we write it as
\begin{equation}\label{2.4}
\cB u=0,\quad\text{where}\quad \cB u=u\quad\text{or}\quad \cB u= \frac{\p u}{\p\boldsymbol{\nu}} + K u,\qquad
\frac{\p u}{\p\boldsymbol{\nu}}:=\sum\limits_{i,j=1}^{d} \nu_i A_{ij}\frac{\p u}{\p x_j}-\sum\limits_{j=1}^{d} \nu_j A_j^- u.
\end{equation}
Here $\nu=(\nu_1,\ldots,\nu_d)$ is the unit outward normal to $\p\Om$ and $K=K(x)$ is some $n\times n$ matrix function defined on $\p\Om$ and belonging to $L_\infty(\p\Om;\mathds{M}_n)$. We shall define the operator $\Op$ rigorously in the next subsection.

Let $\e=(\e_1,\ldots,\e_m)$, $m\geqslant 1$, be a  multi-dimensional small parameter, that is, all $\e_j$ are supposed to be small. By $V^\e=V^\e(x)$, $Q_j^\e=Q_j^\e(x)$, $P_j^\e=P_j^\e(x)$, $j=1,\ldots,d$, we denote arbitrary families of complex-valued $n\times n$ matrix functions belonging to $L_\infty(\Om;\mathds{M}_n)$ and the only assumption we make is that all these families are bounded uniformly in $\e$.

The main object of our study is an operator $\Op^\e$ which can be formally written as
\begin{equation}\label{2.38}
\Op^\e:=\Op + \sum\limits_{j=1}^{d} Q_j^\e \frac{\p\ }{\p x_j} + \sum\limits_{j=1}^{d}  \frac{\p\ }{\p x_j} P_j^\e+ V^\e
\end{equation}
with the boundary condition $\cB^\e u=0$, where
\begin{equation}\label{2.33}
 \cB^\e u =u \quad \text{if}\quad \cB u=u,
\qquad \cB^\e u =\cB u - \sum\limits_{j=1}^{d} P_j^\e \nu_j u\quad \text{if}\quad \cB u= \frac{\p u}{\p\boldsymbol{\nu}} + K u.
\end{equation}
We shall introduce this operator rigorously in the next subsection.

The aim of this paper is to find conditions for the families $Q_j^\e$, $P_j^\e$ and $V^\e$ ensuring the norm resolvent convergence of the operator $\Op^\e$ to some limiting operator of form
\begin{equation*}%\%label{2.28}
\Op^0:=\Op + \sum\limits_{j=1}^{d} Q_j^0 \frac{\p\ }{\p x_j} +\sum\limits_{j=1}^{d}   \frac{\p\ }{\p x_j} P_j^0 + V^0
\end{equation*}
with the boundary condition  $\cB^0 u=0$, where
\begin{equation*}%\%label{2.34}
 \cB^0 u =u \quad \text{if}\quad \cB u=u,
\qquad \cB^0 u =\cB u - \sum\limits_{j=1}^{d} P_j^0 \nu_j u \quad \text{if}\quad \cB u=\frac{\p u}{\p\boldsymbol{\nu}}+ K u.
\end{equation*}
This operator will be also rigorously defined in the next subsection.

\subsection{Operators and multipliers}

In order to formulate our main results, we first define rigorously the operators $\Op$, $\Op^\e$ and $\Op^0$ and introduce some auxiliary spaces.

By $\Ho_2^1(\Om;\mathds{C}^n)$ we denote a subspace of the Sobolev space $W_2^1(\Om;\mathds{C}^n)$ consisting of vector functions with the zero trace on $\p\Om$. In $L_2(\Om;\mathds{C}^n)$ we introduce two sesqulinear forms
\begin{align*}%\l%abel{2.5}
\hf_D(u,v):=&\sum\limits_{i,j=1}^{d} \left(A_{ij}\frac{\p u}{\p x_j}, \frac{\p v}{\p x_i}\right)_{L_2(\Om;\mathds{C}^n)} + \sum\limits_{j=1}^{d}
\left(A_j^+\frac{\p u}{\p x_j}, v \right)_{L_2(\Om;\mathds{C}^n)}
 \\
 &- \sum\limits_{j=1}^{d}
\left(A_j^- u, \frac{\p v}{\p x_j} \right)_{L_2(\Om;\mathds{C}^n)} + (A_0 u, v)_{L_2(\Om;\mathds{C}^n)}
\end{align*}
on the domain $\Ho_2^1(\Om;\mathds{C}^n)$ and
\begin{align*}%\l%abel{2.6}
\hf_R(u,v):=&\sum\limits_{i,j=1}^{d} \left(A_{ij}\frac{\p u}{\p x_j}, \frac{\p v}{\p x_i}\right)_{L_2(\Om;\mathds{C}^n)} + \sum\limits_{j=1}^{d}
\left(A_j^+\frac{\p u}{\p x_j}, v \right)_{L_2(\Om;\mathds{C}^n)}
 \\
 &- \sum\limits_{j=1}^{d}
\left(A_j^- u, \frac{\p v}{\p x_j} \right)_{L_2(\Om;\mathds{C}^n)}  + (A_0 u, v)_{L_2(\Om;\mathds{C}^n)}+(Ku,v)_{L_2(\p\Om;\mathds{C}^n)}
\end{align*}
on the domain $W_2^1(\Om;\mathds{C}^n)$. Owing to conditions (\ref{2.3}), the introduced forms are obviously sectorial and closed. Then by the first representation theorem \cite[Ch. V\!I, Sec. 2.1, Thm. 2.1]{Kato}, for each form there exists a unique associated $m$-sectorial operator in $L_2(\Om;\mathds{C}^n)$. These are the needed operators. If we deal with the Dirichlet condition in (\ref{2.4}), then by $\Op$ we denote the operator associated with the form $\hf_D$. If we deal with the Robin condition in (\ref{2.4}), then the symbol $\Op$ stands for the operator associated with the form $\hf_R$. The form associated with the operator $\Op$ is denoted by $\hf$ and this is the form  $\hf_D$ or $\hf_R$. In what follows, the domain of the chosen form $\hf$ is denoted by $\fV$. This is a Hilbert space with the norm $\|\cdot\|_{\fV}:=\|\cdot\|_{W_2^1(\Om;\mathds{C}^n)}$.

Then we extend the operator $\Op$ to entire $\fV$ as follows. Following \cite[Ch. I, Sect. 2.2]{Kato}, we introduce the adjoint space for $\fV$, which is denoted by $\fV^\ast$. This is the linear space of  continuous antilinear functionals on $\fV$, that is, of functionals $\cF$ on $\fV$ obeying the property
$$
\la \cF, \a_1 v_1 +\a_2 v_2\ra=\overline{\a_1}\la \cF, v_1\ra+ \overline{\a_2}\la \cF, v_2\ra
$$
and being continuous; here $\la \cF, v\ra$ denotes the action of the functional $\cF$ on $v\in\fV$.
The space $\fV^\ast$ is equipped with a natural norm and is a Banach space. Throughout this work we deal with the realization of the adjoint space $\fV^\ast$ with respect to the pairing in $L_2(\Om;\mathds{C}^n)$, in particular,  functions from $L_2(\Om;\mathds{C}^n)$ naturally generate functionals in $\fV^\ast$ and in this sense $L_2(\Om;\mathds{C}^n)\subset\fV^\ast$.

To each $u\in \fV$, we associate a bounded antilinear functional from $\fV^\ast$ acting by the rule $v\mapsto \hf(u,v)$; the boundedness is ensured by the obvious estimate $|\hf(u,v)|\leqslant C \|u\|_{\fV}\|v\|_{\fV}$. This mapping is the extension of the operator $\Op$ on $\fV$;  indeed,   if $u$ belongs to the domain of $\Op$, then for each $v\in \fV$ we have
$\hf(u,v)=(\Op u,v)_{L_2(\Om;\mathds{C}^n)}$. The right hand side of this identity is another representation for the aforementioned functional associated with $u$ and in this sense we mean the above made extension of the operator $\Op$.

The operators $\Op^\e$ and $\Op^0$ are defined in the same way. The corresponding forms are
\begin{align}\label{2.35}
&\hf^\e(u,v):=\hf(u,v) + \sum\limits_{j=1}^{d}
\left(Q_j^\e \frac{\p u}{\p x_j},v\right)_{L_2(\Om;\mathds{C}^n)} - \sum\limits_{j=1}^{d}
\left(P_j^\e u,\frac{\p v}{\p x_j}\right)_{L_2(\Om;\mathds{C}^n)} + (V^\e u,v)_{L_2(\Om;\mathds{C}^n)},
\\
\label{2.36}
&\hf^0(u,v):=\hf(u,v) + \sum\limits_{j=1}^{d}
\left(Q_j^0 \frac{\p u}{\p x_j},v\right)_{L_2(\Om;\mathds{C}^n)} - \sum\limits_{j=1}^{d}
\left(P_j^0 u, \frac{\p v}{\p x_j}\right)_{L_2(\Om;\mathds{C}^n)} + (V^0 u,v)_{L_2(\Om;\mathds{C}^n)}
\end{align}
on the domain $\fV$. We can consider these forms in $L_2(\Om;\mathds{C}^n)$; it will be shown in Section~\ref{sec:th1} that they are closed and sectorial and by the first representation theorem \cite[Ch. V\!I, Sec. 2.1, Thm. 2.1]{Kato} there exist  unique $m$-sectorial operators $\Op^\e$ and $\Op^0$ associated with these forms. These are the operators $\Op^\e$ and $\Op^\e$ considered as  unbounded ones in $L_2(\Om;\mathds{C}^n)$. The domain of these operators are subsets of $\fV$. Then we extend these operators on the entire space $\fV$ as above keeping the same notations.

By $\fM_{1,-1}$ we denote the space of multipliers from $\fV$ into   $\fV^\ast$; this space of multipliers consists of $n\times n$ matrix distributions $V$ defined on $\Om$ such that $V u\in \fV^\ast$ for each $u\in\fV$ and the multiplication by $V$ is a bounded operator from $\fV$ into $\fV^*$. The norm in $\fM_{1,-1}$ is given by the formula
\begin{equation}\label{2.7}
\|V\|_{\fM_{1,-1}}:=\sup\limits_{\substack{u\in \fV \\ u\ne0}}
\frac{\|Vu\|_{\fV^\ast}}{\|u\|_{\fV}}=\sup\limits_{\substack{u,v\in \fV \\ u,v\ne0}}
\frac{|\la Vu,v\ra|}{\|u\|_{\fV}\|v\|_{\fV}}.
\end{equation}
We  introduce a similar space $\fM_{1,0}$ of  multipliers from $\fV$ into $L_2(\Om;\mathds{C}^n)$. The norm in this space is defined as
\begin{equation}\label{2.11}
\|V\|_{\fM_{1,0}}:= \sup\limits_{\substack{u\in \fV\\ u\ne0}}
\frac{\|Vu\|_{L_2(\Om;\mathds{C}^n)}}{\|u\|_{\fV}}.
\end{equation}
And one more space we shall need is denoted by $\fM_{2,0}$ and it consists of multipliers from $\fV\cap W_2^2(\Om;\mathds{C}^n)$ into  $L_2(\Om;\mathds{C}^n)$ with the norm
\begin{equation*}
\|V\|_{\fM_{2,0}}:= \sup\limits_{\substack{u\in \fV\cap W_2^2(\Om;\mathds{C}^n)\\ u\ne0}}
\frac{\|Vu\|_{L_2(\Om;\mathds{C}^n)}}{\|u\|_{W_2^2(\Om;\mathds{C}^n)}}.
\end{equation*}
It is obvious that
\begin{equation}\label{2.39}
\fM_{1,0}\subset \fM_{1,-1},\qquad \fM_{1,0}\subset \fM_{2,0}, \qquad
\|V\|_{\fM_{1,-1}}\leqslant \|V\|_{\fM_{1,0}}, \qquad  \|V\|_{\fM_{2,0}}\leqslant \|V\|_{\fM_{1,0}}.
\end{equation}

By $\|\,\cdot\,\|_{X\to Y}$ we denote the norm of a bounded operator acting from a Banach space $X$ into a Banach space $Y$.
We consider the space $(L_\infty(\Om;\mathds{C}^n))^{2d+1}$ and with each its element $X:=(Q_1,\ldots,Q_d,P_1,\ldots,P_d,V)$ we associate
 a linear bounded operator $\cX$ acting from $\fV$ into $\fV^\ast$ by the following rule:
\begin{equation}\label{2.37}
\la\cX u,v\ra:= \sum\limits_{j=1}^{d}
\left(Q_j\frac{\p u}{\p x_j},v\right)_{L_2(\Om;\mathds{C}^n)} - \sum\limits_{j=1}^{d}
\left(P_j u,\frac{\p v}{\p x_j}\right)_{L_2(\Om;\mathds{C}^n)} + (V u,v)_{L_2(\Om;\mathds{C}^n)}\quad\text{for}\quad v\in\fV.
\end{equation}
%where $\la\cX u,v\ra$ denotes the action of the functional $\cX u$ on an %element $v$.
We observe that if $\cX=0$, this does not imply $P_j=Q_j=V=0$. Indeed, by a simple integration by parts it is easy to see that if the matrix functions $P_j$, $Q_j$ are continuously differentiable, compactly supported and satisfy the identities $Q_j=-P_j$, $V=\sum\limits_{j=1}^{d} \frac{\p Q_j}{\p x_j}$, then $\cX=0$. This means that each operator $\cX$ can be simultaneously represented by different sets of coefficients $Q_j$, $P_j$, $V$. The same situation also holds for the differential expression (\ref{2.38}): we can rewrite the second sum as
\begin{equation*}
\sum\limits_{j=1}^{n}\frac{\p\ }{\p x_j} P_j^\e=\sum\limits_{j=1}^{n} P_j^\e \frac{\p\ }{\p x_j} + \sum\limits_{j=1}^{n} \frac{\p P_j^\e}{\p x_j}
\end{equation*}
and replace then $Q_j^\e$ and $V^\e$ by $Q_j^\e+P_j^\e$ and $V^\e+\sum\limits_{j=1}^{n} \frac{\p P_j^\e}{\p x_j}$ not changing the action of the differential expression.

We denote by $\fM$ the space of all bounded operators from $\fV$ into $\fV^\ast$. If an operator $\cX$ is defined by (\ref{2.37}), then its norm is given by the identity
\begin{equation}\label{2.26}
\|\cX\|_{\fM}=\sup\limits_{\substack{u,v\in
\fV\\ u\ne0,\, v\ne0}}
\frac{\left| \sum\limits_{j=1}^{d}
\left(Q_j\frac{\p u}{\p x_j},v\right)_{L_2(\Om;\mathds{C}^n)} - \sum\limits_{j=1}^{d}
\left(P_j u,\frac{\p v}{\p x_j}\right)_{L_2(\Om;\mathds{C}^n)} + (V u,v)_{L_2(\Om;\mathds{C}^n)}
\right| }{\|u\|_{\fV} \|v\|_{\fV}}.
\end{equation}
Given an element $X\in (L_\infty(\Om;\mathds{C}^n))^{2d+1}$ and the associated operator $\cX$ defined in (\ref{2.37}), it is easy to see that
\begin{equation}\label{2.29}
\|\cX\|_{\fM}\leqslant \sum\limits_{j=1}^{d} \big( \|Q_j^*\|_{\fM_{1,0}}+  \|P_j\|_{\fM_{1,0}}\big)
 + \|V\|_{\fM_{1,-1}}.
\end{equation}
We shall show in Section~\ref{sec53} that
\begin{equation}\label{2.28}
\|Q^*\|_{\fM_{1,0}}\leqslant \sqrt{n} \|Q\|_{\fM_{1,0}}
\end{equation}
and hence, estimate (\ref{2.29}) implies
\begin{equation}\label{2.27}
\|\cX\|_{\fM}\leqslant \sum\limits_{j=1}^{d} \big( \sqrt{n}\|Q_j\|_{\fM_{1,0}}+  \|P_j\|_{\fM_{1,0}}\big)
 + \|V\|_{\fM_{1,-1}}.
\end{equation}

Let $\cX^\e$ and $\cX^0$ be the operators defined by (\ref{2.37}) via the families $P_j^\e$, $Q_j^\e$, $V^\e$ and $P_j^0$, $Q_j^0$, $V^0$. Then  $\Op^\e=\Op+\cX^\e$, $\Op^0=\Op+\cX^0$, where all operators are treated as acting from $\fV$ into $\fV^*$.

\subsection{Main results}

Now we are in position to formulate our main results. The first of them describes a way of approximating  the resolvent of $\Op^\e$.

\begin{theorem}\label{th1}
Assume that the  family  $X^\e:=(Q_1^\e,\ldots,Q_d^\e,P_1^\e,\ldots,P_d^\e,V^\e)$ is bounded uniformly in $(L_\infty(\Om;\mathds{C}^n))^{2d+1}$ and
the family of the corresponding operators $\cX^\e$ given by (\ref{2.37})
converges to some operator $\cX^0$ in  the sense of the norm $\|\,\cdot\,\|_{\fM}$. Then there exists an element
$X^0:=(Q_1^0,\ldots,Q_d^0,P_1^0,\ldots,P_d^0,V^0)\in (L_\infty(\Om;\mathds{C}^n))^{2d+1}$ generating  the operator $\cX^0$   by (\ref{2.37}). In terms of the elements $X^\e$ and $X^0$ we define the operators $\Op^\e$ and $\Op^0$ via the forms in (\ref{2.35}), (\ref{2.36}). Then there exists a real number $\l_0$ independent of $\e$ such that each $\l\in\mathds{C}$ obeying $\RE\l<\l_0$ belongs to the resolvent sets of $\Op^0$ and $\Op^\e$ for all $\e$. For each $\l\in\mathds{C}$ with $\RE\l<\l_0$ the resolvent $(\Op^\e-\l)^{-1}$ converges to $(\Op^0-\l)^{-1}$
in the $\|\,\cdot\,\|_{\fV^\ast\to \fV}$-norm and it is represented by the convergent series
\begin{equation}\label{2.8}
(\Op^\e-\l)^{-1}=\sum\limits_{j=0}^{\infty}(\Op^0-\l)^{-1} \big(-\cL^\e(\Op^0-\l)^{-1}\big)^j, \qquad  \cL^\e:=\cX^\e-\cX^0.
\end{equation}
Series (\ref{2.8}) converges uniformly in $\e$ in the $\|\,\cdot\,\|_{\fV^\ast\to \fV}$-norm and the estimates hold:
\begin{equation}\label{2.9}
%\begin{aligned}
\bigg\|(\Op^\e-\l)^{-1} -(\Op^0-\l)^{-1}%&
 \sum\limits_{j=0}^{N} (-1)^j \big(\cL^\e(\Op^0-\l)^{-1}\big)^j
\bigg\|_{\fV^\ast
\to \fV}
%\\
%&
\leqslant  c_2^{N+2}(\l) \|\cL^\e\|_{\fM}^{N+1}
%\|(\Op^0-\l)^{-1}\|_{\fV^\ast
%\to \fV}
%\end{aligned}
\end{equation}
for each $N\in\mathds{Z}_+$, where $c_2=c_2(\l)$ is some fixed constant independent of $\e$ and $N$.
\end{theorem}

Our second theorem says that if the operator $\Op^\e$ converges to $\Op^0$ in the norm resolvent sense, this  implies a  convergence
of the family  $\cX^\e$ to $\cX^0$ in the norm $\|\,\cdot\,\|_{\fM}$.

\begin{theorem}\label{th5}
Suppose that
\begin{equation}\label{2.22}
\vk(\e):=\sup\limits_{\substack{
f\in \fV^\ast
\\
f\ne0
}}
\frac{\|(\Op^\e-\l)^{-1}f -(\Op^0-\l)^{-1}f\|_{\fV}} {\|f\|_{\fV^*}}
\to0,\qquad \e\to0.
\end{equation}
Then
\begin{equation}\label{2.23}
\|\cX^\e-\cX^0\|_{\fM} \leqslant C \vk(\e)\to 0,\qquad
 \e\to0,
\end{equation}
where
$C$ is some constant independent of $\e$. In particular, if $Q_j^\e\equiv 0$, $P_j^\e\equiv 0$, $j=1,\ldots,d$, then the family $V^\e$ converges to $V^0$ in $\fM_{1,-1}$ and
\begin{equation}\label{2.24}
\|V^\e-V^0\|_{\fM_{1,-1}}\leqslant C \vk(\e)\to 0,\qquad
 \e\to0,
\end{equation}
where $C$ is some constant independent of $\e$.
\end{theorem}

Our third theorem deals with an important case when the domain of the operator $\Op^0$ treated as an unbounded $m$-sectorial operator in $L_2(\Om;\mathds{C}^n)$ is a subset in $W_2^2(\Om;\mathds{C}^n)$.

\begin{theorem}\label{th6}
Assume that the families $V^\e$,  $Q_j^\e$ converge to  $V^0$,  $Q_j^0$
in $\fM_{1,-1}$, the family $P_j^\e$ converges to  $P_j^0$ in $\fM_{2,0}$, the domain of the operator $\Op^0$ considered in $L_2(\Om;\mathds{C}^n)$ is a subset of $W_2^2(\Om;\mathds{C}^n)$ and the resolvent of $\Op^0$ is a bounded operator from $L_2(\Om;\mathds{C}^n)$ into $W_2^2(\Om;\mathds{C}^n)$. Then the operator $\Op^\e$, considered in $L_2(\Om;\mathds{C}^n)$, converges to $\Op^0$ in the norm resolvent sense and for all $\l$ such that $\RE \l<\l_0$ the estimate holds:
\begin{equation}\label{2.10}
\big\|(\Op^\e-\l)^{-1}-(\Op^0-\l)^{-1}\big\|_{L_2(\Om;\mathds{C}^n)\to W_2^1(\Om;\mathds{C}^n)}
\leqslant c_3(\l)\bigg(\|V^\e-V^0\|_{\fM_{1,-1}} + \sum\limits_{j=1}^{d} \|Q_j^\e-Q_j^0\|_{\fM_{1,-1}} + \sum\limits_{j=1}^{d} \|P_j^\e-P_j^0\|_{\fM_{2,0}}\bigg),
\end{equation}
where $\l_0$ is from Theorem~\ref{th1} and $c_3(\l)$ is some constant independent of $\e$.
\end{theorem}

 The second part of our main results is devoted to simple necessary and sufficient conditions ensuring the convergence in the spaces $\fM_{1,-1}$ and $\fM_{1,0}$.

By $\G$ we denote some fixed lattice in $\mathds{R}^d$ with a periodicity cell $\square$. Given $\eta>0$, we let
\begin{equation}\label{2.18}
\square_\g^\eta:=\eta\square+\eta\g=\big\{x\in\mathds{R}^d:\ \eta^{-1} x-\g \in\square\big\}, \qquad \G_\eta:=\big\{\g\in\G:\ \square_\g^\eta\subset\Om\big\}.
\end{equation}
The next theorem gives a criterion of the convergence in the space $\fM_{1,-1}$.

\begin{theorem}\label{th2}
Let $V^\e\in L_\infty(\Om;\mathds{M}_n)$ be a family of matrix functions bounded   uniformly in $\e$ and let there exist a matrix function $V^0\in L_\infty(\Om;\mathds{M}_n)$ and a scalar positive function $\eta=\eta(\e)$ such that
\begin{equation}\label{2.13}
\bigg|\int\limits_{\square_\g^\eta} \eta^{-d}\big(V^\e(x)-V^0(x)\big)\,dx\bigg|\leqslant \rho_1(\e) \quad\text{for each}\quad
\g\in\G_\eta,\qquad \eta(\e)\to 0,\quad \e\to0,
\end{equation}
where $\rho_1=\rho_1(\e)$ is some function independent of $\g$ such that $\rho_1(\e)\to0$ as $\e\to0$. Then $V^\e$ converges to $V^0$ in $\fM_{1,-1}$  as $\e\to0$ and
\begin{equation}\label{2.14}
\|V^\e-V^0\|_{\fM_{1,-1}}\leqslant C\big(\rho_1(\e)+\eta(\e)\big),
\end{equation}
where $C$ is some constant independent of $\e$, $\eta$ and $\rho_1$. And vice versa, if a uniformly bounded family $V^\e\in L_\infty(\Om;\mathds{M}_n)$ converges to $V^0\in L_\infty(\Om;\mathds{M}_n)$ in $\fM_{1,-1}$, then
condition (\ref{2.13}) is satisfied with
\begin{equation*}
\eta(\e)=\|V^\e-V^0\|_{\fM_{1,-1}}^{\frac{1}{2}},\qquad \rho_1(\e)= C\|V^\e-V^0\|_{\fM_{1,-1}}^{\frac{1}{4}},
\end{equation*}
where $C$ is some constant independent of $\e$.
\end{theorem}

The following theorem describes a very natural way of calculating the limit in $\fM_{1,-1}$.

\begin{theorem}\label{th3}
Let $V^\e\in L_\infty(\Om;\mathds{M}_n)$ be a family of matrix functions bounded   uniformly in $\e$ such that  for almost each $x\in\Om$ there exists a finite limit
\begin{equation}\label{2.16}
V^0(x):=\lim\limits_{\e\to0}
\frac{1}{\mu^d(\e)\mes\om} \int\limits_{x+\mu(\e)\om} V^\e(y)\,dy,
\end{equation}
where $\om$ is some fixed bounded domain, $\mu=\mu(\e)$ is some positive function such that $\mu(\e)\to0$ as $\e\to0$.
Assume that the limit in (\ref{2.16}) is uniform in $x\in \Om$,
namely,
\begin{equation}\label{2.17}
 \bigg|V^0(x)-\frac{1}{\mu^d(\e)\mes\om} \int\limits_{x+\mu(\e)\om} V^\e(y)\,dy\bigg|\leqslant \rho_2(\e),\qquad \rho_2(\e)\to0,\quad \e\to0,
\end{equation}
for all $x\in\Om$ such that $x+\mu(\e)\om\subset\Om$,
where the function $\rho_2=\rho_2(\e)$
is independent of the spatial variables. Then $V^0\in L_\infty(\Om;\mathds{M}_n)$, the family $V^\e$ converges to $V^0$ in $\fM_{1,-1}$ and the estimate holds:
\begin{equation}\label{2.19}
\|V^\e-V^0\|_{\fM_{1,-1}}\leqslant C\big(\rho_2(\e)
+\mu^{\frac{1}{2}}(\e)\big),
\end{equation}
where $C$ is some constant independent of $\e$, $\mu$, $\rho_2$.
\end{theorem}

Our final theorem provides a criterion of the convergence in the space $\fM_{1,0}$.

\begin{theorem}\label{th4}
Let $Q^\e\in L_\infty(\Om;\mathds{M}_n)$ be a family of matrix functions bounded   uniformly in $\e$ and let there exists a matrix function $Q^0\in L_\infty(\Om;\mathds{M}_n)$ and a scalar positive function $\eta=\eta(\e)$ such that
\begin{equation}\label{2.20}
 \int\limits_{\square_\g^\eta} \eta^{-d}\big|Q^\e(x)-Q^0(x)\big|^2\,dx \leqslant \rho_3(\e) \quad\text{for each}\quad
\g\in\G_\eta,\qquad \eta(\e)\to 0,\quad \e\to0,
\end{equation}
where $\rho_3=\rho_3(\e)$ is some function independent of $\g$ such that $\rho_3(\e)\to0$ as $\e\to0$. Then $Q^\e$ %and $(Q^\e)^*$
converges % respectively
to $Q^0$ %and $(Q^0)^*$
in $\fM_{1,0}$ as $\e\to0$ and
\begin{equation}\label{2.21}
\|Q^\e-Q^0\|_{\fM_{1,0}}
\leqslant C\Big(\rho_3^\frac{1}{2}(\e)+\eta^\frac{1}{2}(\e)\Big), 
\end{equation}
where $C$ is some constant independent of $\e$, $\eta$ and $\rho_3$. And vice versa, if a uniformly bounded family $Q^\e\in L_\infty(\Om;\mathds{M}_n)$ converges to $Q^0\in L_\infty(\Om;\mathds{M}_n)$ in $\fM_{1,0}$ as $\e\to0$, then
condition (\ref{2.20}) is satisfied with
 \begin{equation*}
 \eta(\e)= \|Q^\e-Q^0\|_{\fM_{1,0}}^{\frac{1}{2}},\qquad \rho_3(\e)= C\|Q^\e-Q^0\|_{\fM_{1,0}}^{\frac{1}{2}},
\end{equation*}
where $C$ is some constant independent of $\e$. Inequality (\ref{2.28}) holds true.
\end{theorem}

\subsection{Discussion}

Here we discuss our problem and the main results. We first mention main features of the unperturbed operator $\Op$. This is a general matrix second order operator with complex-valued coefficients; conditions (\ref{2.3}) ensure the $m$-sectoriality of this operator once we treat as an unbounded one in $L_2(\Om;\mathds{C}^n)$. It is considered on an arbitrary multi-dimensional domain with a sufficiently regular boundary. The boundary condition is classical and can be either Dirichlet condition or Robin condition, see (\ref{2.4}). We extend this operator to the  space $\fV\subseteq  W_2^1(\Om;\mathds{C}^n)$ since this provides a convenient way for studying general perturbations.

The coefficients of the operator $\Op$ belong to $L_\infty(\Om;\mathds{C}^n)$ and this is sufficient to make the operator and the corresponding sesquilinear form $\hf^\e$  well-defined. We note that by using appropriate embedding theorems it is possible to weaken the conditions on the coefficients and assume that they belong to certain spaces $L_p$ or $L_{p,loc}$ with a sufficiently large $p$ or even to appropriate spaces of multipliers. We do not do this to avoid overloading the text by additional conditions.

The operator $\Op^\e$ involves a perturbation described by the families $P_j^\e$, $Q_j^\e$ and $V^\e$,  where $\e$ is a multi-dimensional small parameter. A formal writing of the perturbation is given in (\ref{2.38}) and this is a first order differential operator and also an additional term in boundary condition (\ref{2.33}). A rigorous definition then is given via the sesquilinear form $\hf^\e$ and it explains the choice of the boundary condition (\ref{2.33}): it ensures that the perturbation does not involve any contribution from the boundary $\p\Om$, just ones on $\Om$.

The main feature of our problem is that apriori we  assume  no specific dependence of these families on $\e$ except for the uniform boundedness. The question we study is what are the minimal conditions on the dependence of $P_j^\e$, $Q_j^\e$ and $V^\e$ on $\e$, under which there exists a limiting (homogenized) operator for $\Op^\e$ in the sense of the norm resolvent convergence. This question is effectively solved in two steps. The results of the first step are presented in Theorems~\ref{th1},~\ref{th5} and the main answer is that the studying of the norm resolvent convergence  \textsl{is equivalent} to studying the convergence of the operator $\cX^\e$ defined by (\ref{2.37}) via $P_j^\e$, $Q_j^\e$ and $V^\e$ in the operator norm $\|\,\cdot\,\|_{\fM}$. This explains why it is convenient to treat the operators $\Op^\e$ and $\Op^0$ as acting from $\fV$ into $\fV^\ast$. In this framework, the perturbation becomes \textsl{regular} and this produces complete expansion (\ref{2.8}) and estimate (\ref{2.9}). And this is the first main message of our work: the perturbation of arbitrary coefficients $P_j^\e$, $Q_j^\e$, $V^\e$ in the operator $\Op^\e$ produces the norm resolvent convergence if and only if the perturbation by the operator $\cX^\e$ is regular in the sense of the norm $\|\,\cdot\,\|_{\fV\to \fV^\ast}$. The fact that the supremum in (\ref{2.22}) involves the norms $\|\,\cdot\,\|_{W_2^1(\Om;\mathds{C}^n)}$ is also very natural since the domains of our operators are $\fV\subseteq W_2^1(\Om;\mathds{C}^n)$.

While the convergence in the operator norm $\|\,\cdot\,\|_{\fV\to \fV^\ast}$ looks quite abstract and complicated, inequality (\ref{2.27}) gives a simple way of estimating this norm by means of the norms in the standard spaces of multipliers $\fM_{1,-1}$ and $\fM_{1,0}$. In view of the above discussed results and this inequality, we can state that the convergence of the families $P_j^\e$ and $Q_j^\e$ in the space  $\fM_{1,0}$ and of the family $V^\e$ in the space $\fM_{1,-1}$ \textsl{is sufficient} to ensure the norm resolvent convergence stated in Theorem~\ref{th1} as well as expansion (\ref{2.8}) and estimate (\ref{2.9}). Moreover, if $P_j^\e=Q_j^\e=0$, then the convergence of the operators $\cX^\e$ in the operator norm $\|\,\cdot\,\|_{\fV\to \fV^\ast}$ \textsl{is equivalent} to the convergence of $V^\e$ in the space $\fM_{1,-1}$. In other words, in this important particular case $P_j^\e=Q_j^\e=0$, the norm resolvent convergence of $\Op^\e$ is equivalent to  the convergence of $V^\e$ in $\fM_{1,-1}$.

The latter fact and the definition of the operator $\Op^\e$ via the form $\hf^\e$ could suggest to suppose that $V^\e$ is not an element of $L_\infty(\Om;\mathds{C}^n)$ but an element of the space $\fM_{1,-1}$. At the same time, %Despite it is possible, this does not enlarge the generality. Indeed,
a general form of a multiplier in $\fM_{1,-1}$ is in fact the operator $\cX^\e$ defined in (\ref{2.37}) provided we assume that there is no boundary term and no singular potentials supported by hypersurfaces and in the matrix potential $V^\e$ each entry  $V_{ij}^\e$ should be such that $(V_{ij}^\e)^\frac{1}{2}$ is a scalar multiplier from $W_2^1(\Om)$ into $L_2(\Om)$, see \cite[Thm. 6.1]{MV}. Treating such boundary term and singular potentials is a more gentle problem, which deserves a separate study and this explains our choice of the perturbation. In other words, instead of assuming that $V^\e$ is a multiplier, we use its general form in our perturbation. The matrix 
functions $P_j^\e$, $Q_j^\e$ and $V^\e$ defining such multiplier are not supposed to belong to the most possible general class but instead we  assume that they are uniformly bounded in $L_\infty(\Om;\mathds{C}^n)$: this is a reasonable condition especially in view of the known results in the case of pure periodically fast oscillating coefficients that the presence of a negative power of $\e$ at coefficients influences essentially the norm resolvent convergence.

In view of the above discussed relations between the convergence in the norm $\|\,\cdot\,\|_{\fM}$  and in the spaces of multipliers $\fM_{1,0}$, it is natural to study the convergence in the latter spaces. Theorem~\ref{th2} provides the main result on the convergence in the space $\fM_{1,-1}$: it is equivalent to the uniform in $\g\in\G_\eta$ convergence of the mean values over small cells $\square_\g^\eta$:
\begin{equation}\label{2.30}
\max\limits_{\g\in\G_\eta}\bigg|\int\limits_{\square_\g^\eta} \eta^{-d}\big(V^\e(x)-V^0(x)\big)\,dx\bigg|\to0,\quad \e\to0,
\end{equation}
for some small function $\eta=\eta(\e)$. It is important to say that this criteria works only under the assumption that the family $V^\e$ belong to $L_\infty(\Om;\mathds{M}_n)$ and is bounded uniformly in this space; this boundedness is employed essentially in the proof of Theorem~\ref{th2}. Namely, condition (\ref{2.13}) is not sensitive to the choice of the boundary conditions in (\ref{2.4}) and the definition of the space $\fV$ and it works for both choices of this space. Such situation is exactly due to the assumed uniform boundedness of the family $V^\e$ since it allows us to neglect the influence of the boundary in the proof of Theorem~\ref{th2}.

We also stress that the presence of the lattice $\G$ and cells $\square_\g^\eta$ in condition (\ref{2.13}) \textsl{does not mean} any periodicity of the family $V^\e$ and they are employed just as an auxiliary tool. The only reason of using a lattice and corresponding (rescaled) cells is that these cells cover the most part of the domain $\Om$ being disjoint at the same time.  The results of Theorems~\ref{th1},~\ref{th5},~\ref{th2} yield that   the  norm resolvent convergence under a perturbation by a potential $V^\e$ \textsl{is equivalent} to the convergence of this potential  in the sense of the mean values in (\ref{2.30}). This is a rather explicit criterion, which can be checked for many examples, which we shall discuss later. Moreover, the convergence in (\ref{2.30}) suggests a natural way of finding a limit for $V^\e$ by formula (\ref{2.16}). Theorem~\ref{th3} states that this is indeed the case once the limit is uniform in $x$; in this theorem the uniform boundedness of $V^\e$ in $L_\infty(\Om;\mathds{M}_n)$ is again an important condition, which is used essentially in the proof. This result and the above discussed ones imply, in particular, the following statement: \textit{if the perturbation is made only by the potential $V^\e$, then the norm resolvent convergence under such perturbation is equivalent to the convergence of the mean values in (\ref{2.30}); the limiting potential $V^0$ can be found as the limit in (\ref{2.16}) provided the convergence of this limit is uniform in $x$.}

Theorem~\ref{th4} provides a criterion for the convergence in the space of multipliers $\fM_{1,0}$. It is similar to one in Theorem~\ref{th2}, but here the convergence of the mean values is not enough and instead, we need the uniform convergence of local averaged $L_2$-norms. In this theorem the uniform boundedness of $Q^\e$ in $L_\infty(\Om;\mathds{C}^n)$ is also employed essentially in the proof. The question on finding a formula similar to (\ref{2.16})
remained open.

In   Theorem~\ref{th6} we consider an important case when the domain of  the operator $\Op^0$ treated as an unbounded one in $L_2(\Om;\mathds{C}^n)$ is a subset of $W_2^2(\Om;\mathds{C}^n)$. Then it is sufficient to assume that $Q_j^\e$ converges in $\fM_{1,-1}$, while
$P_j^\e$ converge in $\fM_{2,0}$ to ensure the norm resolvent convergence. We stress that in the considered case no expansion for the resolvent is provided. The norm of the space $\fM_{2,0}$ is estimated by that in $\fM_{1,0}$, see (\ref{2.39}), and we can use sufficient condition from Theorem~\ref{th4} to guarantee the convergence in $\fM_{1,0}$ and then in $\fM_{2,0}$. However, we failed trying to prove some criteria ensuring the convergence in $\fM_{2,0}$ similar to ones in Theorems~\ref{th2},~\ref{th4}.

We return back to Theorem~\ref{th1} and discuss separately expansion (\ref{2.8}). This is a complete expansion for the resolvent, which simultaneously serves as an asymptotic one. It should be stressed at the same time that our perturbation involves only lower order terms and the coefficients at the second derivatives remain unperturbed. Nevertheless, such complete asymptotic expansion is an essential progress and earlier known results provide a few leading terms in the expansions for the resolvents; the known results are true also for the cases when the higher order terms are perturbed as well \cite{BS1}, \cite{DS}, \cite{Gr1}, \cite{Gr2}, \cite{Kenig}, \cite{ZhP2}, \cite{Sen2}, \cite{Sen3}, \cite{Sen4}, \cite{BS2}, \cite{ZhJi}, \cite{Zhi}, \cite{ZhP1}. At the same time, our result on the complete asymptotic expansion does not cover previous results, for instance, for operators with locally periodic fast oscillating coefficients, even in a particular case of perturbation only in lower order terms, but we rather complete them.  Let us dwell on more details here. The next-to-leading term in expansion (\ref{2.8}) is $(\Op^0-\l)^{-1}  \cL^\e(\Op^0-\l)^{-1}$. This expression means that as a first step in approximating the resolvent of the operator $\Op^0$, we should solve the equation
\begin{equation}\label{2.31}
(\Op^0-\l)u=\cL^\e u_0,\qquad u_0:=(\Op^0-\l)^{-1}f.
\end{equation}
This equation is \textsl{simpler} than the original equation $(\Op^\e-\l)u_\e=f$ because  here the perturbation is moved from the operator into the right hand side. The solution to equation (\ref{2.31}) of course depends on $\e$ and we \textsl{provide no information} about the structure of this dependence since it highly depends on similar structure for $\cL^\e$. For instance, in the case of locally periodic fast oscillating coefficients $Q_j^\e=Q_j(x,\frac{x}{\e})$, $P_j^\e=P_j(x,\frac{x}{\e})$, $V^\e(x)=V(x,\frac{x}{\e})$, our approach says that these coefficients should be moved to the right hand side and then the structure of the  solution to corresponding equation (\ref{2.31}) should be studied in more details. For a domain with a boundary this is a famous long-standing problem on finding an appropriate approximation for arising boundary layers and to the best of the author's knowledge, the very recent progress was made in works \cite{Al1}, \cite{Al2}, \cite{Al3}, \cite{Al4}, where nice estimates for the convergence rates were obtained. Our approach \textsl{does not solve} this problem. But instead, it provides  a complete expansion for the resolvent of $\Op^\e$ in terms of the operator $(\Op^0-\l)^{-1}\cL^\e$ mapping the function $u_0$ in the right hand side of equation (\ref{2.31}) into its solution; the information about the structure of such operator can be obtained by applying the results of works \cite{Al1}, \cite{Al2}, \cite{Al3}, \cite{Al4}.

Summarizing all said above, we can say that our results show \textsl{the borderline of the homogenization} describing quite sharp the class of perturbations for which the norm resolvent convergence is present. To demonstrate this statement, we provide a series of examples. We shall discuss these examples in a separate next section, now we just say that we are able to treat the case of regularly perturbed coefficients, the case of sparsely distributed coefficients, the case of non-periodic stabilizing oscillating coefficients, the cases of locally periodic and almost periodic fast oscillating coefficients, as well as the cases of the coefficients with modulated and fractal periodicity and random perturbation. While dealing with fast oscillating coefficients and checking convergence (\ref{2.30}) of the mean values, the main idea is as follows: if the oscillations have a period of order $O(\e)$, then the function $\eta$ should be chosen small but much larger, namely, $\eta=\e^\a$ with an appropriate $\a<1$. Exactly this allows us to deal with non-periodic oscillations: if they are localized on sets of linear sizes of order $O(\e)$,  by calculating the mean values over still small but much larger sets, we in fact make a local intermediate homogenization mollifying the non-periodic microstructure. Apart of possible examples, we also discuss the ways of generating new treatable families $P_j^\e$, $Q_j^\e$ and $V^\e$ from the given ones.  This allows us to treat finally a very wide class of possible families.

Concluding the discussion, we mention the works on stochastic homogenization of the operators with random oscillating coefficients. The classical results can be found in \cite{ZKO}, while a very recent progress is presented in papers \cite{Stoch2}, \cite{Stoch1}. These models include indeed a wide class of operators with non-periodically oscillating coefficients, however, they still require certain restrictions like stationarity and ergodicity of the considered random fields. The results obtained in stochastic homogenizations state the convergence almost surely, that is, for almost all realizations of the random coefficients.
On the contrary, we make \textsl{no apriori} assumptions on the structure of the our perturbation and our results are of deterministic nature: we say whether it is possible to make a homogenization for a given (realization of) perturbation or not. Nevertheless, our approach can be also applied to the stochastic homogenization, see an example in Section~\ref{StochHom}.

\section{Examples}

In this section  we provide various examples of the families $P_j^\e$,
$Q_j^\e$ and $V^\e$, for which our results can be applied. We also discuss how to generate new families from the given ones.  We stress that in our work the homogenization is treated in the sense of the norm resolvent convergence and this is why our conditions are stronger than ones ensuring the strong or weak resolvent convergence in many known works. We also do not try to get the sharp estimates for the convergence rates as well as to provide the most general conditions for the smoothness of the functions describing the perturbations. Our aim is to demonstrate various classes of perturbations, for which our general approach works.

Before proceeding to the examples, we make one remark. In most part of examples we consider the perturbation of form $V^\e(x)=V(x,x\e_1^{-1},\ldots,x\e_m^{-1})$, where $\e_i$ are small parameters and $V(x,\xi^{(1)},\ldots,\xi^{(m)})$ are some given functions with various dependence on $x$ and $\xi^{(i)}$. Without saying this explicitly, in each example we additionally suppose that the functions
$V^\e(x)$ are measurable if this is not implied immediately by the assumed smoothness of $V$.

\subsection{Regular perturbation}

Our first example is a regular perturbation. Here we assume that
the families $P_j^\e$, $Q_j^\e$ and $V^\e$   converge in $L_\infty(\Om;\mathds{M}_n)$ as $\e\to0$ respectively to $P_j^0$, $Q_j^0$ and $V^0$. Then it follows from (\ref{2.7}), (\ref{2.11}) that
\begin{align*}
&\|V^\e-V^0\|_{\fM_{1,-1}}\leqslant \|V^\e-V^0\|_{L_\infty(\Om;\mathds{M}_n)}\to0,
\\
&\| Q_j^\e - Q_j^0 \|_{\fM_{1,0}}\leqslant \|Q_j^\e-Q_j^0\|_{L_\infty(\Om;\mathds{M}_n)}\to 0,
\\
&\|P_j^\e-P_j^0\|_{\fM_{1,0}}\leqslant \|P_j^\e-P_j^0\|_{L_\infty(\Om;\mathds{M}_n)}\to 0
\end{align*}
as $\e\to0$. Then, in view of estimate (\ref{2.27}),
the assumptions of  Theorem~\ref{th1} are satisfied and expansion (\ref{2.8}) reproduces a classical expansion of the resolvent in the case of a usual regular perturbation.

\subsection{Sparse distribution}

Let $M_k^\e$ be some points in $\Om$ indexed by $k\in\mathds{I}^\e$, where $\mathds{I}^\e$ is some at most countable set. Suppose that
\begin{equation}\label{5.5}
\inf\limits_{\substack{k,j\in\mathds{I}^\e
\\
k\ne j}}\dist(M_k^\e,M_j^\e)\geqslant \rho_4(\e)>0,\qquad  \rho_4(\e)\to0,\quad \e\to0.
\end{equation}
By $B_r(a)$ we denote a ball in $\mathds{R}^d$ of radius $r$ centered at a point $a$. For each $k\in\mathds{I}^\e$ we arbitrarily choose a $n\times n$ matrix function $V_k^\e=V_k^\e(x)$ such that it vanishes outside $B_{\rho_4(\e)\rho_5(\e)}(M_k^\e)$, where $\rho_5=\rho_5(\e)$ is some function independent of $k\in\mathds{I}^\e$ such that $\rho_5(\e)\to0$ as $\e\to 0$. We also suppose that the matrix functions $V_k^\e$ belong to $L_\infty(B_{\rho_4\rho_5}(M_k^\e);\mathds{M}_n)$ for each $k\in\mathds{I}^\e$ and are bounded in the norm of this space uniformly in $k\in\mathds{I}^\e$ and $\e$. The function $V^\e$ is defined as $V^\e(x)=V_k^\e(x)$ on $B_{\rho_4\rho_5}(M_k^\e)$, $k\in\mathds{I}^\e$, and $V^\e(x)=0$ otherwise. A sketched graph of the function $V^\e$  in the case of the sparse distribution  is given on Figure~\ref{fig:sparse}.

Then we choose $\G:=\mathds{Z}^d$, $\eta:=\frac{1}{3}\rho_4$ and introduce the sets $\G_\eta$ and $\square_\g^\eta$ by (\ref{2.18}). Condition (\ref{5.5}) then ensures that each cell $\square_\g^\eta$ intersects with at most one of the balls $B_{\rho_4\rho_5}(M_k^\e)$ and in view of the uniform boundedness of $V_k^\e$ we then obtain:
\begin{align*}
\eta^{-d}\bigg|\int\limits_{\square_\g^\eta} V^\e(x)\,dx\bigg|\leqslant C\rho_4^{-d} \mes B_{\rho_4\rho_5}(0)\leqslant C\rho_5^d, \qquad \eta^{-d}\int\limits_{\square_\g^\eta} \big|V^\e(x)\big|^2\,dx\leqslant   C\rho_5^d,
\end{align*}
where $C$ are some constants independent of $\e$, $\rho_4$, $\rho_5$ and $\g$. By Theorems~\ref{th2},~\ref{th4} we then conclude that the family $V^\e$ converges to zero both in $\fM_{1,-1}$ and $\fM_{1,0}$ and
\begin{equation*}
\|V^\e\|_{\fM_{1,-1}}\leqslant C\Big(\rho_5^d(\e) + \rho_4(\e)\Big), \qquad \|V^\e\|_{\fM_{1,0}}\leqslant C\Big(
 \rho_5^{\frac{d}{2}}(\e)    + \rho_4^{\frac{1}{2}}(\e)\Big),
\end{equation*}
with constants $C$ independent of $\e$, $\rho_4$ and $\rho_5$.

\begin{figure}[t]
\begin{center}
\includegraphics[scale=0.5]{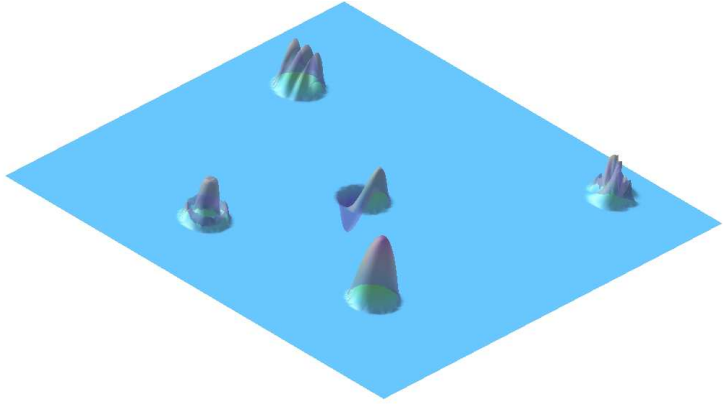}
\caption{Coefficient with sparse distribution}\label{fig:sparse}
\end{center}
\end{figure}

\subsection{Non-periodic stabilizing oscillations}

Let $V=V(x,\xi)$ be a $n\times n$ matrix function defined on $\Om\times\mathds{R}^d$, measurable, uniformly bounded  and having a finite limit as $\xi\to\infty$. Namely, we assume that there exists a $n\times n$ matrix function $V^0=V^0(x)$ belonging to $L_\infty(\Om;\mathds{M}_n)$ such that $V(x,\xi)\to V^0(x)$ as $\xi\to\infty$ uniformly in $x\in\Om$.

We let $\e=(\e_1,\ldots,\e_d)$, where $\e_j$, $j=1,\ldots,d$, are small positive parameters, and we define
\begin{equation*}
V^\e(x):=V\bigg(x,\frac{x_1}{\e_1},\ldots, \frac{x_d}{\e_d}\bigg).
\end{equation*}
Such family is obviously uniformly bounded in $L_\infty(\Om;\mathds{M}_n)$ and we are going to show that it also converges to $V^0$ in $\fM_{1,-1}$  as $\e\to0$.

We denote: $\e_{\max}:=\max\limits_{j=1,\ldots,d} \e_j$, $\eta(\e):=\e_{\max}^\frac{1}{3}$. Then we choose the lattice  $\G$ in $\mathds{R}^d$ as $\G:=(-1,\ldots,-1)+2\mathds{Z}^d$; its periodicity cell is $\square:=(0,2)^d$. Let us show that condition (\ref{2.13}) holds for the family $V^\e$. It is clear that
\begin{equation}\label{5.24}
\frac{x_1^2}{\e_1^2}+\ldots+ \frac{x_d^2}{\e_d^2}\geqslant \frac{|x|^2}{\e_{\max}^2}
\end{equation}
and hence, due to the uniform convergence of $V(x,\xi)$ to its limit $V^0(x)$ as $\xi\to\infty$,
\begin{equation}\label{5.13}
|V^\e(x)-V^0(x)|\leqslant \rho_6(\e)\quad\text{as}\quad |x|>\eta^2(\e),\qquad
\rho_6(\e):=\sup\limits_{
\substack{|\xi|\geqslant \eta^{-1}(\e)
 \\
 x\in\Om}} |V(x,\xi)-V^0(x)|\to 0,\qquad \e\to0.
\end{equation}
Then, in view of the uniform boundedness of the family $V^\e$, we have:
\begin{equation}\label{5.14}
\begin{aligned}
\frac{1}{\eta^d} \bigg|\int\limits_{\eta \square} \big(V^\e(x)-V^0(x)\big) \,dx\bigg|\leqslant &
\frac{1}{\eta^d} \int\limits_{B_{\eta^2}(0)} |V^\e(x)-V^0(x)|\,dx + \frac{1}{\eta^d}
\int\limits_{\eta\square\setminus B_{\eta^2}(0)} |V^\e(x)-V^0(x)|\,dx
\\
\leqslant & C\eta^d + \rho_6,
\end{aligned}
\end{equation}
where $C$ is some constant independent of $\e$, $\eta$ and $\rho_6$. On other cells $\eta\square+\eta\g$ with $\g\ne 0$ estimate (\ref{5.13}) holds and hence,
\begin{equation*}
\frac{1}{\eta^d} \bigg|\int\limits_{\eta \square+\eta\g} \big(V^\e(x)-V^0(x)\big) \,dx\bigg|\leqslant \rho_6.
\end{equation*}
This estimate and (\ref{5.14}) yield that Theorem~\ref{th2} holds with the above introduced $\eta$ and
\begin{equation*}
\|V^\e-V^0\|_{\fM_{1,-1}}\leqslant C\big(\rho_6(\e)+\e_{\max}^{\frac{1}{3}}\big)
\end{equation*}
with a constant $C$ independent of $\e$.
In the same way we confirm that Theorem~\ref{th4} also holds for the family $V^\e$, it converges to $V^0$ in $\fM_{1,0}$ and
\begin{equation*}
\|V^\e-V^0\|_{\fM_{1,0}}\leqslant C\big(\rho_6^\frac{1}{2}(\e)+\e_{\max}^{\frac{1}{6}}\big)
\end{equation*}
with a constant $C$ independent of $\e$.

It is also possible to modify the above example by considering the functions $V(x,\xi)$ with a varying limit at infinity. Namely, let $\mathds{S}^{d-1}$ be the unit sphere in $\mathds{R}^d$. Suppose that we are again given a $n\times n$ matrix function $V=V(x,\xi)$, which is measurable and uniformly bounded and  there exists a $n\times n$ matrix function $V^0=V^0(x,\z)$ defined on $\Om\times \mathds{S}^{d-1}$, measurable and uniformly bounded such that
\begin{equation*}
V(x,t\z)\to V^0(x,\z)\quad\text{as}
\quad t\to+\infty\quad\text{uniformly in} \quad (x,\z)\in\Om\times\mathbb{S}^{n-1}.
\end{equation*}
Then we choose a single small positive parameter $\e$ and we  define $V^\e(x):=V\big(x,\frac{x}{\e}\big)$.  A sketched graph of the function $V^\e$ as $d=2$,
\begin{equation*}
V(x,\xi)=\sin \frac{5\xi_1}{|\xi|} + \cos\frac{10\xi_1}{|\xi|}
+ \frac{\sin 5\xi_2^2}{1+|\xi_1|+2|\xi_2|}+ \frac{\cos 7\xi_1}{1+|\xi_1|+\xi_2^2}
\end{equation*}
is given on Figure~\ref{fig:stab}.

Let us confirm that this family converges to $V^0(x):=V^0(x,\tfrac{x}{|x|})$ as $\e\to0$. We again choose the same lattice $\G$ as above and $\eta:=\e^\frac{1}{3}$. Estimate (\ref{5.24}) remains true with $\e_{\max}=\e_1=\ldots=\e_m=\e$ and relations (\ref{5.13}) are to be modified as follows:
\begin{equation*}
|V^\e(x)-V^0(x)|\leqslant \rho_7(\e)\quad\text{as}\quad |x|>\e^{-\frac{2}{3}},\qquad
\rho_7(\e):=\sup\limits_{
\substack{|\xi|\geqslant \e^{-\frac{1}{3}}
 \\
 x\in\Om}} \bigg|V(x,\xi)-V^0\bigg(x,\frac{\xi}{|\xi|}\bigg) \bigg|\to 0,\qquad \e\to0.
\end{equation*}
All other calculations remain unchanged and we then obtain:
\begin{equation*}%\%label{5.25}
\|V^\e-V^0\|_{\fM_{1,-1}}\leqslant C\big(\rho_7(\e)+\e^{\frac{1}{3}}\big),
\qquad \|V^\e-V^0\|_{\fM_{1,0}}\leqslant C\big(\rho_7^\frac{1}{2}(\e)+\e^{\frac{1}{6}}\big).
\end{equation*}

\subsection{Locally periodic fast oscillations}\label{sec:locper}

Here we consider a classical example of a locally periodic fast oscillating coefficients. Namely, we choose a multi-dimensional small parameter $\e=(\e_1,\ldots,\e_m)$ and $m$ arbitrary fixed lattices $L_j$, $j=1,\ldots,m$, in $\mathds{R}^d$ with periodicity cells $\square_j$.
By $V=V(x,\xi)$, $\xi:=(\xi_1,\ldots,\xi_m)$, we denote an arbitrary $n\times n$ matrix function defined for $x\in\Om$ and $\xi_j\in\mathds{R}^d$, $j=1,\ldots,m$; here each $\xi_j$ is a $d$-dimensional variable and $\xi$ is a $md$-dimensional variable. We assume that the function $V$ is uniformly bounded and $\square_j$-periodic in $\xi_j$ for each $j=1,\ldots,m$.
We also assume the following uniform continuity:
\begin{equation}\label{5.15}
\sup\limits_{\substack{x,\tilde{x}\in \Om,\, \xi',\tilde{\xi'}\in\square_0
\\
|x-\tilde{x}|\leqslant \d,\, |\xi'-\tilde{\xi}'|\leqslant \d}}  \|V(x,\xi',\,\cdot\,)- V(\tilde{x},\tilde{\xi}',\,\cdot\,)\|_{L_\infty(\square_m)}
%\sup\limits_{\xi_m\in\square_m}
%|V(x,\xi)-V(x,\tilde{\xi})|
\leqslant \rho_8(\d),\qquad \rho_8(\d)\to0,\qquad \d\to0,
\end{equation}
where $\rho_8=\rho_8(\d)$ is some   function, $\xi'=(\xi_1,\ldots,\xi_{m-1})$, $\square_0:=\square_1\times\ldots \times\square_{m-1}$. In addition we assume that all $\e_j$ are positive and tend to zero so that
\begin{equation}\label{5.12}
\frac{\e_2}{\e_1}\to0,\qquad \frac{\e_3}{\e_2}\to0,\qquad \ldots \qquad
\frac{\e_m}{\e_{m-1}}\to0.
\end{equation}
The family $V^\e$ is defined as
\begin{equation}\label{5.23}
V^\e(x):=V\bigg(x,\frac{x}{\e_1},\ldots, \frac{x}{\e_m}\bigg).
\end{equation}
A sketched graph of the function $V^\e$ for $d=2$, $\e=(\e_1,\e_1^2)$,  $V(x,\xi,\z)=\sin \xi_1 \cos \xi_2 +\cos \z_1 \sin \z_2 $  is given on Figure~\ref{fig:iter-per}.

\begin{figure}[t]
\begin{center}
\includegraphics[scale=0.5]{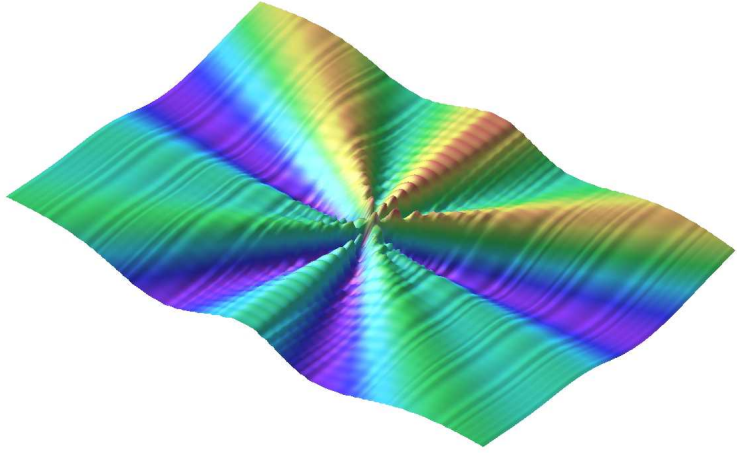}
\caption{Coefficient with stabilizing oscillations}\label{fig:stab}
\end{center}
\end{figure}

We are going to show that the family $V^\e$ converges in $\fM_{1,-1}$ as $\e\to0$ and the limit is
\begin{equation*}
V^0(x):=\frac{1}{\prod\limits_{j=1}^{m}\mes\square_j}
\int\limits_{\square_1}d\xi_1\ldots\int\limits_{\square_m}d\xi_m V(x,\xi_1,\ldots,\xi_m).
\end{equation*}
In order to do this, we shall confirm that the assumptions of Theorem~\ref{th2} are satisfied with  $\eta:=\e_1^\frac{1}{2}$ and $\G:=\mathds{Z}^d$, $\square:=(0,1)^d$.

We introduce the sets $\G_\eta$ and $\square_\g^\eta$ by (\ref{2.18}),
 fix an arbitrary $\g\in\G_\eta$   and denote
\begin{equation}\label{5.21}
L_{m,\e_m}:=\big\{\g_m\in L_m:\ \e_m\square_m+\e_m\g_m\subset \square_\g^\eta\big\}, \qquad \square_m^{(\e_m)}:=\bigcup\limits_{\g_m\in L_{m,\e_m}} (\e_m\square_m+\e_m\g_m).
\end{equation}
It follows from the above definition that
\begin{equation}\label{5.17}
\mes  \square_\g^\eta \setminus\square_m^{(\e_m)} %+ \mes\square_m^{(\e_m)}\setminus \square_\g^\eta
\leqslant C\eta^{d-1}\e_m,
\end{equation}
where $C$ is some constant independent of $\g$, $\eta$ and $\e$. Hence, due to the uniform boundedness of $V$,
\begin{equation}\label{5.16}
\bigg| \frac{1}{\eta^d} \int\limits_{\square_\g^\eta} V^\e(x)\,dx - \frac{1}{\eta^d}\int\limits_{\square_m^{(\e_m)}} V^\e(x)\,dx\bigg| \leqslant C\e_m \eta^{-1}\leqslant C\e_1^\frac{1}{2}
\end{equation}
with a constant $C$ independent of $\g$, $\eta$ and $\e$. For each cell $\subset \square_m^{(\e_m)}$ we obviously have
\begin{equation*}
|x-\e_m \g_m|\leqslant k_m \e_m,\qquad x\in \e_m\square_m + \e_m\g_m,
\end{equation*}
where $k_m$ is the linear size of the cell $\square_m$. Then owing to conditions (\ref{5.12}) we obtain
\begin{equation*}
\bigg|\bigg(x,\frac{x}{\e_1},\ldots,\frac{x} {\e_{m-1}}\bigg)-
\bigg(\e_m\g_m,\frac{\e_m}{\e_1}\g_m,\ldots, \frac{\e_m}{\e_{m-1}}\g_m\bigg)\bigg|\leqslant   k_m \e_m \bigg(1+\frac{1}{\e_1^2}+ \ldots + \frac{1}{\e_{m-1}^2}  \bigg)^{\frac{1}{2}}
\leqslant  \sqrt{m}k_m \frac{\e_m}{\e_{m-1}}
\end{equation*}
as $x\in \e_m \square_m + \e_m\g_m$.
Hence, by condition (\ref{5.15}),
\begin{equation}\label{5.18}
\bigg|V^\e(x) - V\bigg(\e_m\g_m,\frac{\e_m}{\e_1}\g_m,\ldots, \frac{\e_m}{\e_{m-1}}\g_m,\frac{x}{\e_m}\bigg) \bigg| \leqslant \rho_8\bigg(\sqrt{m}k_m \frac{\e_m}{\e_{m-1}}\bigg),
\quad x\in\e_m \square_m + \e_m\g_m.
\end{equation}
\begin{figure}[t]
\begin{center}
\includegraphics[scale=0.5]{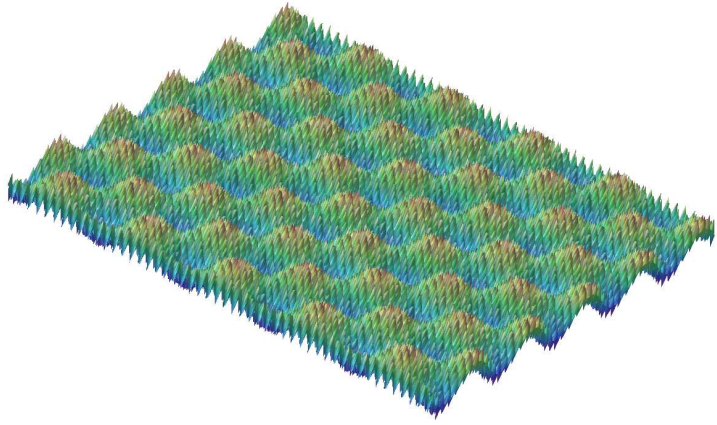}
\caption{Coefficient with multi-scale periodicity}\label{fig:iter-per}
\end{center}
\end{figure}
This estimate implies:
\begin{align*}
\Bigg|\frac{1}{\eta^d}& \int\limits_{\square_m^{(\e_m)}} V^\e(x)\,dx - \frac{1}{\eta^d} \sum\limits_{\g_m\in L_{m,\e_m}} \int\limits_{\e_m\square_m + \e_m\g_m} V\bigg(\e_m\g_m,\frac{\e_m}{\e_1}\g_m,\ldots, \frac{\e_m}{\e_{m-1}}\g_m,\frac{x}{\e_m}\bigg) \,dx
\Bigg|
\\
&\leqslant  \frac{1}{\eta^d} \sum\limits_{\g_m\in L_{m,\e_m}} \int\limits_{\e_m\square_m + \e_m\g_m}
\bigg| V^\e(x)-V\bigg(\e_m\g_m,\frac{\e_m}{\e_1}\g_m,\ldots, \frac{\e_m}{\e_{m-1}}\g_m,\frac{x}{\e_m}\bigg)
\bigg|\,dx
\\
&\leqslant\frac{1}{\eta^d} \rho_8\bigg(\sqrt{m}k_m \frac{\e_m}{\e_{m-1}}\bigg) \sum\limits_{\g_m\in L_{m,\e_m}} \mes (\e_m\square_m + \e_m\g_m)= \frac{\mes \square_m^{(\e_m)}}{\eta^d} \rho_8\bigg(\sqrt{m}k_m \frac{\e_m}{\e_{m-1}}\bigg)
\\
&\leqslant\rho_8\bigg(\sqrt{m}k_m \frac{\e_m}{\e_{m-1}}\bigg)
\end{align*}
and by (\ref{5.16}) we hence get:
\begin{equation}\label{5.22}
\begin{aligned}
\Bigg|\frac{1}{\eta^d}\int\limits_{\square_\g^\eta} V^\e(x)\,dx &- \frac{1}{\eta^d} \sum\limits_{\g_m\in L_{m,\e_m}} \int\limits_{\e_m\square_m + \e_m\g_m} V\bigg(\e_m\g_m,\frac{\e_m}{\e_1}\g_m,\ldots, \frac{\e_m}{\e_{m-1}}\g_m,\frac{x}{\e_m}\bigg) \,dx\Bigg|
\\
&\leqslant C\bigg(\rho_8\bigg(\sqrt{m}k_m \frac{\e_m}{\e_{m-1}}\bigg)+\frac{\e_m}{\eta}\bigg),
\end{aligned}
\end{equation}
where $C$ is some constant independent of $\e$ and the point $\g$ involved in definition (\ref{5.21}) of $\square_m^{(\e_m)}$. At the same time,
\begin{equation}\label{5.19}
\begin{aligned}
\int\limits_{\e_m\square_m+\e_m\g_m} V\bigg(\e_m\g_m,&\frac{\e_m}{\e_1}\g_m,\ldots, \frac{\e_m}{\e_1^{m-1}}\g_m,\frac{x}{\e_m}\bigg) \,dx
\\
=& \e_m^d \int\limits_{\square_m} V\bigg(\e_m\g_m,\frac{\e_m}{\e_1}\g_m,\ldots, \frac{\e_m}{\e_{m-1}}\g_m,\xi_m\bigg)\,d\xi_m
\\
=&\frac{1}{\mes\square_m} \int\limits_{\e_m \square_m+\e_m\g_m} \,dx \int\limits_{\square_m} V\bigg(\e_m\g_m,\frac{\e_m}{\e_1}\g_m,\ldots, \frac{\e_m}{\e_{m-1}}\g_m,\xi_m\bigg)\,d\xi_m
\\
=& \frac{1}{\mes\square_m} \int\limits_{\e_m \square_m+\e_m\g_m} \,dx \int\limits_{\square_m} V\bigg(x,\frac{x}{\e_1},\ldots, \frac{x}{\e_{m-1}},\xi_m\bigg)\,d\xi_m
\\
&+ \frac{1}{\mes\square_m} \int\limits_{\e_m \square_m+\e_m\g_m} \,dx \int\limits_{\square_m} \bigg( V\bigg(\e_m\g_m,\frac{\e_m}{\e_1}\g_m,\ldots, \frac{\e_m}{\e_{m-1}}\g_m,\xi_m\bigg)
\\
&\hphantom{+ \frac{1}{\mes\square_m} \int\limits_{\e_m \square_m+\e_m\g_m} \,dx \int\limits_{\square_m}\bigg(} - V\bigg(x,\frac{x}{\e_1},\ldots, \frac{x}{\e_{m-1}},\xi_m\bigg)
\bigg)
\,d\xi_m
\end{aligned}
\end{equation}
and by (\ref{5.18}) we have:
\begin{align*}
\frac{1}{\mes\square_m}\Bigg|
\int\limits_{\e_m \square_m+\e_m\g_m} \,dx \int\limits_{\square_m} \bigg( V\bigg(\e_m\g_m,\frac{\e_m}{\e_1}\g_m,\ldots, \frac{\e_m}{\e_{m-1}}\g_m,\xi_m\bigg)
&- V\bigg(x,\frac{x}{\e_1},\ldots, \frac{x}{\e_{m-1}},\xi_m\bigg)
\bigg)\,d\xi_m\Bigg|
\\
 &\leqslant C \e_m^d \rho_8\bigg(\sqrt{m}k_m \frac{\e_m}{\e_{m-1}}\bigg)
\end{align*}
with a constant $C$ independent of $\e$ and $\g$. This estimate, (\ref{5.17}), the boundedness of $V$ and relations (\ref{5.19}), (\ref{5.22}) yield:
\begin{equation*}%\%label{5.20}
\Bigg| \frac{1}{\eta^d} \int\limits_{ \square_\g^\eta}V^\e(x)\,dx - \frac{1}{\eta^d} \int\limits_{\square_\g^\eta} \frac{dx}{\mes \square_m} \int\limits_{\square_m} V\bigg(x,\frac{x}{\e_1},\ldots, \frac{x}{\e_{m-1}},\xi_m\bigg)
\,d\xi_m
\Bigg|\leqslant C \bigg(\rho_8\bigg(\sqrt{m}k_m \frac{\e_m}{\e_{m-1}}\bigg) +\frac{\e_m}{\eta}\bigg),
\end{equation*}
where the constant $C$ is independent of $\e$ and $\g$.

We introduce auxiliary function and sets
\begin{align*}
&V_1(x,\xi_1,\ldots,\xi_{m-1}):= \frac{1}{\mes\square_m} \int\limits_{\square_m} V(x,\xi_1,\ldots,\xi_{m-1},\xi_m)\,d\xi_m,
\\
&L_{m-1,\e_{m-1}}:=\big\{\g_{m-1}\in L_{m-1}:\ \e_{m-1}\square_{m-1} + \e_{m-1}\g_{m-1}\subset \square_\g^\eta\big\},
\\
&\square_{m-1}^{(\e_{m-1})}:= \bigcup\limits_{\g_{m-1}\in L_{m-1,\e_{m-1}}} (\e_{m-1}\square_{m-1}+\e_{m-1}\g_{m-1}).
\end{align*}
The function $V_1$ obviously possesses the same properties as $V$, while the set $\square_{m-1}^{(\e_{m-1})}$ satisfies an estimate similar to (\ref{5.17}):
\begin{equation*}
\mes \square_\g^\eta \setminus\square_{m-1}^{(\e_{m-1})}
\leqslant C\eta^{d-1}\e_{m-1}.
\end{equation*}
Arguing then as above, we easily prove that
\begin{align*}
\Bigg|\frac{1}{\eta^d}& \int\limits_{ \square_\g^\eta} \frac{dx}{\mes \square_m} \int\limits_{\square_m} V\bigg( x,\frac{x}{\e_1},\ldots,\frac{x}{\e_{m-1}},\xi_m
\bigg)\,d\xi_m
\\
&- \frac{1}{\eta^d} \int\limits_{\square_\g^\eta} \frac{dx}{\mes \square_m\mes\square_{m-1}} \int\limits_{\square_m}d\xi_m \int\limits_{\square_{m-1}}d\xi_{m-1} V\bigg( x,\frac{x}{\e_1},\ldots,\frac{x}{\e_{m-2}}, \xi_{m-1},\xi_m\bigg)
\Bigg|
\\
=&\Bigg|\frac{1}{\eta^d} \int\limits_{\square_\g^\eta} V_1\bigg(x,\frac{x}{\e_1},\ldots,\frac{x}{\e_{m-1}}
\bigg)\,dx -  \frac{1}{\eta^d} \int\limits_{ \square_\g^\eta} \frac{dx}{\mes\square_{m-1}}  \int\limits_{\square_{m-1}}d\xi_{m-1} V_1\bigg(x,\frac{x}{\e_1},\ldots,\frac{x}{\e_{m-2}}, \xi_{m-1}\bigg)
\Bigg|
\\
\leqslant& C \bigg(\rho_8\bigg(\sqrt{m-1}k_{m-1} \frac{\e_{m-1}}{\e_{m-2}}\bigg) + \frac{\e_{m-1}}{\eta}
\bigg),
\end{align*}
where $C$ is some constant independent of $\g $ and $\e$, while $k_{m-1}$ is the linear size of the cell $\square_{m-1}$. Repeating recurrently this  procedure and summing up the obtained estimates, we finally get:
\begin{equation*}
\Bigg|\frac{1}{\eta^d} \int\limits_{ \square_\g^\eta} V^\e(x)\,dx - \frac{1}{\eta^d} \int\limits_{ \square_\g^\eta} V^0(x)\,dx\Bigg| \leqslant  C \sum\limits_{j=2}^{m} \bigg( \rho_8\bigg( \sqrt{j} k_j \frac{\e_j}{\e_{j-1}}
\bigg) + \frac{\e_j}{\eta}
\bigg)
\leqslant  C \bigg(\sum\limits_{j=2}^{m} \rho_8\bigg( \sqrt{j} k_j \frac{\e_j}{\e_{j-1}}\bigg) + \e_1^\frac{1}{2}\bigg),
\end{equation*}
where $C$ is a constant independent of $\e$ and $\g$. Since $\g$ is arbitrary, the above estimate ensures that the assumptions of Theorem~\ref{th2} are satisfied and the family $V^\e$ converges to $V^0$ in $\fM_{1,-1}$:
\begin{equation*}
\|V^\e-V^0\|_{\fM_{1,-1}}\leqslant  C \Bigg(\sum\limits_{j=2}^{m} \rho_8\bigg( \sqrt{j} k_j \frac{\e_j}{\e_{j-1}}\bigg) + \e_1^\frac{1}{2}\Bigg).
\end{equation*}

In the case $m=1$, function (\ref{5.23}) is a classical example of a locally periodic fast oscillating coefficient. If $m>1$, it covers also the case $\e_i=\e^{\a_i}$ with some exponents $\a_1<\a_2<\ldots<\a_m$. This is also a classical example of the reiterated homogenization. However, our results impose no strict restrictions on $\e_i$, for instance, they are applicable to the case when all $\e_i$ are arbitrary functions of a single small parameter, say, $\e_1=\e$, $\e_2=\e^2$, $\e_3=e^{-\frac{1}{|\e|}}$, $\e_4=e^{-\frac{1}{\e^2}}$, $\e_5=e^{-e^{\frac{1}{\e^2}}}$. Moreover, the parameters $\e_i$ can vary independently and they just should satisfy (\ref{5.12}). To the best of the author's knowledge, such multi-parametric perturbation have never been studied before.

\subsection{Locally almost periodic fast oscillations}

This example is devoted to a locally almost periodically fast oscillating coefficient. The homogenization of operators with almost periodic fast oscillating coefficients is a  classical example, which was intensively studied, see, for instance, \cite{Ko}, \cite{ZhJi}, and the references therein.

It is introduced similar to the case of locally periodically fast oscillating coefficient, but for simplicity, we suppose now that $\e$ is a one-dimensional small parameter and $V=V(x,\xi)$. The $n\times n$  matrix function $V$ is   defined for $x\in \Om$, $\xi\in\mathds{R}^d$, is measurable, uniformly bounded and uniformly continuous in $x$:
\begin{equation}\label{5.26}
\sup\limits_{\substack{x,\tilde{x}\in\Om
\\
|x-\tilde{x}|\leqslant \d}} \sup\limits_{\xi\in\mathds{R}^d}|V(x,\xi)-V(\tilde{x},\xi)| \leqslant \rho_9(\d),\qquad \rho_9(\d)\to0,\quad \d\to0.
\end{equation}
The dependence of $V$ on $\xi$ is supposed to be almost periodic. Namely, for each $x$ the matrix function $V(x,\,\cdot\,)$ belongs to the Weyl class of almost periodic matrix functions. This class is denoted by $\mathds{W}(\mathds{R}^d;\mathds{M}_n)$ and is introduced as follows.
We first define a seminorm
\begin{equation}\label{5.27}
\|f\|_{\mathds{W}(\mathds{R}^d;\mathds{M}_n)}:=\lim\limits_{r\to+\infty}
\sup\limits_{\g\in\G} \frac{1}{r^d\mes\square}\int\limits_{ r\square+ r\g} |f(\xi)|\,d\xi
\end{equation}
on an appropriate subspace in $L_\infty(\mathds{R}^d;\mathds{M}_n)$, where $\G$ is some fixed lattice with a periodicity cell $\square$. The Weyl class $\mathds{W}(\mathds{R}^d;\mathds{M}_n)$ is introduced as a completion, in the above seminorm, of $n\times n$ matrix functions, the entries of which are finite trigonometric polynomials.  The family $V^\e$ is defined as
\begin{equation*}
V^\e(x):=V\left(x,\frac{x_1}{\e},\ldots,\frac{x_d}{\e}\right).
\end{equation*}
A sketched graph of the function $V^\e$ as $d=2$, $V(x,\xi)=\sin\xi_1\cos\xi_2+\cos \sqrt{2}\xi_1 \cos\sqrt{3}\xi_2$  is given on Figure~\ref{fig:alm-per}.

The definition of the Weyl almost periodic function ensures that the considered functions possess a mean value:
\begin{equation}\label{5.29}
f^0:=\lim\limits_{r\to+\infty} \frac{1}{r^d\mes\square} \int\limits_{r\square+r\g } f(\xi)\,d\xi,
\end{equation}
which turns out to be independent of $\g\in\G$. Indeed, for each $\d>0$ there exists a trigonometric polynomial $T_\d=T_\d(\xi)$ such that $\|f-T_\d\|_{\mathds{W}(\mathds{R}^d;\mathds{M}_n)}<\d$. Then definition (\ref{5.27}) implies the existence of $R(\d)$ such that for all $r>R(\d)$ and $\g\in\G$ we have
\begin{equation}\label{5.27a}
\begin{aligned}
\bigg|\frac{1}{r^d \mes\square}\int\limits_{r\square+r\g}f(\xi)\,d\xi - \frac{1}{r^d \mes\square}\int\limits_{r\square+r\g}T_\d(\xi)\,d\xi
\bigg| \leqslant & \frac{1}{r^d\mes \square} \int\limits_{r\square+r\g} |f(\xi)-T_\d(\xi)|\,d\xi
\\
\leqslant & \|f-T_\d\|_{\mathds{W}(\mathds{R}^d;\mathds{M}_n)}<\d.
\end{aligned}
\end{equation}
Hence,
\begin{equation*}
\bigg|\frac{1}{r_1^d \mes\square}\int\limits_{r_1\square+r_1\g}f(\xi)\,d\xi - \frac{1}{r_2^d \mes\square}\int\limits_{r_2\square+r_2\g}f(\xi)\,d\xi
\bigg| \leqslant 2\d\quad\text{for}\quad r_1,r_2>R(\d)
\end{equation*}
and this proves the existence of the limit in (\ref{5.29}) for each $\g\in\G$. Since $T_\d$ is a finite trigonometric polynomial, that is, $T_\d(\xi)=\sum\limits_{\a} e^{\iu \a\cdot \xi} T_{\d,\a}$, where the sum is finite, $\a\in\mathds{R}^d$ are some vectors, $T_{\d,\a}$ are some constant $n\times n$ matrices, we see that the limit
\begin{equation*}
T_\d^0:=\lim\limits_{r\to+\infty} \frac{1}{r^d\mes\square} \int\limits_{r\square+r\g } T_\d(\xi)\,d\xi
\end{equation*}
is well-defined for each $\g\in\G$ and is independent of $\g$. Then it follows from (\ref{5.27a}), (\ref{5.29}) that $|f^0-T_\d^0|<\d$ for each $\g\in\G$. Therefore, the limit in (\ref{5.29}) is independent of $\g$ and it converges uniformly in $\g$.

\begin{figure}[t]
\begin{center}
\includegraphics[scale=0.5]{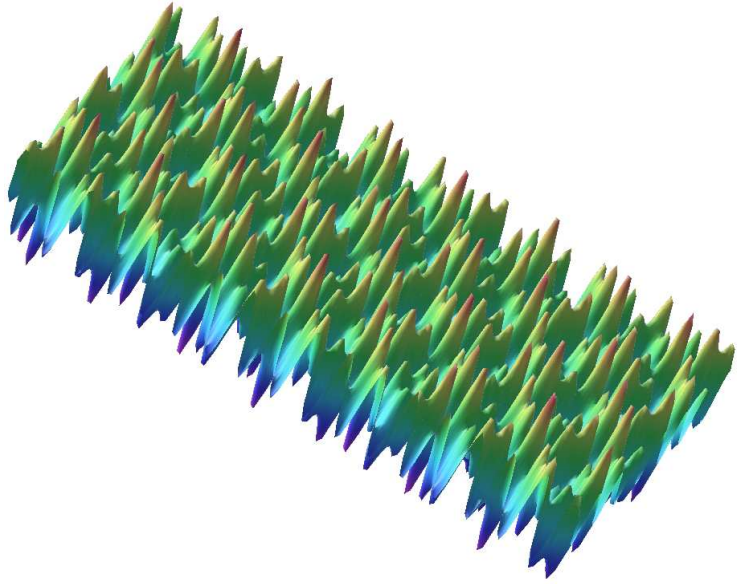}
\caption{Coefficient with almost periodic oscillations}\label{fig:alm-per}
\end{center}
\end{figure}

Since $V(x,\,\cdot\,)\in \mathds{W}(\mathds{R}^d;\mathds{M}_n)$, for each $x\in \Om$ the mean value
\begin{equation}\label{5.30}
V^0(x):= \lim\limits_{r\to+\infty} \frac{1}{r^d\mes\square} \int\limits_{r\square+r\g } V(x,\xi)\,d\xi
\end{equation}
is well-defined, independent of $\g$ and the limit converges uniformly in $\g$. In addition, we assume that the convergence in the above limit is uniform in $x\in\Om$:
\begin{equation}\label{5.28}
\bigg|\frac{1}{r^d\mes\square} \int\limits_{r\square+ r\g} V(x,\xi)\,d\xi - V^0(x)\bigg| \leqslant \rho_{10} (r)\quad\text{for all}\quad x\in\Om,\quad \g\in\G,
\end{equation}
where $\rho_{10}=\rho_{10}(r)$ is some function independent of $x$ and $\g$ and $\rho_{10}(r)\to0$ as $r\to+\infty$.

We are going to show that under the above assumptions, condition (\ref{2.13}) is satisfied with the matrix function $V^0$ introduced in (\ref{5.30}), the lattice $\G$ used in (\ref{5.27}) and $\eta=\e^\frac{1}{2}$. First of all we observe that since the matrix function $V$ is uniformly bounded, definition (\ref{5.30}) implies the same for $V^0$. Hence, $V^\e, V^0\in L_\infty(\Om;\mathds{M}_n)$ and $V^\e$ is bounded uniformly in $\e$.

We introduce the sets $\G_\eta$ and $\square_\g^\eta$ by (\ref{2.18}),
 fix an arbitrary $\g\in\G_\eta$ and by (\ref{5.26}), (\ref{5.30}) we obtain:
\begin{equation}\label{5.31}
\begin{aligned}
& \bigg|\frac{1}{\eta^d \mes\square} \int\limits_{\square_\g^\eta} V^\e(x)\,dx - \frac{1}{\eta^d \mes\square} \int\limits_{\square_\g^\eta} V\bigg(\eta\g,\frac{x}{\e}\bigg)\,dx \bigg|\leqslant \rho_9(\kappa\eta),
\\
& \bigg|\frac{1}{\eta^d \mes\square} \int\limits_{\square_\g^\eta} V^0(x)\,dx - \frac{1}{\eta^d \mes\square} \int\limits_{\square_\g^\eta} V^0(\eta\g)\,dx \bigg|\leqslant \rho_9(\kappa\eta),
\end{aligned}
\end{equation}
where $\kappa:=\max\limits_{y,z\in\overline{\square}}|y-z|$.
In view of  (\ref{5.28}) we find:
\begin{align*}
\Bigg| \frac{1}{\eta^d \mes\square} \int\limits_{\square_\g^\eta} V\left(\eta\g,\frac{x}{\e}\right)\,dx - V^0(\eta\g)\Bigg| = \Bigg| \frac{1}{\big(\frac{\eta}{\e}\big)^d\mes\square} \int\limits_{\square_\g^{\frac{\eta}{\e}}}
 V(\eta\g,\xi) \,d\xi-V^0(\eta\g)\Bigg|\leqslant \rho_{10}\left(\frac{\eta}{\e}\right).
\end{align*}
This estimate and (\ref{5.31}) imply:
\begin{equation*}
\Bigg| \frac{1}{\eta^d \mes\square} \int\limits_{\square_\g^\eta} V^\e(x)\,dx - \frac{1}{\eta^d \mes\square} \int\limits_{\square_\g^\eta}V^0(x)\,dx \Bigg|\leqslant 2\rho_9\big(\kappa\e^\frac{1}{2}\big) + \rho_{10}\big(\e^{-\frac{1}{2}}\big)\to0,\quad \e\to0.
\end{equation*}
Hence, the assumptions of Theorem~\ref{th2} are satisfied and
\begin{equation*}
\|V^\e-V^0\|_{\fM_{1,-1}} \leqslant C\big(\rho_9\big(\kappa\e^\frac{1}{2}\big) + \rho_{10}\big(\e^{-\frac{1}{2}}\big)+\e^{\frac{1}{2}}\big),
\end{equation*}
where $C$ is a constant independent of $\e$.

\subsection{Oscillation with modulated periodicity}

Here we consider one more case of a locally periodically fast oscillating coefficient similar to Subsection~\ref{sec:locper} but instead of using the variables $\frac{x}{\e}$ like in (\ref{5.23}), we replace them by $\frac{\phi(x)}{\e}$, where $\phi=\phi(x)$ is some mapping of the domain $\Om$ into the space $\mathds{R}^d$. Namely, let $\e$ be a one-dimensional small parameter, $L_0$ be a fixed lattice with a periodicity cell $\square_0$ and $V=V(x,\xi)$, $(x,\xi)\in\overline{\Om}\times\mathds{R}^d$, be an arbitrary $n\times n$ matrix function, uniformly bounded, $\square_0$-periodic in $\xi$.
We again assume that condition   (\ref{5.15}) holds true; since $m=1$, in this condition the supremum over $\xi'$ and $\tilde{\xi}'$ is in fact absent. By $\phi=\phi(x)=(\phi_1(x),\ldots,\phi_d(x))$ we define a vector-valued function on $\overline{\Om}$ with values in $\mathds{R}^d$; for simplicity we assume that $\phi\in C^2(\overline{\Om})$. We introduce the perturbation as
\begin{equation}\label{5.34}
V^\e(x):=V\left(x,\frac{\phi(x)}{\e}\right).
\end{equation}

We consider two main cases. In the first case we suppose that $\phi$ is diffeomorphism of $\overline{\Om}$ onto the closure of some domain $\tilde{\Om}\subset \mathds{R}^d$. The corresponding  Jacobi matrix
\begin{equation}\label{5.42}
\mathrm{J}(x):=
\begin{pmatrix}\displaystyle
\frac{\p \phi_1}{\p x_1}(x) & \ldots & \displaystyle\frac{\p \phi_1}{\p x_d}(x)
\\
\vdots & & \vdots
\\
\displaystyle\frac{\p \phi_d}{\p x_1}(x) & \ldots & \displaystyle\frac{\p \phi_d}{\p x_d}(x)
\end{pmatrix}
\end{equation}
then is invertible and we assume that the inverse matrix is uniformly bounded together with its  first derivatives on $\overline{\Om}$. We denote the entries of the inverse matrix $\mathrm{J}^{-1}$ by $\Phi_{ij}$:
\begin{equation*}
\mathrm{J}^{-1}(x)=
\begin{pmatrix}
\Phi_{11}(x) & \ldots & \Phi_{1d}(x)
\\
\vdots & & \vdots
\\
\Phi_{d1}(x) & \ldots & \Phi_{dd}(x)
\end{pmatrix}.
\end{equation*}
%
%vector function
%\begin{equation}\label{5.46}
%\begin{pmatrix}
%\Phi_1(x)
%\\
%\vdots
%\\
%\Phi_d(x)
%\end{pmatrix} = \mathrm{J}^{-1}(x)
%\begin{pmatrix}
%1
%\\
%\vdots
%\\
%1
%\end{pmatrix}
%\end{equation}
%is uniformly bounded together with its  first derivatives on $\overline{\Om}$.
A typical example of such vector-valued function $\phi(x)$ is a growing at infinity function, for instance,
$\phi_i(x)=x_i^{2\a_i+1}$ for $|x_i|>1$, where $\a_i$ are some natural numbers. In this case the quotient $\frac{\phi}{\e}$ in (\ref{5.34}) leads to increasing oscillations as $x$ grows, that is, the variable $\frac{\phi}{\e}$ grows simultaneously with $x$ and $\e$. A one-dimensional example demonstrating such situation is  $V(x,\xi)=\sin\xi$, $\phi(x)=x^3$, $\Om=[1,+\infty)$, and then
\begin{equation*}
V^\e(x)=\sin \frac{x^3}{\e}.
\end{equation*}
A sketched graph of the function $V^\e$ as $d=2$, $V(x,\xi)=\sin^2 \xi_1 \cos \xi_2$, $\phi(x)=(x_1^3,x_2^5)$ is given on Figure~\ref{fig:mod-gr}.

In the second case we consider a periodic $\phi$. Namely, let $\phi$ be periodic with respect to some fixed lattice $L_1$ with a periodicity cell $\square_1$.  Here a model one-dimensional example is
\begin{equation*}
V^\e(x)=\sin \frac{\cos x}{\e}.
\end{equation*}
A sketched graph of the function $V^\e$ as $d=2$, $V(x,\xi)=\sin^2 \xi_1 \cos \xi_2 $, $\phi(x)=(\sin x_1,\cos^2 x_2)$ is given on Figure~\ref{fig:mod-per}.

We introduce the corresponding Jacobi matrix by formula (\ref{5.42}) and assume that the determinant $\det \mathrm{J}(x)$ vanishes on some set $S$ and for all sufficiently small $\eta$ and all $\g\in \mathds{Z}^d$ the estimate holds:
\begin{equation}\label{5.36}
\mes\big\{x\in\mathds{R}^d:\ x\in \eta\g+\eta(0,1)^d,\ \dist(x,S)<\eta^2\big\} \leqslant C\eta^{d+1},
\end{equation}
where $C$ is some fixed constant independent of $\g$, $\eta$.
This condition is very natural and means the following. Consider a layer of width $2\eta^2$ along $S$ and consider its piece contained in the cube $\eta\g+\eta (0,1)^d$ for an arbitrary $\g\in\mathds{Z}^d$.  In the general situation, $S$ is some surface in $\Om$ of codimension one and if it is sufficiently smooth and regular, then condition (\ref{5.36}) is obviously satisfied.

%
%In the considered case the quotient $\frac{\phi(x)}{\e}$ is $\square_1$-periodic with respect to $x$ % %and the fast oscillations of the perturbation defined by (\ref{5.34}) are produced by  fast %oscillations of $V^\e$.

We are going to show that in both  cases condition (\ref{2.13}) is satisfied with $\G=\mathds{Z}^d$, $\square:=(0,1)^d$ and
\begin{equation*}%\%label{5.35}
V^0(x):=\frac{1}{\mes \square_0} \int\limits_{\square_0} V(x,\xi)\,d\xi
\end{equation*}
under an appropriate choice of $\eta$. We observe also that the made assumptions ensure that $V^\e $ and $V^0$  are uniformly bounded.

\begin{figure}[t]
\begin{center}
\includegraphics[scale=0.5]{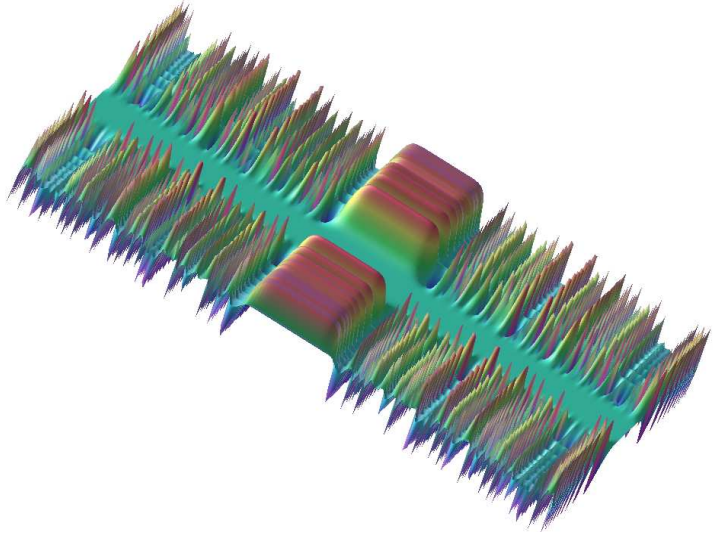}
\caption{Coefficient with modulated periodicity: growing $\phi$}\label{fig:mod-gr}
\end{center}
\end{figure}

\begin{figure}[t]
\begin{center}
\includegraphics[scale=0.5]{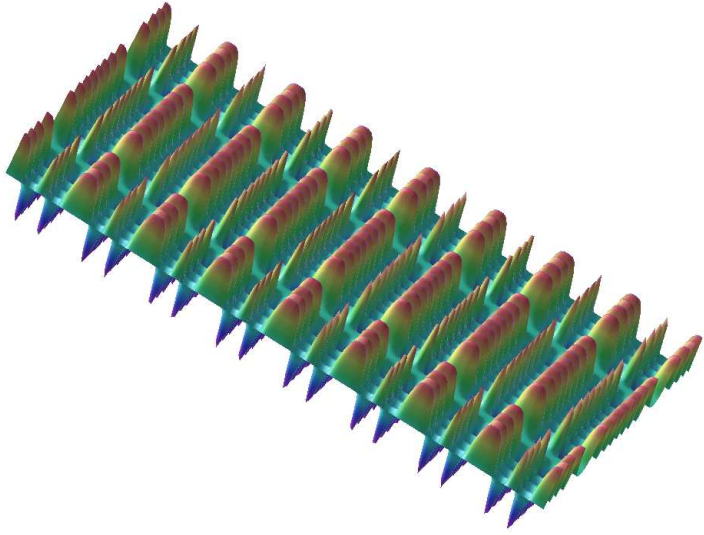}
\caption{Coefficient with modulated periodicity: periodic $\phi$}\label{fig:mod-per}
\end{center}
\end{figure}

We begin with the first case. We let $\eta:=\e^\frac{1}{2}$,  introduce the sets $\G_\eta$ and $\square_\g^\eta$ by (\ref{2.18}), choose an arbitrary $\g\in\G_\eta$ and by condition~(\ref{5.15}) we obtain:
\begin{equation}\label{5.40}
 \frac{1}{\eta^d}\Bigg| \int\limits_{\square_\g^\eta} \bigg( V^\e(x) - V\bigg(\eta\g,\frac{\phi(x)}{\e}\bigg)\bigg)\,dx
\Bigg| + \frac{1}{\eta^d}\Bigg| \int\limits_{\square_\g^\eta} \bigg( V^0(x) - V^0(\eta\g)\bigg)\,dx
\Bigg|\leqslant 2\rho_8\big(\sqrt{d}\e^\frac{1}{2}\big).
\end{equation}
Since
\begin{equation*}
\int\limits_{\square_0} \big(V(\eta\g,\xi)-
V^0(\eta\g)\big) \,d\xi=0,
\end{equation*}
there exists a solution to the equation
\begin{equation*}
\D_\xi V_1(\eta\g,\xi)=V(\eta\g,\xi)-
V^0(\eta\g)\quad\text{in}\quad\square_0
\end{equation*}
subject to the periodic boundary conditions on $\p\square_0$. Applying then the results from \cite[Ch. I\!I\!I, Sect. 15]{Ld}, in view of the assumed smoothness of $V$ we see that $V_1$ is differentiable in $\xi$ and
\begin{equation}\label{5.41}
|\nabla_\xi V_1(\eta\g,\xi)|\leqslant C
\end{equation}
uniformly in $\eta\g$ and $\xi$. Then, letting $\xi=\frac{\phi(x)}{\e}$, it is straightforward to confirm that
\begin{equation}\label{5.43}
V\bigg(\eta\g,\frac{\phi(x)}{\e}\bigg)- V^0(\eta\g) = \sum\limits_{j=1}^{d} \frac{\p^2 V_1}{\p\xi_j^2}\bigg(\eta\g,\frac{\phi(x)}{\e}\bigg)=\e\sum\limits_{i,j=1}^{d} \Phi_{ij}(x) \frac{\p\ }{\p x_i} \frac{\p V_1}{\p\xi_j}\bigg(\eta\g,\frac{\phi(x)}{\e}\bigg)
\end{equation}
and this allows us to integrate by parts as follows:
\begin{equation}\label{5.44}
\begin{aligned}
\frac{1}{\eta^d} \int\limits_{\square_\g^\eta} \bigg(V\bigg(\eta\g,\frac{\phi(x)}{\e}\bigg)- V^0(\eta\g)\bigg)\,dx= &
\frac{\e}{\eta^d} \int\limits_{\p \square_\g^\eta} \sum\limits_{i,j=1}^{d} \Phi_{ij}(x)  \frac{\p V_1}{\p\xi_j}\bigg(\eta\g,\frac{\phi(x)}{\e}\bigg)\nu_i(x)\,ds
\\
& -
\frac{\e}{\eta^d} \int\limits_{\square_\g^\eta} \sum\limits_{i,j=1}^{d} \frac{\p\Phi_{ij}}{\p x_i}(x) \frac{\p V_1}{\p\xi_j}\bigg(\eta\g,\frac{\phi(x)}{\e}\bigg)\,dx,
\end{aligned}
\end{equation}
where $\nu=(\nu_1,\ldots,\nu_d)$ is the outward normal to $\p \square_\g^\eta$. Applying then estimate (\ref{5.41}) and using the assumed uniform boundedness of $\Phi_{ij}$ and their first derivatives, we get:
\begin{align*}
\frac{1}{\eta^d} \Bigg|\int\limits_{\square_\g^\eta} \bigg(V\bigg(\eta\g,\frac{\phi(x)}{\e}\bigg)- V^0(\eta\g)\bigg)\,dx\Bigg|
\leqslant C\bigg(\frac{\e}{\eta}+\e\bigg)\leqslant C \e^\frac{1}{2},
\end{align*}
where $C$ are some constants independent of $\e$. In view of (\ref{5.40}) this proves (\ref{2.13}) with $$\rho_1(\e)=C\big(\e^\frac{1}{2}+\rho_8\big(\sqrt{d}\e^\frac{1}{2}\big)\big)$$ and hence,
\begin{equation*}
\|V^\e-V^0\|_{\fM_{1,-1}}\leqslant C \big(\e^\frac{1}{2}+\rho_8\big(\sqrt{d}\e^\frac{1}{2}\big)\big).
\end{equation*}

We proceed to the second case, when $\phi$ is periodic. We denote
\begin{equation*}%\%label{5.33}
p_0(r):=\inf\limits_{\substack{x\in\overline{\Om}
\\
\dist(x,S)\geqslant r}} \det J(x),\qquad p_1(r):=\min\big\{r p_0(r^2),p_0^2(r^2)\big\},\qquad r>0.
\end{equation*}
It is clear that $p_0(r)$ is a continuous function positive for $r>0$ and
$p_0(r)\to+0$ as $r\to+0$. Hence, the same is true for $p_1(r)$. We then choose $\eta(\e)$ as
\begin{equation*}
\eta(\e):=\inf\big\{r>0:\, p_1(r)\geqslant \e^\frac{1}{2}\big\}.
\end{equation*}
In view of the aforementioned properties of the function $p_1$, the function $\eta(\e)$ is positive for $\e>0$ and
\begin{equation}\label{5.48}
\eta(\e)\to+0,\quad \e\to+0, \qquad p_1(\eta(\e))\geqslant \e^\frac{1}{2}.
\end{equation}
%
% We then choose $\eta=\eta(\e)$ so that
%%\begin{equation*}%\%label{5.48}
%$\eta^d(\e) p^2\big(\eta^2(\e)\big)\geqslant \e^\frac{1}{2}$.
%%\end{equation*}

We again use estimate (\ref{5.40}) and representation (\ref{5.43}), which is now valid only outside the set $S$. In order to estimate the contribution of the vicinity of $S$, we observe that
by condition (\ref{5.36}) %, (\ref{5.15})
and the uniform boundedness of $V^\e$ and $V^0$ the inequality holds:
\begin{equation}\label{5.39}
\frac{1}{\eta^d} \Bigg| \int\limits_{S_\g} \bigg( V^\e(x) - V\bigg(\eta\g,\frac{\phi(x)}{\e}\bigg)\bigg)\,dx
\Bigg| \leqslant C \eta,\qquad S_\g:=\big\{x\in \square_\g^\eta:\ \dist(x,S)\leqslant \eta^2\big\},
\end{equation}
where $C$ is some constant independent of $\e$,  $\eta$ and $\g$. As in (\ref{5.44}), we obtain:
\begin{equation}\label{5.45}
\begin{aligned}
\frac{1}{\eta^d} \int\limits_{\square_\g^\eta \setminus S_\g} \bigg(V\bigg(\eta\g,\frac{\phi(x)}{\e}\bigg)- V^0(\eta\g)\bigg)\,dx= &
\frac{\e}{\eta^d} \int\limits_{\p (\square_\g^\eta\setminus S_\g)} \sum\limits_{i,j=1}^{d} \Phi_{ij}(x)  \frac{\p V_1}{\p\xi_j}\bigg(\eta\g,\frac{\phi(x)}{\e}\bigg)\nu_i(x)\,ds
\\
& -
\frac{\e}{\eta^d} \int\limits_{\square_\g^\eta\setminus S_\g} \sum\limits_{i,j=1}^{d} \frac{\p\Phi_{ij}}{\p x_i}(x) \frac{\p V_1}{\p\xi_j}\bigg(\eta\g,\frac{\phi(x)}{\e}\bigg)\,dx.
\end{aligned}
\end{equation}
Since the functions $\phi_j$ are $\square_1$-periodic and smooth, it follows that the entries of the inverse matrix $\mathrm{J}^{-1}$ satisfy the estimates
\begin{equation*}
|\Phi_{ij}(x)|\leqslant \frac{C}{p_0(\eta^2)},\qquad |\nabla_x \Phi_{ij}(x)|\leqslant \frac{C}{p_0^2(\eta^2)},\qquad x\in \square_\eta^\g\setminus S_\g,
\end{equation*}
where $C$ are some fixed constants independent of $x$,  $\eta$ and $\g$. These estimates and (\ref{5.48}), (\ref{5.45}) give:
\begin{equation*}
\frac{1}{\eta^d}\Bigg|\int\limits_{\square_\eta^\g\setminus S_\g} \bigg(V\bigg(\eta\g,\frac{\phi(x)}{\e}\bigg)- V^0(\eta\g)\bigg)\,dx\Bigg| \leqslant C \left(\frac{\e}{\eta p_0(\eta^2)} + \frac{\e}{p_0^2(\eta^2)}\right)\leqslant C\e^\frac{1}{2},
\end{equation*}
where $C$ are some constants independent of $\e$, $\eta$ and $\g$. Hence, by (\ref{5.39}), (\ref{5.40}),
\begin{equation*}%\l%abel{5.49}
\frac{1}{\eta^d}\Bigg| \int\limits_{\square_\eta^\g} \big( V^\e(x) - V^0(x))\big)\,dx\Bigg| \leqslant C\Big(\e^\frac{1}{2} + \eta + \rho_8\big(\sqrt{d}\eta\big) \Big).
\end{equation*}
Therefore,  by Theorem~\ref{th2}, the matrix function $V^\e$ converges to $V^0$ in $\fM_{1,-1}$ and
\begin{equation*}%\%label{5.50}
\|V^\e - V^0\|_{\fM_{1,-1}}\leqslant C \Big(\e^\frac{1}{2} + \eta + \rho_8\big(\sqrt{d}\eta\big) \Big),
\end{equation*}
where $C$ is some constant independent of $\e$ and $\eta$.

\subsection{Oscillations with fractal periodicity}

Here we consider fast oscillating coefficients of a fractal structure suggested in \cite{NS}. We make the  assumptions similar to the previous example: $\e$ is a one-dimensional small parameter and $V=V(x,\xi)$, $\xi=(\xi_1,\ldots,\xi_d)\in\mathds{R}^d$,  $(x,\xi)\in\overline{\Om}\times\mathds{R}^d$, is an arbitrary $n\times n$ matrix function, uniformly bounded, $a_i$-periodic in each variable $\xi_i$, $i=1,\ldots,d$;
condition (\ref{5.15}) is supposed to be satisfied.

The oscillating coefficient is defined as follows:
\begin{equation*}%\l%abel{5.52}
V^\e(x):=V\left(x,\frac{x_1}{\e},\frac{x_1 x_2}{\e^2},\frac{x_1 x_2 x_3}{\e^3},\ldots,\frac{x_1\cdots x_d}{\e^d}\right).
\end{equation*}
The main feature of such coefficients is that it is periodic in $\frac{x_1}{\e}$ with a fixed period, while the periodicity in $\frac{x_2}{\e}$ is with a period proportional to $x_1^{-1}$ and as $x_1$ increases, this period becomes smaller. A graph of the function $V^\e$ in the case $d=2$,  $V(x,\xi)=\sin \xi_1 \sin \xi_2$ is sketched on Figure~\ref{fig:fractal}.

We are going to show that the family $V^\e$ converges in $\fM_{1,-1}$ to
\begin{equation*}%\%label{5.53}
V^0(x):=\frac{1}{\mes\square_0} \int\limits_{\square_0} V(x,\xi)\,d\xi,
\qquad \square_0:=(0,a_1)\times(0,a_2)\times\ldots\times(0,a_d).
\end{equation*}

We let $\eta:=\e^\frac{1}{2}$, $\G:=2\mathds{Z}^d+(-1,\ldots,-1)$, $\square:=(0,2)^d$ and introduce the sets $\G_\eta$ and $\square_\g^\eta$ by (\ref{2.18}) and let us confirm that assumptions of Theorem~\ref{th2} are satisfied.  For an arbitrary $\g\in\G_\eta$ we obtain an inequality similar to (\ref{5.40}):
\begin{equation}\label{5.54}
 \frac{1}{\eta^d}\Bigg| \int\limits_{\square_\g^\eta} \bigg( V^\e(x) - V\bigg(\eta\g,\frac{x_1}{\e},\ldots,\frac{x_1\cdots x_d}{\e^d}\bigg)\bigg)\,dx
\Bigg| + \frac{1}{\eta^d}\Bigg| \int\limits_{\square_\g^\eta} \bigg( V^0(x) - V^0(\eta\g)\bigg)\,dx
\Bigg|\leqslant 2\rho_8\big(2\sqrt{d}\eta\big).
\end{equation}
Letting $\mu:=\eta/\e$, we have obvious identities
\begin{gather}\label{5.75}
\frac{1}{\eta^d} \int\limits_{\square_\g^\eta} \bigg( V\bigg(\eta\g,\frac{x_1}{\e},\ldots,\frac{x_1\cdots x_d}{\e^d}\bigg) - V^0(\eta\g)\bigg)\,dx
=\frac{1}{\mu^d} \int\limits_{\square_\g^\mu} \big( V(\eta\g,\z(\xi)) - V^0(\eta\g)\big)\,d\xi,
\\
\g=(\g_1,\ldots,\g_d), %\square_\g^\mu:=(\mu\g_1,\mu\g_1+2\mu)\times\cdots\times(\mu\g_{d-1},\mu\g_{d-1}+2\mu),\nonumber
\qquad
 \z(\xi):=(\z_1(\xi),\ldots,\z_d(\xi)),\qquad \z_i(\xi):=\xi_1\cdots\xi_i.
\nonumber
\end{gather}
We define
\begin{equation*}
V_d(x,\z'):=\frac{1}{a_d} \int\limits_{0}^{a_d} V(x,\z',s)\,ds,\qquad \z':=(\z_1,\ldots,\z_{d-1}),
\end{equation*}
and we immediately see that
\begin{equation}\label{5.56}
\Bigg|\int\limits_a^{a+t} \big( V(x, \z',s) - V_d(x,\z')\big)\,d s\Bigg|= \Bigg|\int\limits_{a+a_d\lceil\frac{t}{a_d}\rceil}
^{a+t} \big( V(x, \z',s) - V_d(x,\z')\big)\,ds\Bigg|
\leqslant C \left\{\frac{t}{a_d}\right\},
\end{equation}
for all $a\in\mathds{R}$, $x\in\Om$, $t>0$, $\z'\in\mathds{R}^{d-1}$, where $\lceil\,\cdot\,\rceil$ is the ceiling function, $\{\,\cdot\,\}$ is the fractional part of a number and $C$ is a constant independent of $a$, $t$, $\z'$ and $x$. We denote  $\tht:=\bigcup\limits_{i=1}^{d-1}\{\xi':\ |\xi_i|\leqslant \d\}$ with some $\d\in(0,\mu)$ and observe that as $\xi'\notin \tht$, we have  $\z_{d-1}(\xi)\ne0$ and
\begin{align*}
\Bigg|\int\limits_{\mu\g_d}^{\mu\g_d+2\mu}  \big( V(\eta\g,\z(\xi)) - V_d(\eta\g,\z'(\xi))\big)\,d\xi_d\Bigg|
&= \frac{1}{|\z_{d-1}(\xi)|}
\Bigg|\int\limits_{ \z_{d-1}(\xi) \mu\g_d}^{ \z_{d-1}(\xi) (\mu\g_d+2\mu)} \big( V(\eta\g,\z'(\xi),s) -  V_d(\eta\g,\z'(\xi))\big)\,ds\Bigg|
\\
&\leqslant C \frac{\left\{\frac{2\mu |\z_{d-1}(\xi)|}{a_d}\right\}}{|\z_{d-1}(\xi)|}
\leqslant \frac{C}{|\z_{d-1}(\xi)|},
\end{align*}
where $C$ is a constant independent of $\mu$, $\xi$ and $\g$. Hence,
\begin{equation}\label{5.76}
\Bigg|\int\limits_{\square_\g^\mu\setminus\{\xi:\,\xi'\in\tht\}}
 \big( V(\eta\g,\z(\xi)) - V_d(\eta\g,\z'(\xi)\big)\,d\xi\Bigg|\leqslant C \int\limits_{\tilde{\square}_\g^\mu\setminus\tht}
 \frac{d\xi'}{|\z_{d-1}(\xi)|}
 \leqslant C \bigg|\ln\frac{\mu}{\d}\bigg|^{d-1},
\end{equation}
where $\tilde{\square}_\g^\mu:=(\mu\g_1,\mu\g_1+2\mu)\times\ldots\times(\mu\g_{d-1},\mu\g_{d-1}+2\mu)$. Due to the uniform boundedness of $V$ and $V_d$ we also have a simple estimate
\begin{equation*}
\Bigg|\int\limits_{\square_\g^\mu\cap\{\xi:\,\xi'\in\tht\}}
 \big( V(\eta\g,\z(\xi)) - V_d(\eta\g,\z'(\xi)\big)\,d\xi\Bigg|\leqslant C \mes \square_\g^\mu\cap\{\xi:\,\xi'\in\tht\} \leqslant C
 \mu^{d-1}\d.
\end{equation*}

\begin{figure}[t]
\begin{center}
\includegraphics[scale=0.5]{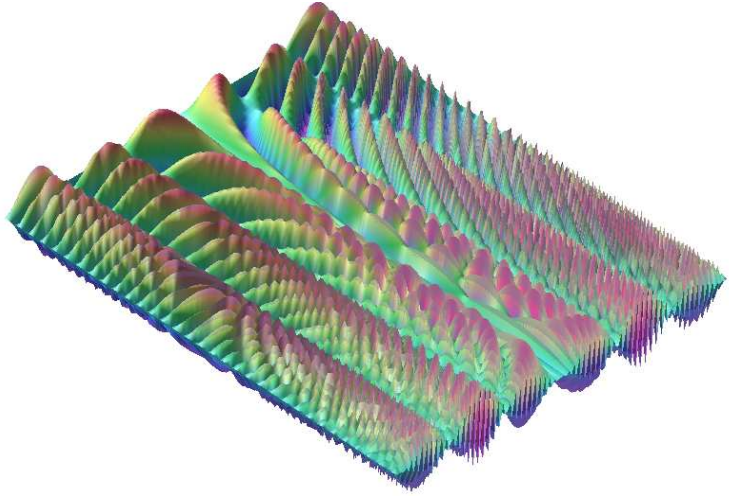}
\caption{Coefficient with fractal periodicity}\label{fig:fractal}
\end{center}
\end{figure}

\noindent
Choosing then $\d:=\mu^{-d+1}|\ln\mu|^{d-1}$, by the above estimate and (\ref{5.76}) we obtain:
\begin{equation}\label{5.77}
\begin{aligned}
\frac{1}{\mu^d}\Bigg|\int\limits_{\square_\g^\mu}
 \big( V(\eta\g,\z(\xi)) - V_d(\eta\g,\z'(\xi))\big)\,d\xi\Bigg|\leqslant C |\ln\mu|^{d-1}\mu^{-d}.
\end{aligned}
\end{equation}
where $C$ is some constant independent of $\mu$, $\g$, $\eta$.
This estimate means that we can replace the function $V(\eta\g,\xi_1,\ldots,\xi_1\cdots \xi_d)$ in the right hand side of (\ref{5.75}) by $V_d(\eta\g,\z'(\xi))$ making an error at most $C |\ln\mu|^{d-1}\mu^{-d}$; then we are led to the need of estimating the integral
\begin{equation*}
\frac{1}{\mu^{d-1}} \int\limits_{\tilde{\square}_\g^\mu} \big( V_d(\eta\g,\z'(\xi)) - V^0(\eta\g)\big)\,d\xi'.
\end{equation*}
Introducing then an auxiliary function
\begin{equation*}
V_{d-1}(x,\z''):=\frac{1}{a_{d-1}} \int\limits_{0}^{a_{d-1}} V_d(x,\z'',s)\,ds,\qquad \z'':=(\z_1,\ldots,\z_{d-2}),
\end{equation*}
as above we see that  the function $V_d(\eta\g,\z'(\xi))$ can be replaced by $V_{d-1}(\eta\g,\z''(\xi))$ up to an error at most $C |\ln\mu|^{d-2}\mu^{-d+1}$.  Reproducing then recurrently the above calculations, we finally see that
\begin{equation*}
\frac{1}{\eta^d} \Bigg|\int\limits_{\square_\g^\eta} \bigg( V\bigg(\eta\g, \frac{x_1}{\e},\ldots,\frac{x_1\cdots x_d}{\e^d}\bigg)-V^0(\eta\g)
\bigg)\,dx\Bigg|\leqslant \frac{C}{\mu},
\end{equation*}
where $C$ is some constant independent of $\e$, $\mu$, $\eta$ and $\g$. This estimate and (\ref{5.54}) yield:
\begin{equation*}
 \frac{1}{\eta^d}\Bigg| \int\limits_{\square_\g^\eta} \big( V^\e(x) - V^0(x)\big)\,dx
\Bigg| \leqslant 2\rho_8\big(2\sqrt{d}\e^{\frac{1}{2}}\big)+C\e^\frac{1}{2}.
\end{equation*}
Therefore, by Theorem~\ref{th2}, the function $V^\e$ converges to $V^0$ in $\fM_{1,-1}$ and
\begin{equation*}
\|V^\e-V^0\|_{\fM_{1,-1}}\leqslant C\big(
\rho_8\big(2\sqrt{d}\e^{\frac{1}{2}}\big)+\e^\frac{1}{2}
\big),
\end{equation*}
where $C$ is some constant independent of $\e$.

\subsection{Random perturbation}\label{StochHom}

Here we consider an example of random perturbation. We follow a general scheme, see, for instance, \cite[Ch. 7]{ZKO}. Namely, let $(\Ups,\cF,\vs)$ be a probability space with a sigma-algebra $\cF$ and a probability measure $\vs$.  By $\Tht=\Tht(x)$, $x\in\mathds{R}^d$, $\Tht:\, \Ups\to\Ups$ we denote a transformation obeying the standard properties:

1) group property: $\Tht(0)=\cI$, $\Tht(x+y)=\Tht(x)\Tht(y)$, $x,y\in\mathds{R}^d$.

2) measure preservation: the transformation $\Tht$ preserves the measure $\vs$ on $\Ups$, that is, for each $x\in\mathds{R}^d$ and each $\vs$-measurable set $F\in\cF$ the set $\Tht(x)F$ is $\vs$-measurable and $\vs(\Tht(x) F)=\vs(F)$;

3) measurability: for each measurable function $f=f(\varpi)$ on $\Ups$, the function $f(\Tht(x)\varpi)$ defined  on the Cartesian product $\Ups\times\mathds{R}^d$ is also measurable; here $\mathds{R}^d$ is endowed with the Lebesgue measure;

4) ergodicity: if a measurable function $f$ on $\Ups$ satisfies the identity $f(\Tht(x)\varpi)=f(\varpi)$ almost everywhere in $\Ups$ for each $x\in\mathds{R}^d$, then $f$ is constant almost everywhere on $\Ups$.

According to the Birkhoff Ergodic Theorem, if $f\in L_p(\Upsilon)$, $p\geqslant 1$, then for almost all $\varpi\in\Ups$ the function $f(\Tht(x)\varpi)$ possesses a mean value, namely,
\begin{equation*}%\%label{5.59}
\lim\limits_{t\to+\infty} \frac{1}{t^d\mes S} \int\limits_{t S} f(\Tht(x)\varpi)\,dx = \int\limits_{\Ups} f(y)\, \vs(dy)=:\mathbb{E}(f)\quad\text{for almost all}\quad \varpi\in\Ups,
\end{equation*}
where the limit exists for each bounded subset $S\subset \mathds{R}^d$ of a positive Lebesgue measure.

We proceed to defining a random perturbation. Let $V=V(x,\varpi)$ be a matrix function such that $V\in L_\infty(\Om\times\Ups;\mathds{C}^n)$. In addition, we assume that
\begin{equation}\label{5.60}
\sup\limits_{\substack{x,\tilde{x}\in
\Om
\\
|x-\tilde{x}|\leqslant \d} } \|V(x,\,\cdot\,)-V(\tilde{x},\,\cdot\,)\|_{L_\infty(\Ups;\mathds{C}^n)} \leqslant \rho_{10}(\d),\qquad \rho_{10}(\d)\to+0,\quad \d\to+0.
\end{equation}

Given a transformation $\Tht$ obeying the above properties, we define the random perturbation as
\begin{equation}\label{5.61}
V^\e(x,\varpi):=V\left(x,\Tht(x_\e)\varpi\right), \qquad x\in\Om,\quad \varpi\in\Ups,
\end{equation}
where $\e=(\e_1,\ldots,\e_d)\to0$ and for each $x\in\mathds{R}^d$ we have denoted $x_\e:=\big(\frac{x_1}{\e_1},\ldots,\frac{x_d}{\e_d}\big)$.
 The assumed properties of $V$ then easily imply that the introduced family $V^\e$ is bounded uniformly in $\e$. For each $x\in\Om$ we define
\begin{equation}\label{5.62}
V^0(x):=\mathbb{E}(V(x,\,\cdot\,))=\int\limits_{\Ups} V(x,y)\,\vs(dy).
\end{equation}
According to the Birkhoff Ergodic Theorem, for each $x\in \Om$
\begin{equation*}%\%label{5.63}
V^0(x)=\lim\limits_{t\to+\infty} \frac{1}{t^d} \int\limits_{(0,t)^d} V(x,\Tht(y)\varpi)\,dy \quad\text{for almost all}\quad \varpi\in\Ups.
\end{equation*}
Let $R\subseteq\Ups$ be some $\vs$-measurable set 
 and suppose that
the limit in the above identity converges uniformly in the following sense: there exists a function $\eta=\eta(\e)$ such that
\begin{equation}\label{5.74}
\eta(\e)\to0,\qquad  \frac{\e_{\max}}{\eta^{1+\a}(\e)}\to 0,\qquad \e\to0,\qquad \e_{\max}:=\max\limits_{j=1,\ldots,d}\{\e_j\},
\end{equation}
with some fixed $\a>0$ and the convergence holds:
\begin{equation}\label{5.64}
\esssup\limits_{\varpi\in R}
\sup\limits_{x,z\in\Om
}\Bigg|V^0(x)- \eta^{\a d}(\e) \int\limits_{\big(0,\eta^{-\a}(\e)\big)^d} V(x,\Tht(y)\Tht(z_\e)\varpi)\,dy\Bigg|\leqslant \rho_{11}(\e),\qquad \rho_{11}(\e)\to+0,\quad \e\to+0.
\end{equation}
As an example, when the above convergence is satisfied,  we can suppose that $R=\Ups$ and
\begin{equation*}
\esssup\limits_{\varpi\in\Ups}\sup\limits_{x\in\Om
}\Bigg|V^0(x)- \eta^{\a d}(\e) \int\limits_{\big(0,\eta^{-\a}(\e)\big)^d} V(x,\Tht(y)\varpi)\,dy\Bigg|\leqslant \rho_{11}(\e),\qquad \rho_{11}(\e)\to+0,\quad \e\to+0.
\end{equation*}
It is clear that the above convergence ensures (\ref{5.64}); at the same time, the above uniform convergence is a rather rare situation. Nevertheless, under certain assumptions it can be ensured, for instance, by applying the technique from the proof of Theorem~2.6 in \cite[Th. 2.6, Ch. I, Sect. 1.2]{Kr}. Namely, it is sufficient to assume that $\Ups$ is a compact metric space, each ball in this space has a non-zero measure, and the matrix functions $V(x,\Tht(y)\varpi)$ are uniformly continuous in $(y,\varpi)\in\mathds{R}^d\times\Ups$ uniformly in $x$, that is,
\begin{equation*}
\sup\limits_{\substack{x\in\Om,
\,
\varpi_1,\varpi_2\in\Ups
\\
\dist(\varpi_1,\varpi_2)\leqslant \d}} |V(x,\varpi_1)-V(x,\varpi_2)|\leqslant \rho_{12}(\d),
\qquad \sup\limits_{\substack{x\in\mathds{R}^d,\, \varpi_1,\varpi_2\in\Ups
\\ \dist(\varpi_1,\varpi_2)\leqslant \d}}
\dist(\Tht(x)\varpi_1,\Tht(x)\varpi_2)\leqslant\rho_{13}(\d),
\end{equation*}
where $\rho_{12}(\d)\to0$, $\rho_{13}(\d)\to0$ as $\d\to0$.

We are going to show that under the above assumptions, the family $V^\e$ satisfies condition (\ref{2.13}) for each $\varpi\in R$. We let $\G:=\mathds{Z}^d$, $\square:=(0,1)^d$ and choose the function $\eta$ introduced in (\ref{5.74}), (\ref{5.64}).
% $\eta(\e):=\e_0^{\frac{1}{3}}$, $\e_0:=\max\{\e_1,\ldots,\e_d\}$.
It follows immediately from (\ref{5.60}) that the function $V^0$ defined in (\ref{5.62}) is uniformly continuous in $\Om$:
\begin{equation}\label{5.66}
\sup\limits_{\substack{x,\tilde{x}\in
\Om
\\
|x-\tilde{x}|\leqslant \d} } |V^0(x)-V^0(\tilde{x})|\leqslant \rho_{10}(\d).
\end{equation}
We choose   arbitrary $\g\in\G_\eta$, $\varpi\in R$ and by (\ref{5.60}) we again get:
\begin{equation}\label{5.65}
\Bigg|\frac{1}{\eta^d}\int\limits_{\eta\g+\eta\square} \Big(V\big(x,\Tht(x_\e)\varpi\big) - V\big(\eta\g,\Tht(x_\e)\varpi\big)\Big)\,dx\Bigg| \leqslant \rho_{10}(\sqrt{d}\eta).
\end{equation}
We denote
%\begin{equation*}
$\Xi(a_1,\ldots,a_d):=(0,a_1)\times \ldots \times (0,a_d)$. 
By a simple change of variables and the group property we obtain:
\begin{equation}\label{5.67}
\begin{aligned}
\frac{1}{\eta^d} \int\limits_{\eta\g+\eta\square} &V\big(\eta\g,\Tht(x_\e)\varpi\big)\,dx= \frac{1}{\eta^d} \int\limits_{\eta\square} V\big(\eta\g,\Tht(x_\e)
\Tht(\g_\e)\varpi\big)\,dx
\\
=&\frac{\prod\limits_{j=1}^{d}\e_j}{\eta^d} \int\limits_{\Xi\big(\frac{\eta}{\e_1},\ldots,\frac{\eta}{\e_d}\big)}
V\big(\eta\g,\Tht(x) \Tht(\g_\e)\varpi\big)\,dx
\\
=& \frac{\prod\limits_{j=1}^{d}\e_j}{\eta^d} \sum\limits_{\substack{z=(z_1,\ldots,z_d)\in\mathds{Z}_+^d
\\
0\leqslant z_j\leqslant \big[\frac{\eta^{1+\a}}{\e_j}\big]-1}} \;
\int\limits_{\eta^{-\a}(z+\square)}
V\big(\eta\g,\Tht(x) \Tht(\g_\e)\varpi\big)\,dx
\\
&+ \frac{\prod\limits_{j=1}^{d}\e_j}{\eta^d} \int\limits_{\Xi\big(\frac{\eta}{\e_1},\ldots,\frac{\eta}{\e_d}\big) \setminus
\Xi\big(\eta^{-\a}\big[\frac{\eta^{1+\a}}{\e_1}\big],
\ldots,\eta^{-\a}\big[\frac{\eta^{1+\a}}{\e_d}\big]\big)} V\big(\eta\g,\Tht(x) \Tht(\g_\e)\varpi\big)\,dx
\\
=& \prod\limits_{j=1}^{d}\frac{\e_j}{\eta^{1+\a}}  \sum\limits_{\substack{z=(z_1,\ldots,z_d)\in\mathds{Z}_+^d
\\
0\leqslant z_j\leqslant \big[\frac{\eta^{1+\a}}{\e_j}\big]-1}} \eta^{\a d} \int\limits_{\eta^{-\a} \square }
V\big(\eta\g,\Tht(x) \Tht(\eta^{-\a}z+\g_\e)\varpi\big)\,dx
\\
&+ \frac{\prod\limits_{j=1}^{d}\e_j}{\eta^d} \int\limits_{\Xi\big(\frac{\eta}{\e_1},\ldots,\frac{\eta}{\e_d}\big) \setminus
\Xi\big(\eta^{-\a}\big[\frac{\eta^{1+\a}}{\e_1}\big],
\ldots,\eta^{-\a}\big[\frac{\eta^{1+\a}}{\e_d}\big]\big)} V\big(\eta\g,\Tht(x) \Tht(\g_\e)\varpi\big)\,dx,
\end{aligned}
\end{equation}
where $[\,\cdot\,]$ denotes the integer part of a number.
It is clear that
\begin{align*}
\frac{\prod\limits_{j=1}^{d}\e_j}{\eta^d} \mes \Xi\left(\frac{\eta}{\e_1},\ldots,\frac{\eta}{\e_d}\right) \setminus &
\Xi\left(\eta^{-\a}\left[\frac{\eta^{1+\a}}{\e_1}\right],\ldots, \eta^{-\a}\left[\frac{\eta^{1+\a}}{\e_d}\right]\right)= 1-
\prod\limits_{j=1}^{d} \frac{\e_j}{\eta^{1+\a}}\left[\frac{\eta^{1+\a}}{\e_j}\right]
\\
&=1-\prod\limits_{j=1}^{d} \left(1-\frac{\e_j}{\eta^{1+\a}} \left\{\frac{\eta^{1+\a}}{\e_j}\right\} \right)\leqslant C \sum\limits_{j=1}^{d} \frac{\e_j}{\eta^{1+\a}} \leqslant C \frac{\e_{\max}}{\eta^{1+\a}},
\end{align*}
where $\{\,\cdot\,\}$ is the fractional part of a number and $C$ is a constant independent of $\e$. Hence, in view of the uniform boundedness of $V$,
\begin{equation}\label{5.69}
\begin{aligned}
&
\frac{\prod\limits_{j=1}^{d}\e_j}{\eta^d} \Bigg| \int\limits_{\Xi\big(\frac{\eta}{\e_1},\ldots,\frac{\eta}{\e_d}\big) \setminus
\Xi\big(\frac{1}{\eta}\big[\frac{\eta^{1+\a}}{\e_1}\big],
\ldots,\frac{1}{\eta}\big[\frac{\eta^{1+\a}}{\e_d}\big]\big)} V\big(\eta\g,\Tht(x) \Tht(\g_\e)\varpi\big)\,dx\Bigg|\leqslant C\frac{\e_{\max}}{\eta^{1+\a}},
\\
& \left|\prod\limits_{j=1}^{d}\frac{\e_j}{\eta^{1+\a}}-  \prod\limits_{j=1}^{d}\left[\frac{\eta^{1+\a}}{\e_j}\right]^{-1}\right| \Bigg|\sum\limits_{\substack{z=(z_1,\ldots,z_d)\in\mathds{Z}_+^d
\\
0\leqslant z_j\leqslant \big[\frac{\eta^{1+\a}}{\e_j}\big]-1}}  \eta^{\a d}
\int\limits_{\eta^{-\a} \square }
V\big(\eta\g,\Tht(x) \Tht(\eta^{-\a}z+\g_\e) \varpi\big)\,dx \Bigg|\leqslant C\frac{\e_{\max}}{\eta^{1+\a}},
\end{aligned}
\end{equation}
where $C$ is a constant independent of $\e$, $\g$ and $\varpi$. It follows from (\ref{5.64}) that
\begin{align*}
\Bigg|\prod\limits_{j=1}^{d}\left[\frac{\eta^{1+\a}}{\e_j}\right]^{-1}   & \sum\limits_{\substack{z=(z_1,\ldots,z_d)\in\mathds{Z}_+^d
\\
0\leqslant z_j\leqslant \big[\frac{\eta^{1+\a}}{\e_j}\big]-1}} \eta^{\a d}
\int\limits_{\eta^{-\a} \square }
V\big(\eta\g,\Tht(x) \Tht\big(\eta^{-\a}z+\g_\e\big) \varpi\big)\,dx - V^0(\eta\g)\Bigg|
\\
&= \Bigg|\prod\limits_{j=1}^{d} \left[\frac{\eta^{1+\a}}{\e_j}\right]^{-1} \sum\limits_{\substack{z=(z_1,\ldots,z_d)\in\mathds{Z}_+^d
\\
0\leqslant z_j\leqslant \big[\frac{\eta^{1+\a}}{\e_j}\big]-1}}  \eta^{\a d}\int\limits_{\eta^{-\a} \square }
\left(V\big(\eta\g,\Tht(x) \Tht(\eta^{-\a}z+\g_\e)\varpi\big)- V^0(\eta\g)\right)\,dx\Bigg|
\leqslant \rho_{11}(\e).
\end{align*}
This estimate and (\ref{5.69}), (\ref{5.67}), (\ref{5.65}), (\ref{5.66}) imply
\begin{equation*}
\Bigg|  \frac{1}{\eta^d}\int\limits_{\eta\g+\eta\square} \Big(V\big(x,\Tht(x_\e)\varpi\big)- V_0(x)\Big)\,dx \Bigg| \leqslant  C \left( \rho_{10}(\sqrt{d}\eta) + \frac{\e_{\max}}{\eta^{1+\a}} + \rho_{11}(\e)\right),
\end{equation*}
where $C$ is a constant independent of $\e$, $\eta$, $\g$ and $\varpi$.
This estimate ensures condition (\ref{2.13}) and hence, by Theorem~\ref{th2}, the convergence of $V^\e$ to $V^0$ in the norm of the space $\fM_{1,-1}$ for each $\varpi\in \Ups$.

Another way of studying the random perturbation concerns the convergence in the sense of the mathematical expectation. Namely, let $P_j^\e=Q_j^\e=0$ and $V^\e$ be given by (\ref{5.61}). For simplicity, assume that $V$ is independent of $x$, that is, $V=V(\varpi)$, and $V\in L_\infty(\Ups;\mathds{M}_n)$. Then in the considered case  the limiting potential $V^0$ defined by (\ref{5.62}) is constant. We define the space $L_2(\Om;\mathds{C}^n;\Ups)$ as consisting of the vector functions $u=u(x,\varpi)$ such that $u=u(\,\cdot\,,\varpi)$ is an element of $L_2(\Om;\mathds{C}^n)$ for almost each $\varpi\in\Ups$ with a finite norm
\begin{equation*}
\|u\|_{L_2(\Om;\mathds{C}^n;\Ups)}:=\bigg(\int\limits_{\Ups}
\|u(\,\cdot\,,\varpi)\|_{L_2(\Om;\mathds{C}^n)}^2\vs(d\varpi)
\bigg)^\frac{1}{2}.
\end{equation*}
In the same way we define the spaces $W_2^1(\Om;\mathds{C}^n;\Ups)$ and $\Ho_2^1(\Om;\mathds{C}^n;\Ups)$; a corresponding of these spaces serves as $\mathfrak{V}$ and is denoted by $\fV_\Ups$. The corresponding spaces $\mathfrak{V}_\Ups^\ast$  consist of the bounded antilinear functionals $f(\varpi)$,  $\varpi\in\Ups$, on $W_2^1(\Om;\mathds{C}^n;\Ups)$ or, respectively, on $\Ho_2^1(\Om;\mathds{C}^n;\Ups)$, such that the function $\la f(\varpi), u(\,\cdot\,,\varpi)\ra$ is $\vs$-integrable on $\Ups$. The action of such functional is defined as
\begin{equation*}
\la f,u\ra= \int\limits_{\Ups} \la f(\varpi), u(\,\cdot\,,\varpi)\ra \vs(d\varpi).
\end{equation*}
A corresponding analog of the space $\fM_{1,-1}$, denoted here by $\fM_{1,-1,\Ups}$, is introduced as the space of the matrix functions $V=V(x,\varpi)$ with a finite norm
\begin{equation*}
\|V\|_{\fM_{1,-1,\Ups}}:=\sup\limits_{u,v\in \mathfrak{V}} \frac{|(Vu,v)_{L_2(\Om;\mathds{C}^n;\Ups)}|} {\|u\|_{\mathfrak{V}}\|v\|_{\mathfrak{V}}}.
\end{equation*}
In this framework, the statement of Theorem~\ref{th1} remains true, just the operator norm should be replaced by $\|\,\cdot\,\|_{\fV_\Ups^\ast\to\fV_\Ups}$.
Then we can reproduce the proof of Theorem~\ref{th2}. Namely, in its proof we can  employ relations (\ref{5.71}) and then identity (\ref{5.70})
is to be rewritten as
\begin{equation}\label{5.73}
\begin{aligned}
((V^\e-V^0)u,v)_{L_2(\square_\g^\eta;\mathds{C}^n;\Ups)}= & \eta^d \mes \square  \int\limits_{\Ups} \big(V_\g(\varpi)\la u(\,\cdot\,,\varpi)\ra_\g, \la v(\,\cdot\,,\varpi)\ra_\g\big)_{\mathds{C}^n}\vs(d\varpi)
\\
&+ \big((V^\e-V^0)\la u(\,\cdot\,,\varpi)\ra_\g, v(\,\cdot\,,\varpi)\big)_{L_2(\square_\g^\eta;\mathds{C}^n;\Ups)}
\\
&+\big((V^\e-V^0)u_\g,
v(\,\cdot\,,\varpi)\big)_{L_2(\square_\g^\eta;\mathds{C}^n;\Ups)},
\end{aligned}
\end{equation}
where $\eta$ is chosen as above in this subsection. By the second and third estimates in (\ref{5.72}) we can get similar bounds for the second and the third terms in (\ref{5.73}). However, the integral in the first term is to be estimated then by $\eta^d\mes\square\|V_\g(\varpi)\la u(\,\cdot\,,\varpi)\ra_\g\|_{L_2(\Ups;\mathds{C}^n)}\|\la v(\,\cdot\,,\varpi)\ra_\g\|_{L_2(\Ups;\mathds{C}^n)}$ times some small factor independent of $\g$. In the particular case $u=v$ this is equivalent to estimating the integral $\int\limits_{\Ups} V_\g(\varpi) h(\varpi)\vs(d\varpi)$ by $\|h\|_{L_1(\Ups;\mathds{C}^n)}$ uniformly in $h$ and the best possible estimate here is by $\|V_\g\|_{L_\infty(\Ups;\mathbb{M}_n)} \|h\|_{L_1(\Ups;\mathds{C}^n)}$. Then we necessary need $\|V_\g\|_{L_\infty(\Ups;\mathbb{M}_n)}\to0$ as $\eta\to0$ and this is equivalent to the uniform convergence in $\varpi$ discussed above. Thus, we conclude that the norm resolvent convergence as in Theorem~\ref{th1} even in the case of the perturbation only by a random potential  and even in the sense of the expectation value
requires uniform convergence (\ref{5.64}) and fails otherwise. This   agrees with the results in \cite{Stoch2}, \cite{Stoch1} on stochastic homogenization, where all convergence rates were estimated either in the strong sense, that is, for each fixed right hand side, or by using certain weighted norm of the right hand sides.

\subsection{Generating new perturbations}

Here we discuss how to produce new families $P_j^\e$, $Q_j^\e$ and $V^\e$ converging in appropriate norms  once we are given
some families already converging in these spaces. In other words, we discuss possible operations on such families, which preserve the needed convergence.

Assume that we are given families $V^\e$ and $Q_j^\e$ converging to $V^0$ and $Q_j^0$ respectively in $\fM_{1,-1}$ and $\fM_{1,0}$. Let $\Psi=\Psi(x)$ be a $n\times n$ matrix function independent of $\e$ and belonging to $L_\infty(\Om;\mathds{C}^n)$. Then it follows immediately from definition  (\ref{2.11}) of the norms in $\fM_{1,0}$ that the family $\Psi Q_j^\e$ converge to $\Psi Q_j^0$ in  $\fM_{1,0}$ as $\e\to0$ and
\begin{align*}
\|\Psi Q_j^\e -\Psi Q_j^0\|_{\fM_{1,0}}= &
\sup\limits_{\substack{u\in \fV\\ u\ne0}}
\frac{\|\Psi(Q_j^\e-Q_j^0)u\|_{L_2(\Om;\mathds{C}^n)}}
{\|u\|_{\fV}}
\\
\leqslant & \|\Psi\|_{L_\infty(\Om;\mathds{M}_n)} \sup\limits_{\substack{u\in \fV\\ u\ne0}}
\frac{\|(Q_j^\e-Q_j^0)u\|_{L_2(\Om;\mathds{C}^n)}}
{\|u\|_{\fV}}=
 \|\Psi\|_{L_\infty(\Om;\mathds{M}_n)} \|Q_j^\e - Q_j^0\|_{\fM_{1,0}}.
\end{align*}
If $\Psi\in W_\infty^1(\Om;\mathds{C}^n)$,  a right multiplication is also possible:
\begin{align*}
\|Q_j^\e\Psi - Q_j^0\Psi\|_{\fM_{1,0}}= &
\sup\limits_{\substack{u\in \fV\\ u\ne0}}
\frac{\|(Q_j^\e-Q_j^0)\Psi u\|_{L_2(\Om;\mathds{C}^n)}}
{\|u\|_{\fV}}
\\
\leqslant &
\sup\limits_{\substack{\Psi u\in \fV
\\ \Psi u\ne0}}
\frac{\|(Q_j^\e-Q_j^0)\Psi u\|_{L_2(\Om;\mathds{C}^n)}}
{\|\Psi u\|_{\fV}} \sup\limits_{\substack{u\in
\fV\\ u\ne0}}
\frac{\|\Psi u\|_{\fV}}{\|u\|_{\fV}}
\\
\leqslant &
2\|\Psi\|_{W_\infty^1(\Om;\mathds{M}_n)}\|(Q_j^\e-Q_j^0)\|_{\fM_{1,0}}.
\end{align*}
A similar multiplication is also possible for the family $V^\e$ but here we again should assume that $\Psi$ is an element of $W_\infty^1(\Om;\mathds{M}_n)$. Then, as above,
\begin{equation*}
\|\Psi (V^\e - V^0)\|_{\fM_{1,-1}}\leqslant 2\|\Psi\|_{W_\infty^1(\Om;\mathds{M}_n)}\| (V^\e - V^0)\|_{\fM_{1,-1}},
\qquad \| (V^\e - V^0)\Psi\|_{\fM_{1,-1}}\leqslant 2\|\Psi\|_{W_\infty^1(\Om;\mathds{M}_n)}\| (V^\e - V^0)\|_{\fM_{1,-1}}.
\end{equation*}

One more possibility is to sum two families. Namely, if we are given two families $V_1^\e$ and $V_2^\e$ converging respectively to $V_1^0$ and $V_2^0$ in $\fM_{1,-1}$, then their sum also converges in $\fM_{1,-1}$, the limit is $V_1^0+V_2^0$ and
\begin{equation*}
\|(V_1^\e+V_2^\e)-(V_1^0+V_2^0)\|_{\fM_{1,-1}}\leqslant
\|V_1^\e-V_1^0\|_{\fM_{1,-1}} + \|V_2^\e-V_2^0\|_{\fM_{1,-1}}.
\end{equation*}
Similarly, if two families $Q_{j,1}^\e$ and $Q_{j,2}^\e$ converge to $Q_{j,1}^0$ and $Q_{j,2}^0$ in $\fM_{1,0}$, the sum $Q_{j,1}^\e+Q_{j,2}^\e$ converges in $\fM_{1,0}$ as well and
\begin{equation*}
\|(Q_{j,1}^\e+Q_{j,2}^\e)-(Q_{j,1}^0+Q_{j,2}^0)\|_{\fM_{1,0}} \leqslant
\|Q_{j,1}^\e-Q_{j,1}^0\|_{\fM_{1,0}} + \|Q_{j,1}^\e-Q_{j,1}^0\|_{\fM_{1,0}}.
\end{equation*}

The next option is two replace one family by another, which is in certain sense equivalent. Namely, given a fixed lattice $\G$ with a periodicity cell $\square$, we introduce the sets $\square_\g^\eta$ by formula (\ref{2.18}).
Let $V^\e\in L_\infty(\Om;\mathds{M}_n)$ be an uniformly bounded family converging to $V^0\in L_\infty(\Om;\mathds{M}_n)$ in $\fM_{1,-1}$. Suppose we are given another uniformly bounded family $\tilde{V}^\e\in L_\infty(\Om;\mathds{M}_n)$ such that
\begin{equation*}%\%label{5.6}
\bigg|\frac{1}{\eta^d\mes\square} \int\limits_{\square_\g^\eta} \big(V^\e(x)-\tilde{V}^\e(x)\big)\,dx
\bigg|\leqslant \rho(\e)\quad \text{for each} \quad \g\in\G_\eta,
\end{equation*}
where $\rho=\rho(\e)$ and $\eta=\eta(\e)$ are some functions independent of $\g$ and $\rho(\e)\to0$, $\eta(\e)\to0$ as $\e\to0$. Then by Theorem~\ref{th2} the family $V^\e-\tilde{V}^\e$ converges to zero in $\fM_{1,-1}$ as $\e\to0$ and
\begin{equation*}%\l%abel{5.7}
\|V^\e-\tilde{V}^\e\|_{\fM_{1,-1}}\leqslant C\big(\rho(\e)+\eta (\e)\big),
\end{equation*}
where $C$ is some constant independent of $\e$, $\eta$ and $\rho$. Then the family $\tilde{V}^\e$ also converges to $V^0$ in $\fM_{1,-1}$ and the estimate holds:
\begin{equation*}%\%label{5.8}
\|\tilde{V}^\e-V^0\|_{\fM_{1,-1}}\leqslant \|V^\e-V^0\|_{\fM_{1,-1}}+ C\big(\rho(\e)+\eta (\e)\big).
\end{equation*}
This property means that on each cell $\square_\g^\eta$ we can replace the family $V^\e$ by another one with the same mean value over $\square_\g^\eta$ and such replacing does not spoil the convergence in $\fM_{1,-1}$.

The latter fact gives rise to an interesting example. If in some cells $\square_\g^\eta$ the family $V^\e$ is non-zero only on some balls contained in these cells, then it is possible to rotate arbitrarily each of these balls changing locally the family $V^\e$. Such rotations however do not change the integrals $\int_{\square_\g^\eta} V^\e \,dx$ and therefore, the rotations do not spoil the convergence in $\fM_{1,-1}$.

A similar property holds also for the convergence in $\fM_{1,0}$. Namely, let $Q^\e\in L_\infty(\Om;\mathds{M}_n)$ be a given uniformly bounded family converging to $Q^0\in L_\infty(\Om;\mathds{M}_n)$ in $\fM_{1,0}$. Let  $\tilde{D}^\e\in L_\infty(\Om;\mathds{M}_n)$ be another uniformly bounded family such that
\begin{equation*}%\%label{5.9}
\int\limits_{\square_\g^\eta} \eta^{-d}(\e)\big|Q^\e(x)- \tilde{D}^\e(x)\big|^2\,dx \leqslant \rho(\e) \quad\text{for each}\quad \square_\g^\eta\subset\Om,\quad
\g\in\G_\eta,
\end{equation*}
where $\rho=\rho(\e)$ and $\eta=\eta(\e)$ are some functions independent of $\g$ and $\rho(\e)\to0$, $\eta(\e)\to0$ as $\e\to0$. Then the family $\tilde{D}^\e$ also converges to $Q^0$ in $\fM_{1,0}$ and
\begin{equation*}%\%label{5.10}
\|\tilde{Q}^\e-Q^0\|_{\fM_{1,0}}\leqslant \|Q^\e-Q^0\|_{\fM_{1,0}}+ C\big(\rho^\frac{1}{2}(\e) +\eta^\frac{1}{2}(\e)\big),
\end{equation*}
where $C$ is a constant independent of $\e$, $\eta$ and $\rho$.

A final possibility of generating new perturbations is related with gluing domains. Namely, suppose we are given two domains $\Om_i$, $i=1,2$.  We suppose that $\Om_1\cap\Om_2=\emptyset$ but $\p\Om_1\cap\p\Om_2\ne\emptyset$ and hence, $\Om_1$ and $\Om_2$ are a partition of some domain into two disjoint parts; this domain is denoted by $\Om$ and it is the interior of $\overline{\Om_1\cup\Om_2}$. We assume that on the domains $\Om_i$, $i=1,2$,  some perturbations
$V^\e_i$, $Q_{j,i}^\e$, $i=1,2$, are defined and  they converge  to $V^0_i$ and $Q_{j,i}^0$ respectively  in $\fM_{1,-1}$ and  $\fM_{1,0}$. Then we define
\begin{equation*}
V^\e:=V^\e_i,\qquad Q_j^\e:=Q_{j,i}^\e,\qquad V^0:=V^0_i,\qquad Q_j^0:=Q_{j,i}^0\qquad
\text{on}\quad \Om_i.
\end{equation*}
The criterions of convergences in $\fM_{1,-1}$ and  $\fM_{1,0}$ stated in Theorems~\ref{th2} and~\ref{th4} are independent of the boundary conditions chosen in (\ref{2.4}) and influenced then the definition of the space $\fV$. This is why for the Dirichlet condition and $\fV=\Ho_2^1(\Om;\mathds{C}^n)$ the convergences (\ref{2.13}) and (\ref{2.20}) are satisfied and hence, the families $V^\e_i$ and $Q_{j,i}^\e$ then converge also in $\fM_{1,-1}$ and  $\fM_{1,0}$ to the same limits for $\fV=W_2^1(\Om;\mathds{C}^n)$.

The restrictions of each function $u\in W_2^1(\Om_i;\mathds{C}^n)$ and each function $v\in L_2(\Om_i;\mathds{C}^n)$ on $\Om_i$ are the elements of respectively $W_2^1(\Om_i;\mathds{C}^n)$ and  $L_2(\Om_i;\mathds{C}^n)$. Hence, in view of definition (\ref{2.7}), (\ref{2.11}) of the norms in $\fM_{1,-1}$ and $\fM_{1,0}$, we conclude that $V^\e$ and $Q_j^\e$ converge to $V^0$ and $Q_j^0$ in these spaces  and
\begin{align*}
&\|V^\e-V^0\|_{\fM_{1,-1}(\Om)}\leqslant \|V_1^\e-V_1^0\|_{\fM_{1,-1}(\Om_1)} + \|V_2^\e-V_2^0\|_{\fM_{1,-1}(\Om_2)},
\\
&\|Q_j^\e-Q_j^0\|_{\fM_{1,0}(\Om)}\leqslant \|Q_{j,1}^\e-Q_{j,1}^0\|_{\fM_{1,0}(\Om_1)} + \|Q_{j,2}^\e-Q_{j,2}^0\|_{\fM_{1,0}(\Om_2)},
\end{align*}
where $\fM_\flat(\Om)$, $\fM_\flat(\Om_i)$, $\flat\in\{(1,-1),(1,0)\}$, $i=1,2$, are the spaces $\fM_\flat$ associated with the domains $\Om$ and $\Om_i$ and the spaces $\fV$ correspond to the Robin condition in (\ref{2.4}).

The examples of possible families $Q_j^\e$ and $V^\e$ provided in previous subsections and the ways of generating new families discussed in this subsection produce a very wide class of possible matrix functions $V^\e$ and $Q_j^\e$, to which our results can be applied. This class includes in particular many non-periodically fast oscillating functions, for which the homogenization is possible in the sense of Theorem~\ref{th1}.

\section{Approximation of resolvent}

In this section we prove Theorems~\ref{th1},~\ref{th5},~\ref{th6}.

\subsection{Proof of Theorem~\ref{th1}}\label{sec:th1}

Let the operators $\cX^\e$ be defined by (\ref{2.37}) and  converge to some operator $\cX^0$ in the norm $\|\,\cdot\,\|_{\fM}$ and the elements $X^\e$ are bounded uniformly in $(L_\infty(\Om;\mathds{C}^n))^{2d+1}$. The space $L_\infty(\Om)$ is dual to the space $L_1(\Om)$. Hence, each entry in the matrix functions $V^\e$, $P_j^\e$, $Q_j^\e$   regarded as a functional  on $L_1(\Om)$, is bounded uniformly in $\e$. Hence, there exists a subsequence $\e'\to+0$ such that the sequence  $X^{\e'}$ converges weakly  as $\e'\to+0$, namely, there exist $V^0, Q_j^0, P_j^0\in L_\infty(\Om;\mathds{C}^n)$ such that
\begin{equation}\label{3.32}
\int\limits_{\Om}(V^{\e'}-V^0)_{ik}(x)\phi(x)\,dx\to0,\qquad
\int\limits_{\Om}(P_j^{\e'}-P_j^0)_{ik}(x)\phi(x)\,dx\to0,\qquad
\int\limits_{\Om}(Q_j^{\e'}-Q_j^0)_{ik}(x)\phi(x)\,dx\to0,
\end{equation}
as $\e'\to0$
for each function $\phi\in L_1(\Om)$, where $(\,\cdot\,)_{ik}$ denotes the entry of a matrix at the intersection of $i$th row and $k$th column.
Given vector functions $u,v\in \fV$ with entries $u_j$ and $v_j$,  the functions $u_j v_j$, $u_j\frac{\p v_j}{\p x_i}$ and $v_j\frac{\p u_j}{\p x_i}$ belong to $L_1(\Om)$. Hence, convergences (\ref{3.32})
yield
%\begin{equation}\%label{3.33}
\begin{align*}
\sum\limits_{j=1}^{d}&
\left(Q_j^{\e'}\frac{\p u}{\p x_j},v\right)_{L_2(\Om;\mathds{C}^n)}- \sum\limits_{j=1}^{d}
\left(P_j^{\e'} u,\frac{\p v}{\p x_j}\right)_{L_2(\Om;\mathds{C}^n)}  + (V^{\e'} u,v)_{L_2(\Om;\mathds{C}^n)}
\\
& \to \sum\limits_{j=1}^{d}
\left(Q_j^0\frac{\p u}{\p x_j},v\right)_{L_2(\Om;\mathds{C}^n)} - \sum\limits_{j=1}^{d}
\left(P_j^0 u,\frac{\p v}{\p x_j}\right)_{L_2(\Om;\mathds{C}^n)} + (V^0 u,v)_{L_2(\Om;\mathds{C}^n)},\quad\e'\to+0,
\end{align*}
%\end{equation}
for all $u,v\in\fV$. At the same time, it follows from the convergence of $\cX^\e$ to $\cX^0$  that
\begin{equation*}%\%label{3.31}
\sum\limits_{j=1}^{d}
\left(Q_j^{\e'}\frac{\p u}{\p x_j},v\right)_{L_2(\Om;\mathds{C}^n)} - \sum\limits_{j=1}^{d}
\left(P_j^{\e'} u,\frac{\p v}{\p x_j}\right)_{L_2(\Om;\mathds{C}^n)} + (V^{\e'} u,v)_{L_2(\Om;\mathds{C}^n)} \to \la \cX_0 u,v\ra,\quad\e\to+0,
\end{equation*}
for all $u,v\in\fV$. These two convergences imply that
\begin{equation*}
\la\cX_0 u,v\ra= \sum\limits_{j=1}^{d}
\left(Q_j^0\frac{\p u}{\p x_j},v\right)_{L_2(\Om;\mathds{C}^n)} - \sum\limits_{j=1}^{d}
\left(P_j^0 u,\frac{\p v}{\p x_j}\right)_{L_2(\Om;\mathds{C}^n)} + (V^0 u,v)_{L_2(\Om;\mathds{C}^n)},
\end{equation*}
which is the desired statement.

We proceed to  proving the existence of the number $\l_0$ mentioned in the formulation of the theorem. The assumed regularity of the boundary $\p\Om$ ensures that the trace of the functions from $W_2^1(\Om)$ is well-defined and belongs to $L_2(\p\Om)$, see, for instance,
\cite[Ch. V, Sec. 5.21, 5.22]{Adams}. Let $\chi_1=\chi_1(\tau)$ be an infinitely differentiable cut-off function equalling to one as $\tau\leqslant \frac{\tau_0}{3}$ and vanishing as $\tau\geqslant \frac{2\tau_0}{3}$. Then for each $u\in W_2^1(\Om)$ we have
\begin{equation*}
|u|^2=-\int\limits_{\tau}^{\tau_0} \frac{\p\ }{\p\tau} (|u|^2\chi_1) \,d\tau\qquad\text{for}\qquad \tau<\frac{\tau_0}{3}
\end{equation*}
and by Cauchy-Schwarz inequality we immediately get
\begin{equation}\label{3.24}
 |u|^2 \leqslant\int\limits_{\tau}^{\tau_0} (2|u||\nabla u|+C|u|^2)\,d\tau\leqslant \int\limits_{0}^{\tau_0} (\d|\nabla u|^2+C(\d)|u|^2)\,d\tau,
\end{equation}
where $\d>0$ is an arbitrary fixed number and   $C(\d)$ is some constant independent of $u$ and $\tau$. Integrating the obtained estimate over $\p\Om$ and employing the assumed condition on the corresponding Jacobians, we get:
\begin{equation}\label{3.2}
\|u\|_{L_2(\p\Om)}\leqslant \d\|\nabla u\|_{L_2(\Om)}+C(\d)\|u\|_{L_2(\Om)},\qquad u\in W_2^1(\Om),
\end{equation}
where  $\d>0$ is an arbitrary fixed number and  $C(\d)$ is some constant independent of $u$.

Using condition (\ref{2.3}) and the uniform boundedness of the families $V^\e$ and $Q_j^\e$ in $L_\infty(\Om)$ we have obvious estimates
%\begin{equation}\%label{3.6}
\begin{align*}
&\RE\sum\limits_{i,j=1}^{d} \left(A_{ij}\frac{\p u}{\p x_j}, \frac{\p u}{\p x_i}\right)_{L_2(\Om)}\geqslant c_1  \|\nabla u\|_{L_2(\Om;\mathds{C}^n)}^2,
\\
& \left|\sum\limits_{j=1}^{d}
\left(A_j\frac{\p u}{\p x_j}, u \right)_{L_2(\Om;\mathds{C}^n)} + (A_0 u, u)_{L_2(\Om;\mathds{C}^n)}\right|\leqslant \d \|\nabla u\|_{L_2(\Om;\mathds{C}^n)}^2 + C(\d) \|u\|_{L_2(\Om;\mathds{C}^n)}^2,
\\
& \left|\sum\limits_{j=1}^{d}
\left(Q_j^\e\frac{\p u}{\p x_j}, u \right)_{L_2(\Om;\mathds{C}^n)} + (V^\e u, u)_{L_2(\Om;\mathds{C}^n)}\right|\leqslant \d \|\nabla u\|_{L_2(\Om;\mathds{C}^n)}^2 + C(\d) \|u\|_{L_2(\Om;\mathds{C}^n)}^2,
\\
& \left|\sum\limits_{j=1}^{d}
\left(Q_j^0\frac{\p u}{\p x_j}, u\right)_{L_2(\Om;\mathds{C}^n)} + (V^0 u, u)_{L_2(\Om;\mathds{C}^n)}\right|\leqslant \d \|\nabla u\|_{L_2(\Om;\mathds{C}^n)}^2 + C(\d) \|u\|_{L_2(\Om;\mathds{C}^n)}^2,
\\
& \Bigg|\sum\limits_{j=1}^{d}
\left(Q_j^\e \frac{\p u}{\p x_j},u\right)_{L_2(\Om;\mathds{C}^n)} - \sum\limits_{j=1}^{d}
\left(P_j^\e u,\frac{\p u}{\p x_j}\right)_{L_2(\Om;\mathds{C}^n)}
\\
&\hphantom{\Bigg|\sum\limits_{j=1}^{d}
\left(Q_j^\e \frac{\p u}{\p x_j},u\right)_{L_2(\Om;\mathds{C}^n)} } + (V^\e u, u)_{L_2(\Om;\mathds{C}^n)}\Bigg|\leqslant \d \|\nabla u\|_{L_2(\Om;\mathds{C}^n)}^2 + C(\d) \|u\|_{L_2(\Om;\mathds{C}^n)}^2,
\\
& \Bigg|\sum\limits_{j=1}^{d}
\left(Q_j^0\frac{\p u}{\p x_j}, u\right)_{L_2(\Om;\mathds{C}^n)}
- \sum\limits_{j=1}^{d}
\left(P_j^0 u,\frac{\p u}{\p x_j}\right)_{L_2(\Om;\mathds{C}^n)}
\\
&\hphantom{\Bigg|\sum\limits_{j=1}^{d}
\left(Q_j^0\frac{\p u}{\p x_j}, u\right)_{L_2(\Om;\mathds{C}^n)}
}+ (V^0 u, u)_{L_2(\Om;\mathds{C}^n)}\Bigg|\leqslant \d \|\nabla u\|_{L_2(\Om;\mathds{C}^n)}^2 + C(\d) \|u\|_{L_2(\Om;\mathds{C}^n)}^2,
\end{align*}
%\end{equation}
where $\d$ is arbitrary and fixed, while $C(\d)$ is some constant independent of $\e$ and $u\in W_2^1(\Om;\mathds{C}^n)$.

Using these estimates  and (\ref{3.2}) with a sufficiently small fixed $\d$,
we obtain:
\begin{align*}
& \RE \hf^\e(u,u)\geqslant \frac{c_1}{2}\|\nabla u\|_{L_2(\Om;\mathds{C}^n)}^2 - C\|u\|_{L_2(\Om;\mathds{C}^n)}^2,
 && \RE \hf^0(u,u)\geqslant \frac{c_1}{2}\|\nabla u\|_{L_2(\Om;\mathds{C}^n)}^2 - C\|u\|_{L_2(\Om;\mathds{C}^n)}^2,
\\
&|\IM \hf^\e(u,u)|\leqslant C\RE \hf^\e(u,u) + C\|u\|_{L_2(\Om;\mathds{C}^n)}^2,
&&|\IM \hf^0(u,u)|\leqslant C \RE \hf^0(u,u) + C\|u\|_{L_2(\Om;\mathds{C}^n)}^2,
\end{align*}
where $C$ are some inessential fixed constants independent of $\e$ and $u\in \fV$ but depending on $\l$ and $\hf^\e$, $\hf^0$ are sesqulinear forms defined in (\ref{2.35}), (\ref{2.36}).  These estimates mean that the numerical ranges of both forms $\hf^0$ and $\hf^\e$ are located in a fixed cone in the complex plane, which lies in the complex half-plane $\RE\l\geqslant \l_0$ with some fixed $\l_0$ independent of $\e$. Hence, all $\l$ outside this half-plane are in the resolvent sets of the operators $\Op^\e$ and $\Op^0$ considered as unbounded ones in $L_2(\Om;\mathds{C}^n)$. This is why the resolvents of both these operators are well-defined for $\RE\l<\l_0$. Moreover, it is clear that there exists $\l_0$ independent of $\e$ such that for $\l$ with $\RE \l<\l_0$ the estimates
\begin{align}\label{3.7}
&\RE \hg^\e_\l(u,u)
\geqslant c_4(\l)\|u\|_{\fV}^2,
%W_2^1(\Om;\mathds{C}^n)}^2,
&& \RE  \hg^0_\l(u,u)
\geqslant c_4(\l)\|u\|_{\fV}^2,
%W_2^1(\Om;\mathds{C}^n)}^2,
\\
&\hg^\e_\l(u,v):=\hf^\e(u,v) -\l(u,v)_{L_2(\Om;\mathds{C}^n)}, &&
\hg^0_\l(u,v):=\hf^0(u,v) -\l(u,v)_{L_2(\Om;\mathds{C}^n)},
\nonumber
\end{align}
hold for all $u\in \fV$, where $c_4=c_4(\l)$ is some fixed positive constant  independent of $u$ and $\e$. In what follows we assume that $\l$ is chosen exactly in such way.

Now we consider the operators $\Op^\e$ and $\Op^0$ as defined on $\fV$ and we are going to show that then for the chosen values of $\l$ the resolvents $(\Op^\e-\l)^{-1}$ and $(\Op^0-\l)^{-1}$ are well-defined and bounded.   We first of all observe that both operators $(\Op^\e-\l)$ and $(\Op^0-\l)$ are bounded as acting from $\fV$ into $\fV^*$ and therefore owing to estimates (\ref{3.7}) the kernels of both operators $(\Op^\e-\l)$ and $(\Op^0-\l)$ are trivial and their images are closed. Hence, in view of the mentioned boundedness of these operators, it is sufficient to prove that their images coincide with $\fV^*$; then the standard Banach theorem on the inverse operator implies the existence and boundedness of the needed resolvents. If the image of the operator $(\Op^\e-\l)$ or $(\Op^0-\l)$ does not coincide with $\fV^*$, due the reflexivity of the space $\fV$, there exists a non-zero element $u\in\fV$ such that $\la f,u\ra=0$ for all $f$ from this image. Letting then $f=(\Op^\e-\l)u$ or $f=(\Op^0-\l)u$, we get  $0=\la f,u\ra=\hg^\e_\l(u,u)$ or $0=\la f,u\ra=\hg^0_\l(u,u)$ and this contradicts estimates (\ref{3.7}) and the inequality $u\ne 0$.

We proceed to proving representation (\ref{2.8}) and estimate (\ref{2.9}). Both operators $(\Op^0-\l)^{-1}$ and $(\Op^\e-\l)^{-1}$ are well-defined and it follows from inequalities (\ref{3.7}) that
\begin{equation}\label{3.36}
\|(\Op^0-\l)^{-1}\|_{\fV^\ast\to \fV} \leqslant c_4^{-1}(\l),\qquad \|(\Op^\e-\l)^{-1}\|_{\fV^\ast\to \fV} \leqslant c_4^{-1}(\l).
\end{equation}
We represent the operator $\Op^\e$ as
$\Op^\e=\Op^0+\cL^\e$ and  hence
\begin{equation}\label{3.34}
\Op^\e-\l=\big(\cI+\cL^\e (\Op^0-\l)^{-1}\big) (\Op^0-\l),
\end{equation}
where $\cI$ is the identity mapping in $\fV^\ast$. The norm of the operator $\cL^\e:\, \fV\to \fV^\ast$ is small and the operator $(\Op^0-\l)^{-1}:\, \fV^\ast\to \fV$ is bounded. Therefore, the operator
\begin{equation*}
\big(\cI+\cL^\e (\Op^0-\l)^{-1}\big)^{-1}:\ \fV^\ast \to \fV^\ast
\end{equation*}
is well-defined for sufficiently small $\e$ and can be represented by the Neumann series
\begin{equation*}
\big(\cI+\cL^\e (\Op^0-\l)^{-1}\big)^{-1}=\sum\limits_{j=0}^{\infty}
\big(-\cL^\e (\Op^0-\l)^{-1}\big)^j
\end{equation*}
uniformly converging in the sense of operator norm $\|\,\cdot\,\|_{\fV^\ast \to \fV^\ast}$. Then it follows from (\ref{3.34}) that
\begin{equation*}%\%label{3.35}
(\Op^\e-\l)^{-1}=(\Op^0-\l)^{-1}\big(\cI+\cL^\e (\Op^0-\l)^{-1}\big)^{-1}=(\Op^0-\l)^{-1}\sum\limits_{j=0}^{\infty}
\big(-\cL^\e (\Op^0-\l)^{-1}\big)^j
\end{equation*}
and this proves representation (\ref{2.8}). Estimate (\ref{2.9}) can be confirmed by direct computations. Indeed, it follows from representation (\ref{2.8}) that
\begin{equation*}
(\Op^\e-\l)^{-1}=(\Op^0-\l)^{-1}\sum\limits_{j=0}^N
\big(-\cL^\e (\Op^0-\l)^{-1}\big)^j + \big(-(\Op^0-\l)^{-1}\cL^\e\big)^{N+1}(\Op^\e-\l)^{-1}.
\end{equation*}
Estimating then the second term in the right hand side of the obtained identity  by means of (\ref{3.36}), we immediately obtain (\ref{2.9}).

\subsection{Proof of Theorem~\ref{th5}}

We choose an arbitrary $u_0\in \fV$ %$f\in \fV^\ast$
and define
\begin{equation} \label{3.9a}
f:=(\Op^0-\l)u_0,\qquad %:=(\Op^0-\l)^{-1}f
u_\e:=(\Op^\e-\l)^{-1}f,\qquad %,\qquad
v_\e:=u_\e-u_0.
\end{equation}
The function $v_\e$ solves the equation
\begin{equation}\label{3.9}
(\Op^\e-\l)v_\e=f_\e,\qquad f_\e:=-\cL^\e(\Op^0-\l)^{-1}f=-\cL^\e u_0,
\end{equation}
obeys the estimate
\begin{equation}\label{3.29}
\|v_\e\|_{\fV}\leqslant \vk(\e) \|f\|_{\fV^*}
\end{equation}
due to (\ref{2.22}) and
satisfies the integral identity
\begin{equation}\label{3.26}
\hg^\e_\l(v_\e,v)
= - \la\cL^\e u_0,v\ra\quad
\text{for each} \quad v\in \fV.
\end{equation}
The function $u_0$ satisfies a similar integral identity
\begin{equation}\label{3.1}
\hg^0_\l(u_0,v)
=\la f,v\ra\quad
\text{for each} \quad v\in \fV.
\end{equation}
It follows from the definition of the forms $\hg^0_\l$ and $\hg^\e_\l$ that
\begin{equation}\label{3.3}
|\hg^0_\l(u_0,v)|\leqslant C \|u_0\|_{\fV} \|v\|_{\fV}, \qquad
|\hg^\e_\l(v_\e,v)|\leqslant C \|v_\e\|_{\fV} \|v\|_{\fV},
\end{equation}
where $C=C(\l)$ is some constant independent of $\e$, $f$, $u_0$, $v_\e$ and $v$.
The above estimate for $|\hg^0_\l(u_0,v)|$ and integral identity (\ref{3.1}) yield the following bound for the norm of $\|f\|_{\fV^*}$:
\begin{equation*}
\|f\|_{\fV^*} \leqslant C\|u_0\|_{\fV}
\end{equation*}
with the same constant $C$ as in (\ref{3.3}).
Hence, by (\ref{3.29}), (\ref{3.26}) and the second estimate in (\ref{3.3}),
\begin{equation}\label{3.28}
\sup\limits_{\substack{u_0\in \fV
\\
u_0\ne0}}
\frac{|\la\cL^\e u_0,v\ra|} {\|u_0\|_{\fV}\|v\|_{\fV}}
\leqslant
C
\sup\limits_{\substack{f\in \fV^\ast
\\
f\ne0
}}
 \frac{\|v_\e\|_{\fV}} {\|u_0\|_{\fV}} \leqslant
C^2
\sup\limits_{\substack{f\in \fV^\ast
\\
f\ne0
}}
 \frac{\|v_\e\|_{\fV}} {\|f\|_{\fV^*}}
 \leqslant C^2 \vk(\e)
 \end{equation}
and this proves (\ref{2.23}). If $P_j^\e=Q_j^\e\equiv0$, $j=1,\ldots,d$, then $P_j^0=Q_j^0\equiv 0$, $j=1,\ldots,d$ and the left hand side in (\ref{2.23}) becomes exactly definition (\ref{2.7}) of the norm in the space of multipliers $\fM_{1,-1}$  and we arrive at (\ref{2.24}).

\subsection{Proof of Theorem~\ref{th6}}

%We follow the notations in the proof of Theorem~\ref{th5}. We assume that
Let $f$ be an arbitrary element of $L_2(\Om;\mathds{C}^n)$ and  $u_0:=(\Op^0-\l)^{-1} f$. Then we define $u_\e$ and $v_\e$ as in (\ref{3.9a}).  Then identity (\ref{3.26}) with $v=v_\e$ becomes
\begin{equation}\label{3.23}
\begin{aligned}
\hg^\e_\l(v_\e,v_\e)
=&
-\la\cL^\e u_0,v_\e\ra
\\
=&\sum\limits_{j=1}^{N}\left((P_j^\e -P_j^0) u_0, \frac{\p v_\e}{\p x_j}\right)_{L_2(\Om;\mathds{C}^n)}-
 \sum\limits_{j=1}^{N}\left((Q_j^\e -Q_j^0) \frac{\p u_0}{\p x_j}, v_\e\right)_{L_2(\Om;\mathds{C}^n)}
\\
&- \big((V^\e-V^0)u_0, v_\e\big)_{L_2(\Om;\mathds{C}^n)}.
\end{aligned}
\end{equation}
Since in the considered case we assume that the resolvent of the operator $\Op^0$ is bounded as an operator from $L_2(\Om;\mathds{C}^n)$ into $W_2^2(\Om;\mathds{C}^n)$, this means the validity of the estimate
\begin{equation}\label{3.25}
\|u_0\|_{W_2^2(\Om;\mathds{C}^n)} \leqslant C\|f\|_{L_2(\Om;\mathds{C}^n)}
\end{equation}
with some constant independent of $f$ but depending on $\l$. This estimate and definition of the norm in $\fM_{1,-1}$
allow us to estimate the right hand side in (\ref{3.23}):
\begin{align*}
\bigg|& \sum\limits_{j=1}^{N}\left((P_j^\e -P_j^0) u_0, \frac{\p v_\e}{\p x_j}\right)_{L_2(\Om;\mathds{C}^n)}-
\sum\limits_{j=1}^{N}\left((Q_j^\e -Q_j^0) \frac{\p u_0}{\p x_j}, v_\e\right)_{L_2(\Om;\mathds{C}^n)}
+\big((V^\e-V^0)u_0, v_\e\big)_{L_2(\Om;\mathds{C}^n)}
\bigg|
\\
&\leqslant \bigg(\sum\limits_{j=1}^{N} \|P_j^\e-P_j^0\|_{\fM_{2,0}}+
\sum\limits_{j=1}^{N} \|Q_j^\e-Q_j^0\|_{\fM_{1,-1}}+  \|V^\e-V^0\|_{\fM_{1,-1}} \bigg)\|u_0\|_{W_2^2(\Om;\mathds{C}^n)}\|v_\e\|_{W_2^1(\Om;\mathds{C}^n)}.
\end{align*}
Taking then the real part of identity (\ref{3.23}) and using (\ref{3.25}) and the first inequality in (\ref{3.7}), we obtain:
\begin{align*}
\|v_\e\|_{W_2^1(\Om)}\leqslant & c_4^{-1}(\l)\bigg(
\sum\limits_{j=1}^{N} \|P_j^\e-P_j^0\|_{\fM_{2,0}}+
\sum\limits_{j=1}^{N} \|Q_j^\e-Q_j^0\|_{\fM_{1,-1}}+  \|V^\e-V^0\|_{\fM_{1,-1}} \bigg)\|u_0\|_{W_2^2(\Om)}
\\
\leqslant & C(\l) \bigg(
\sum\limits_{j=1}^{N} \|P_j^\e-P_j^0\|_{\fM_{2,0}}+
\sum\limits_{j=1}^{N} \|Q_j^\e-Q_j^0\|_{\fM_{1,-1}}+  \|V^\e-V^0\|_{\fM_{1,-1}} \bigg)\|f\|_{L_2(\Om)},
\end{align*}
where $C(\l)$ is some constant independent of $\e$ and $f$. This proves estimate (\ref{2.10}) and completes the proof of Theorem~\ref{th6}.

\section{Convergence in spaces of multipliers}

In this section we prove criterions ensuring the convergence in the spaces $\fM_{1,-1}$ and $\fM_{1,0}$.

\subsection{Proof of Theorem~\ref{th2}}

Given $u\in W_2^1(\square_\g^\eta;\mathds{C}^n)$, we define
\begin{equation}\label{5.71}
\begin{aligned}
&\la u\ra_\g:=\frac{1}{\eta^d\mes \square} \int\limits_{\square_\g^\eta}
u(x)\,dx,\qquad u_\g:=u-\la u\ra_\g,\qquad \int\limits_{\square_\g^\eta} u_\g(x)\,dx=0,
\\
& \|u\|_{L_2(\square_\g^\eta)}^2=\eta^d
|\la u\ra_\g|_{\mathds{C}^n}^2\mes\square + \|u_\g\|_{L_2(\square_\g^\eta;\mathds{C}^n)}^2.
\end{aligned}
\end{equation}
Since all sets $\square_\g^\eta$ are obtained by shifts of $\eta\square$, the Poincar\'e inequality holds:
\begin{equation}\label{4.1}
\|u_\g\|_{L_2(\square_\g^\eta;\mathds{C}^n)}\leqslant C\eta
\|\nabla u\|_{L_2(\square_\g^\eta;\mathds{C}^n)},
\end{equation}
where $C$ is some fixed constant independent of $\g$, $\eta$ and $u$.
We denote
\begin{equation*}%\%label{4.2}
Z_\g^\e:=\frac{1}{\eta^d\mes \square}\int\limits_{\square_\g^\eta} (V^\e(x)-V^0(x))\,dx,\qquad \g\in\G_\eta.
\end{equation*}
By (\ref{2.13}) we have:
\begin{equation}\label{4.3}
|Z_\g^\e|\leqslant \frac{\rho_1(\e)}{\mes\square}.
\end{equation}

Let $u,v\in W_2^1(\Om;\mathds{C}^n)$ be arbitrary functions. Then on each $\square_\g^\eta$ we represent them as $u=\la u\ra_\g + u_\g$, $v=\la v\ra_\g + v_\g$ and we obtain:
\begin{equation}\label{5.70}
\begin{aligned}
\big((V^\e-V^0)u,v\big)_{L_2(\square_\g^\eta;\mathds{C}^n)} =&
(Z_\g^\e\la u\ra_\g,\la v\ra_\g)_{\mathds{C}^n}\eta^d\mes\square
+ \big((V^\e-V^0) \la u\ra_\g,v_\g\big)_{L_2(\square_\g^\eta;\mathds{C}^n)}
\\
&+ \big((V^\e-V^0) u_\g,v\big)_{L_2(\square_\g^\eta;\mathds{C}^n)}.
\end{aligned}
\end{equation}
Using the uniform in $\e$ boundedness of $V^\e-V^0$ in $L_\infty(\Om;\mathds{M}_n)$, we estimate the terms in the right hand side of the above identity by means of (\ref{4.3}), (\ref{4.1}):
\begin{equation}\label{5.72}
\begin{aligned}
&\big|(V_\g\la u\ra_\g,\la v\ra_\g)_{\mathds{C}^n}\eta^d\mes\square
\big|\leqslant C\rho_1 \|u\|_{L_2(\square_\g^\eta;\mathds{C}^n)} \|v\|_{L_2(\square_\g^\eta;\mathds{C}^n)},
\\
&\big| \big((V^\e-V^0) \la u\ra_\g,v_\g\big)_{L_2(\square_\g^\eta;\mathds{C}^n)}
\big|\leqslant C\eta^{\frac{d}{2}}|\la u\ra_\g|\|v_\g\|_{L_2(\square_\g^\eta;\mathds{C}^n)}
\leqslant C\eta \|u\|_{L_2(\square_\g^\eta;\mathds{C}^n)} \|\nabla v\|_{L_2(\square_\g^\eta;\mathds{C}^n)},
\\
&\big| \big((V^\e-V^0) u_\g,v\big)_{L_2(\square_\g^\eta;\mathds{C}^n)}
\big| \leqslant C \|u_\g\|_{L_2(\square_\g^\eta;\mathds{C}^n)} \|v\|_{L_2(\square_\g^\eta;\mathds{C}^n)} \leqslant C \eta \|\nabla u\|_{L_2(\square_\g^\eta;\mathds{C}^n)} \|v\|_{L_2(\square_\g^\eta;\mathds{C}^n)},
\end{aligned}
\end{equation}
where $C$ is some constant independent of $\eta$, $\g$, $u$ and $v$. Hence,
\begin{equation}\label{4.4}
\big|\big((V^\e-V^0)u,v\big)_{L_2(\square_\g^\eta;\mathds{C}^n)} \big| \leqslant C(\rho_1+\eta)\| u\|_{W_2^1(\square_\g^\eta;\mathds{C}^n)} \|v\|_{W_2^1(\square_\g^\eta;\mathds{C}^n)},
\end{equation}
where $C$ is some constant independent of $\eta$, $\g$, $u$ and $v$.

 We recall definition (\ref{2.0}) of the set $\Pi(\,\cdot\,)$ and
we integrate estimate (\ref{3.24}) over $\Pi(\eta)$ to we see that
\begin{equation}\label{4.9}
\|u\|_{L_2(\Pi(\eta);\mathds{C}^n)}\leqslant C \eta^{\frac{1}{2}} \|u\|_{W_2^1(\Om;\mathds{C}^n)}\quad\text{for each}\quad
u\in W_2^1(\Om;\mathds{C}^n),
\end{equation}
with a constant $C$ independent of $u$ and $\eta$. It follows from the definition of the set $\G_\eta$ that
\begin{equation}\label{4.11}
\Om\setminus\Om_\eta\subset \Pi(C\eta),\qquad \Om_\eta:=\bigcup\limits_{\g\in\G_\eta} \square_\g^\eta,
\end{equation}
where $C$ is some constant independent of $\eta$ and $\e$.
Then estimates (\ref{4.9}) and (\ref{4.4}) yield
\begin{equation*}
\big|\big((V^\e-V^0)u,v\big)_{L_2(\Om;\mathds{C}^n)} \big| \leqslant C(\rho_1+\eta)\| u\|_{W_2^1(\Om;\mathds{C}^n)} \|v\|_{W_2^1(\Om;\mathds{C}^n)}
\end{equation*}
and this proves (\ref{2.14}) as well  as the convergence of $V^\e$ to $V^0$ in $\fM_{1,-1}$.

Assume now that $V^\e$ converges to some $V^0\in L_\infty(\Om;\mathds{M}_n)$ in $\fM_{1,-1}$. We denote $\eta:=\|V^\e-V^0\|_{\fM_{1,-1}}^\frac{1}{4}$.
Then we choose a fixed infinitely differentiable vector function $U_j=U_j(\xi,\eta)$ defined on $\overline{\square}$ such that
$U_j(\xi,\eta)$ vanishes on $\p\square$ and
$U_j(\xi,\eta)\equiv E_j$ as $\xi\in\square$, $\dist(\xi,\p\square)\geqslant \eta$, where $E_j$ is $j$th vector in the standard basis in $\mathds{C}^n$. Without loss of generality we also suppose that
\begin{equation}\label{5.1}
|\nabla_\xi U_j(\xi,\eta)|\leqslant C\eta^{-1}\quad\text{on}\quad \overline{\square},
%, \qquad \int\limits_{\square} |\nabla_\xi U(\xi)|^2\,d\xi \geqslant C\eta^{-1},
\end{equation}
where $C$ is some fixed positive constants independent of $\eta$.
%
% We also suppose that $|\nabla_\xi U_j(\xi,\eta)|\leqslant C\eta^{-1}$.
Then we choose an arbitrary $\g\in\G_\eta$ and we define a function on $\Om$:
\begin{equation*}
u_j(x):=U_j\big(x\eta^{-1}-\g,\eta\big) \quad\text{on}\quad \square_\g^\eta, \qquad u_j(x)=0\quad\text{outside}\quad \square_\g^\eta.
\end{equation*}
It is straightforward to see that
\begin{equation}\label{4.12}
\big((V^\e-V^0)u_i,u_j\big)_{L_2(\Om;\mathds{C}^n)}=
\int\limits_{\square_\g^\eta} \big((V^\e(x)-V^0(x))E_i, E_j\big)_{\mathds{C}^n}\,dx +O(\eta^{d+1}),
\end{equation}
where the $O$-term is uniform in $\e$ and $\g$, namely, it is estimated by $C\eta^{d+1}$, where $C$ is some constant independent of $\e$, $\g$,  $\eta$. Similarly, the $L_2(\Om;\mathds{C}^n)$-norm of $u_i$ satisfies the identities
\begin{equation}\label{4.13}
\|u_i\|_{L_2(\Om;\mathds{C}^n)}^2 =\|u_i\|_{L_2(\square_\g^\eta;\mathds{C}^n)}^2
=\mes  \square_\g^\eta  + O(\eta^{d+1})=
\eta^d\mes \square  + O(\eta^{d+1}),
\end{equation}
where the $O$-term is again uniform in $\g$ and $\e$. The main contribution to the  norm of the gradient of $u_i$ is made by a thing layer along the boundary $\p\square_\g^\eta$, see (\ref{5.1}):
\begin{equation*}%\%label{4.14}
\|\nabla u_i\|_{L_2(\Om;\mathds{C}^n)}^2= \|\nabla u_i\|_{L_2(\square_\g^\eta;\mathds{C}^n)}^2\leqslant C\eta^{d-3},
\end{equation*}
where $C$ is some positive constant independent of $\e$, $\eta$ and $\g$.
This estimate and (\ref{4.12}), (\ref{4.13}) then imply:
%\begin{equation}\%label{4.15}
\begin{align*}
\frac{\big|\big((V^\e-V^0)u_i, u_j\big)_{L_2(\Om;\mathds{C}^n)}\big|} {\|u_i\|_{W_2^1(\Om;\mathds{C}^n)} \|u_j\|_{W_2^1(\Om;\mathds{C}^n)}} \geqslant & \frac{\bigg|\int\limits_{\square_\g^\eta} \big((V^\e(x)-V^0(x))E_i, E_j\big)_{\mathds{C}^n}\,dx \bigg|+ O(\eta^{d+1})}{C\eta^{d-3} +O(\eta^{d})}
\\
\geqslant & \frac{\eta^3}{c\eta^d} \Bigg|\int\limits_{\square_\g^\eta} \big((V^\e(x)-V^0(x))E_i, E_j\big)_{\mathds{C}^n}\,dx \Bigg| + O(\eta^4).
\end{align*}
%\end{equation}
Comparing this inequality with definition (\ref{2.7}) of the norm in $\fM_{1,-1}$, we see that
\begin{equation*}
\frac{\eta^3}{c\eta^d} \Bigg|\int\limits_{\square_\g^\eta} \big((V^\e(x)-V^0(x))E_i, E_j\big)_{\mathds{C}^n}\,dx \Bigg| + O(\eta^4)\leqslant \|V^\e-V^0\|_{\fM_{1,-1}}
\end{equation*}
and hence,
\begin{equation*}%\%label{4.16}
\frac{1}{\eta^d} \Bigg|\int\limits_{\square_\g^\eta} \big((V^\e(x)-V^0(x))E_i, E_j\big)_{\mathds{C}^n}\,dx \Bigg|\leqslant
C(\eta^{-3}\|V^\e-V^0\|_{\fM_{1,-1}} + \eta)\leqslant C\|V^\e-V^0\|_{\fM_{1,-1}}^\frac{1}{4},
\end{equation*}
where $C$ is some constant independent of $\e$, $\g$, $\eta$. Summing up the above estimate over $i$ and $j$, we see that condition (\ref{2.13}) holds with the above chosen function $\eta$ and $\rho_1=C\|V^\e-V^0\|_{\fM_{1,-1}}^\frac{1}{4}$. The proof is complete.

\subsection{Proof of Theorem~\ref{th3}}

Since the family $V^\e$ is uniformly bounded in $L_\infty(\Om;\mathds{M}_n)$, it follows from   definition (\ref{2.16}) of the function $V^0$ that
\begin{equation*}
|V^0(x)|\leqslant \|V^\e\|_{L_\infty(\Om;\mathds{M}_n)} \leqslant C
\end{equation*}
for almost each $x\in\Om$, where $C$ is some fixed constant independent of $x\in\Om$ and $\e$. Hence, $V^0\in L_\infty(\Om;\mathds{C}^n)$.

We choose an arbitrary fixed lattice $\G\subset\mathds{R}^d$
and introduce sets $\G_\eta$ and $\square_\g^\eta$  by formula (\ref{2.18}) with $\eta:=\mu^{\frac{1}{2}}$. Then for each $\g\in   \G_\eta$ we have
\begin{equation}\label{4.6}
\begin{aligned}
\frac{1}{\eta^d\mes\square} \int\limits_{\square_\g^\eta} \big(V^\e(x)-V^0(x)\big)\,dx= &\frac{1}{\eta^d\mes\square} \int\limits_{\square_\g^\eta} \bigg(V^\e(x) - \frac{1}{\mu^d\mes\om}
\int\limits_{\mu\om+x} V^\e (y)\,dy\bigg)\,dx
\\
&+
\frac{1}{\eta^d\mes\square} \int\limits_{\square_\g^\eta} \bigg( \frac{1}{\mu^d\mes\om}
\int\limits_{\mu\om+x} V^\e (y)\,dy-V^0(x)\bigg)\,dx.
\end{aligned}
\end{equation}
By condition (\ref{2.17}) we then immediately get:
\begin{equation}\label{4.7}
\begin{aligned}
\bigg|\frac{1}{\eta^d\mes\square} \int\limits_{\square_\g^\eta} \frac{dx}{\mu^d\mes\om}\bigg(
\int\limits_{\mu\om+x} V^\e (y)\,dy-V^0(x)\bigg)
\bigg|\leqslant \rho_2.
\end{aligned}
\end{equation}

We rewrite the first term in the right hand side of (\ref{4.6}) as
\begin{equation}\label{4.8}
\begin{aligned}
\frac{1}{\eta^d\mes\square}\int\limits_{\square_\g^\eta} \bigg(V^\e(x) - \frac{1}{\mu^d\mes\om}
\int\limits_{\mu\om+x} V^\e (y)\,dy\bigg)\,dx = & \frac{1}{\eta^d\mu^d\mes\square\mes\om} \int\limits_{\square_\g^\eta}dx \int\limits_{\mu\om} \big(V^\e(x) -   V^\e (x+y)\big) dy
\\
=&\frac{1}{\eta^d\mu^d\mes\square\mes\om}\int\limits_{\mu\om} dy \bigg(\int\limits_{\square_\g^\eta}V^\e(x) dx - \int\limits_{\square_\g^\eta+y}V^\e(x) dx
 \bigg).
\end{aligned}
\end{equation}
Since the family $V^\e$ is uniformly bounded in $L_\infty(\Om;\mathds{M}_n)$, we have the estimate:
\begin{align*}
\bigg|\int\limits_{\square_\g^\eta}V^\e(x) dx - \int\limits_{\square_\g^\eta+y}V^\e(x) dx\bigg|\leqslant & \|V^\e\|_{L_\infty(\Om;\mathds{M}_n)} \big(\mes \square_\g^\eta\setminus (\square_\g^\eta+y)+ \mes (\square_\g^\eta+y)\setminus\square_\g^\eta \big)
\\
\leqslant & C |y| \mes \p\square_\g^\eta \leqslant C \mu\eta^{d-1},
\end{align*}
where $C$ are some constants independent of $\e$, $\eta$, $\mu$.
Substituting this estimate into (\ref{4.8}), we obtain:
\begin{align*}
\bigg|\frac{1}{\eta^d\mes\square}\int\limits_{\square_\g^\eta} \bigg(V^\e(x) - \frac{1}{\mu^d\mes\om}
\int\limits_{\mu\om+x} V^\e (y)\,dy\bigg)\,dx\bigg| \leqslant \frac{C\mu\eta^{d-1}}{\mu^d\eta^d \mes\square\mes\om} \int\limits_{\mu\om} dy \leqslant C\frac{\mu}{\eta}=C\mu^\frac{1}{2},
\end{align*}
where $C$ are some constants independent of $\e$, $\eta$, $\mu$. This estimate and (\ref{4.7}) ensure that the assumptions of Theorem~\ref{th2} are satisfied with $\eta=\mu^{\frac{1}{2}}$ and $\rho_1=\rho_2+\mu^\frac{1}{2}$. Rewriting then estimate (\ref{2.14}) in terms of these functions, we arrive at (\ref{2.19}). The proof is complete.

\subsection{Proof of Theorem~\ref{th4}}\label{sec53}

We introduce the sets $\G_\eta$ and $\Om_\eta$ by formulae (\ref{2.18}), (\ref{4.11}). Given a function $u\in W_2^1(\Om)$, on each $\square_\g^\eta$, $\g\in\G_\eta$, we represent $u$ as $u=\la u\ra_\g+u_\g$ and then we have the estimate
\begin{equation*}
\int\limits_{\square_\g^\eta} |(Q^\e-Q^0)u|^2\,dx\leqslant
 \int\limits_{\square_\g^\eta} |Q^\e-Q^0|^2|u|^2\,dx
\leqslant
 C |\la u\ra_\g|^2\int\limits_{\square_\g^\eta} |Q^\e-Q^0|^2\,dx+C
\int\limits_{\square_\g^\eta} |Q^\e-Q^0|^2|u_\g|^2\,dx,
\end{equation*}
where $C$ is some function independent of $\e$, $\g$, $\eta$, $u$, $Q^\e$ and $Q^0$. Then by the uniform boundedness of $Q^\e$ in $L_\infty(\Om;\mathds{M}_n)$, the belonging $Q^0\in L_\infty(\Om;\mathds{M}_n)$, condition (\ref{2.20}) and estimate (\ref{4.1}) we obtain:
\begin{equation*}
\int\limits_{\square_\g^\eta} |(Q^\e-Q^0)u|^2\,dx\leqslant  C\rho_3\|u\|_{L_2(\square_\g^\eta)}^2 + C\eta^2 \|\nabla u\|_{L_2(\square_\g^\eta)}^2,
\end{equation*}
where $C$ is a constant independent of $\e$, $\g$, $\eta$ and $u$. Summing up the above estimate over $\g\in\G_\eta$, we find:
\begin{equation}\label{4.10}
\|(Q^\e-Q^0)u\|_{L_2(\Om_\eta)}\leqslant C\big(\rho_3^\frac{1}{2}+\eta\big) \|u\|_{W_2^1(\Om_\eta)},
\end{equation}
with a constant $C$ independent of $\e$, $\eta$, $\rho_3$ and $u$.
It follows from estimate (\ref{4.9}) that
\begin{equation*}
\|(Q^\e-Q^0)u\|_{L_2(\Om\setminus\Om_\eta)}\leqslant C\eta^\frac{1}{2} \|u\|_{W_2^1(\Om_\eta)}
\end{equation*}
with a constant $C$ independent of $\e$, $\eta$ and $u$.
This inequality and (\ref{4.10}) imply the first estimate in (\ref{2.21}). The second estimate in (\ref{2.21}) can be proved in the same way.

The proof of the second part of theorem follows the same lines as a similar part in the proof of Theorem~\ref{th2}. We let $\eta:= \|Q^\e-Q^0\|_{\fM_{1,0}}^\frac{1}{2}$.
Then we again choose a   fixed infinitely differentiable vector function $U=U(\xi)$ defined on $\overline{\square}$ such that
$U(\xi)$ vanishes on $\p\square$ and
$U(\xi)\equiv U_0$ as $\xi\in\square$, $\dist(\xi,\p\square)\geqslant \eta$, where $U_0$ is an arbitrary vector in $\mathds{C}^n$ with the unit norm.  Without loss of generality we   suppose that
\begin{equation*}
 |\nabla_\xi U(\xi)|\leqslant C\eta^{-1}\quad\text{on}\quad \overline{\square},
 %, \qquad \int\limits_{\square} |\nabla_\xi U(\xi)|^2\,d\xi\geqslant C\eta^{-1},
\end{equation*}
where $C$ is a fixed positive constant independent of $\eta$.
For an arbitrary $\g\in\G_\eta$ we define a function on $\Om$:
\begin{equation*}
u(x):=U\big(x\eta^{-1}-\g\big) \quad\text{on}\quad \square_\g^\eta, \qquad u(x)=0\quad\text{outside}\quad \square_\g^\eta.
\end{equation*}
We again see that
\begin{equation}\label{4.17}
\|(Q^\e-Q^0)u\|_{L_2(\Om;\mathds{C}^n)}^2=
\int\limits_{\square_\g^\eta} |(Q^\e-Q^0)U_0|_{\mathds{C}^n}^2  \,dx +O(\eta^{d+1}),
\end{equation}
where the $O$-term is uniform in $\e$, $\g$ and $U_0$. We also have similar relations
\begin{align*}
&
\|u\|_{L_2(\Om;\mathds{C}^n)}^2 =\|u\|_{L_2(\square_\g^\eta;\mathds{C}^n)}^2
=\mes \square_\g^\eta  + O(\eta^{d+1})=
\eta^d\mes  \square  + O(\eta^{d+1}),
\\
&\|\nabla u\|_{L_2(\Om;\mathds{C}^n)}^2= \|\nabla u\|_{L_2(\square_\g^\eta;\mathds{C}^n)}^2 \leqslant c\eta^{d-3},
\end{align*}
where the $O$-term is again uniform in $\g$, $\e$, $U_0$, while  $c$ is some positive constant independent of $\e$, $\eta$, $U_0$ and $\g$. By these estimates and (\ref{4.17}) we then get:
\begin{equation*}
\frac{\big\|(Q^\e-Q^0)u \big\|_{L_2(\Om;\mathds{C}^n)}^2 } {\|u\|_{W_2^1(\Om;\mathds{C}^n)}^2} \geqslant   \frac{\int\limits_{\square_\g^\eta}
\big|(Q^\e-Q^0)U_0\big|^2
 \,dx + O(\eta^{d+1})}{c\eta^{d-3} +O(\eta^{d})}
\geqslant   \frac{\eta^3}{c\eta^d} \int\limits_{\square_\g^\eta} \big|(Q^\e-Q^0)U_0\big|^2\,dx  + O(\eta^4).
\end{equation*}
And comparing this inequality with definition (\ref{2.11}) of the norm in $\fM_{1,0}$, we find:
\begin{equation*}
\frac{\eta^3}{c\eta^d} \int\limits_{\square_\g^\eta} \big|(Q^\e-Q^0)U_0\big|^2\,dx  + O(\eta^4)\leqslant \|Q^\e-Q^0\|_{\fM_{1,0}}^2.
\end{equation*}
Therefore,
\begin{equation*}
\frac{1}{\eta^d} \int\limits_{\square_\g^\eta} \big|(Q^\e-Q^0)U_0\big|^2\,dx \leqslant
C(\eta^{-3}\|Q^\e-Q^0\|_{\fM_{1,0}}^2 + \eta),
\end{equation*}
where $C$ is some constant independent of $\e$, $\g$, $\eta$ and $U_0$.  This yields:
\begin{equation*}
\frac{1}{\eta^d} \int\limits_{\square_\g^\eta} \big|(Q^\e-Q^0)\big|^2\,dx \leqslant
C(\eta^{-3}\|Q^\e-Q^0\|_{\fM_{1,0}}^2 + \eta)
\end{equation*}
and we see that condition (\ref{2.20}) holds with the above chosen function $\eta$ and $\rho_3=C \|Q^\e-Q^0\|_{\fM_{1,0}}^\frac{1}{2}$.

Let us prove inequality (\ref{2.28}). If $\fV=W_2^1(\Om;\mathds{C}^n)$, by $\mathfrak{v}$ we denote the Sobolev space $W_2^1(\Om)$ of scalar functions, while in the case $\fV=\Ho^1(\Om;\mathds{C}^n)$ we let $\mathfrak{v}:=\Ho^1(\Om)$. In definition (\ref{2.11}) of the norm in $\fM_{1,0}$ we choose the vector function $u$, all components of which are zero except for the $j$th one. Then we see that the function $\Big(\sum\limits_{i=1}^{n} |Q_{ij}|^2\Big)^\frac{1}{2}$, where $Q_{ij}$ are the entries of the matrix $Q$, is a scalar multiplier from $\mathfrak{v}$ into $L_2(\Om)$ and its norm is estimated by  $\|Q\|_{\fM_{1,0}}$:
\begin{equation*}
\sup\limits_{\substack{u\in\mathfrak{v}\\ u\ne 0}} \frac{\bigg\|\bigg(\sum\limits_{i=1}^{n} |Q_{ij}|^2\bigg)^\frac{1}{2}u\bigg\|_{L_2(\Om)}}{\|u\|_{\mathfrak{v}}}\leqslant \|Q\|_{\fM_{1,0}}.
\end{equation*}
Using this estimate and the Cauchy-Schwartz inequality, for each $u\in\fV$ we easily find
\begin{equation*}
\|Q^*u\|_{L_2(\Om;\mathds{C}^n)}^2\leqslant \sqrt{n} \|Q\|_{\fM_{1,0}} \|u\|_{\fV}
\end{equation*}
and in view of (\ref{2.11}) this proves (\ref{2.28}).
The proof is complete.

\section*{Acknowledgments}

The author thanks A.L.~Piatnitski,  A.A.~Fedotov, E.A.~Zhizhina  for many useful discussions and valuable remarks. The author is deeply grateful to T.A.~Suslina and N.N.~Senik for useful discussions, a careful reading of the paper and many valuable remarks and comments.

\section*{Declarations of interest}

None.

\end{document}